\documentclass[a4paper,10pt,oneside]{article}

\usepackage{t1enc}
\usepackage[latin1]{inputenc}
\usepackage{graphics}
\usepackage{fix-cm}

\usepackage{graphicx}
\usepackage{amsfonts}
\usepackage{eufrak}
\usepackage[all]{xypic}
\usepackage{amssymb}
\usepackage{amsmath}
\usepackage{amsthm}
\usepackage{mathrsfs}
\usepackage{makeidx}

\theoremstyle{plain}
\newtheorem{thm}{Theorem}[section]
\newtheorem{lemma}[thm]{Lemma}
\newtheorem{prop}[thm]{Proposition}
\newtheorem{cor}[thm]{Corollary}
\newtheorem*{convention}{Convention}

\theoremstyle{definition}
\newtheorem{df}[thm]{Definition}
\newtheorem{rem}[thm]{Remark}
\newtheorem{ex}[thm]{Example}

\newtheorem{standing}[thm]{Standing Assumption}

\renewcommand{\a}{\mathcal{A}}
\renewcommand{\b}{\mathcal{B}}

\providecommand{\gm}{\Gamma}
\providecommand{\h}{\mathcal{H}}

\title{\bf Crossed products by Hecke pairs I: $^*$-algebraic theory}
\author{Rui Palma\\ {\footnotesize \emph{Department of Mathematics, University of Oslo,}}\\  {\footnotesize \emph{P.O. Box 1053 Blindern, NO-0316 Oslo, Norway}}\\ {\footnotesize\emph{E-mail: ruip@math.uio.no}}}
\date{}

\makeindex

\begin{document}
\maketitle

\begin{abstract}
We develop a theory of crossed products by ``actions'' of Hecke pairs $(G, \gm)$, motivated by applications in non-abelian $C^*$-duality. Our approach gives back the usual crossed product construction whenever $G / \gm$ is a group and retains many of the aspects of crossed products by groups. In this first of two articles we lay the $^*$-algebraic foundations of these crossed products by Hecke pairs and we explore their representation theory.\\
\end{abstract}

\tableofcontents

{\renewcommand{\thefootnote}{}
\footnotetext{\emph{Date:} \today}}

{\renewcommand{\thefootnote}{}
\footnotetext{Research supported by the Research Council of Norway and the Nordforsk research network ``Operator Algebra and Dynamics''}}

\section*{Introduction}
\addcontentsline{toc}{section}{Introduction}

The present work is the first of two articles whose goal is the development of a theory of crossed products by Hecke pairs with a view towards application in non-abelian $C^*$-duality.

A \emph{Hecke pair} $(G, \gm)$ consists of a group $G$ and a subgroup $\gm \subseteq G$ for which every double coset $\gm g \gm$ is the union of finitely many left cosets. In this case $\gm$ is also said to be a \emph{Hecke subgroup} of $G$. Examples of Hecke subgroups include finite subgroups, finite-index subgroups and normal subgroups. It is in fact many times insightful to think of this definition as a generalization of the notion of normality of a subgroup.

Given a Hecke pair $(G,\gm)$ the \emph{Hecke algebra} $\h(G, \gm)$ is a $^*$-algebra of functions over the set of double cosets $\gm \backslash G / \gm$, with a suitable convolution product and involution. It generalizes the definition of the group algebra $\mathbb{C}(G / \gm)$ of the quotient group when $\gm$ is a normal subgroup.

Heuristically, a crossed product of an algebra $A$ by a Hecke pair $(G,\gm)$ should be thought of as a crossed product (in the usual sense) of $A$ by an ``action'' of $G / \gm$. The quest for a sound definition of crossed products by Hecke pairs may seem hopelessly flawed since $G / \gm$ is not necessarily a group and thus it is unclear how it should ``act'' on $A$. It is the goal of this article and its coming sequel to show that in some circumstances such a definition can be given in a meaningful way, recovering the original one whenever $G/\gm$ is a group.

The term ``crossed product by a Hecke pair'' was first used by Tzanev \cite{tzanev talk} in order to give another perspective on the work of Connes and Marcolli \cite{connes marcolli}. This point of view was later formalized by Laca, Larsen and Neshveyev in \cite{phase}, where they defined a $C^*$-algebra which can be interpreted as a reduced $C^*$-crossed product of a commutative $C^*$-algebra by a Hecke pair.

It seems to be a very difficult task to define crossed products of \emph{any} given algebra $A$ by a Hecke pair, and for this reason we set as our goal to define a crossed product by a Hecke pair in a generality that will cover the following aspects:
\begin{itemize}
\item existence of a canonical spanning set of elements in the crossed product;
 \item possibility of defining covariant representations;
 \item the Hecke algebra must be a trivial example of a crossed product by a Hecke pair;
 \item the classical definition of a crossed product must be recovered whenever $G / \gm$ is a group;
 \item our construction should agree with that of Laca, Larsen and Neshveyev, whenever they are both definable;
 \item our definition should be suitable for applications in non-abelian $C^*$-duality.
\end{itemize}

In this first article on this subject we focus on defining such crossed products on a purely $^*$-algebraic level and on developing their representation theory. The subjects of $C^*$-completions and relations with non-abelian $C^*$-duality will be further explored in the second article on the subject.

     We develop a theory of crossed products of certain algebras $A$ by Hecke pairs which takes into account the above requirements. Our approach makes sense when $A$ is a certain algebra of sections of a Fell bundle over a discrete groupoid. To summarize our set up: we start with a Hecke pair $(G, \gm)$, a Fell bundle $\a$ over a discrete groupoid $X$ and an action $\alpha$ of $G$ on $\a$  satisfying some ``nice'' properties. From this we naturally give the space $\a / \gm$ of $\gm$-orbits of $\a$ a Fell bundle structure over the orbit space $X / \gm$, which under our assumptions on the action $\alpha$ is in fact a groupoid. We can then define a $^*$-algebra
\begin{align*}
 C_c(\a / \gm) \times^{alg}_{\alpha} G / \gm\,,
\end{align*}
which can be thought of as the crossed product of $C_c(\a / \gm)$ by the Hecke pair $(G, \gm)$. We should point out that a \emph{crossed product} for us is simply a $^*$-algebra, which we can then complete with different $C^*$-norms or an $L^1$-norm. Hence, and so that no confusion arises, the symbol $\times^{alg}$ will always be used when talking about the (uncompleted) $^*$-algebraic crossed product.

Our construction gives back the usual crossed product construction when $\gm$ is a normal subgroup of $G$. Moreover, given any action of the group $G / \gm$ on a Fell bundle $\b$ over a groupoid $Y$, the usual crossed product $C_c(\b) \times^{alg} G / \gm$ can be obtained via our setup as a crossed product by the Hecke pair $(G, \gm)$.

 Many of the features present in crossed products by discrete groups carry over to our setting. For instance, the role of the group $G / \gm$ is played by the Hecke algebra $\h(G, \gm)$, which embeds in a natural way in the multiplier algebra of $C_c(\a / \gm) \times^{alg} G / \gm$. Additionally, just like a crossed product $A \times G$ by a discrete group is spanned by elements of the form $a *g$, with $a \in A$ and $g \in G$, our crossed products by Hecke pairs also admit a canonical spanning set of elements.

The representation theory of crossed products by Hecke pairs also has many similarities with the group case, but some distinctive new features arise. For instance, as it is well-known in the group case, there is a bijective correspondence between nondegenerate representations of a crossed product $A \times G$ and the so-called covariant representations of $A$ and $G$, which are certain pairs of unitary representations of $G$ and representations of $A$. We will show that something completely analogous occurs for Hecke pairs, but in this case one is obliged to consider \emph{pre-representations} of the Hecke algebra, i.e. representations of $\h(G, \gm)$ as (possibly) unbounded operators. This consideration was unnecessary in the group case because unitary operators are automatically bounded.

As stated before, this theory of crossed products by Hecke pairs is intended for applications in non-abelian $C^*$-duality theory. One of the main motivations is the establishment of a Stone-von Neumann theorem for Hecke pairs that encompasses the work of an Huef, Kaliszewski and Raeburn \cite{cov} and expresses their results in the language of crossed products. Additionally, we envisage for future work a form of Katayama duality with respect to Echterhoff-Quigg's   ``crossed product'' \cite{full} (a terminology used in \cite{cov}). In a succinct and non-rigorous way this would mean that there is a canonical isomorphism of the type:
\begin{align*}
A \times_{\delta} G / \gm \times_{\widehat{\delta}, \omega} G / \gm \cong A \otimes \mathcal{K}(\ell^2(G / \gm))\,,
\end{align*}
where $A \times_{\delta} G / \gm$ is a crossed product by a coaction of the homogeneous space $G / \gm$, while the second crossed product should be by the ``dual action'' of the Hecke pair $(G, \gm)$ in our sense. Such a result would bring insight into the emerging theory of crossed products by coactions of homogeneous spaces (\cite{echt kal rae}, \cite{full}). We explain in Example \ref{katayama duality example} how our construction adapts well to the settings of \cite{full}.

This article is organized as follows. In Section \ref{preliminaries section} we set up the conventions and prelimanry results to be used in the rest of the article.

Section \ref{groupoids Fell bundles chapter} is dedicated to the development of the required set up for defining crossed products by Hecke pairs. Here we explain what type of actions are involved, how to define the orbit space groupoids $X / H$ and the orbit bundles $\a / H$ out of $\a$, and how all the algebras $C_c(\a / H)$ are related with each other for different subgroups $H \subseteq G$.

Lastly, in Section \ref{algebraic crossed products chapter} we introduce the notion of a crossed product by a Hecke pair, explore some of its algebraic aspects and develop its representation theory. In the last part of this section we show how many of the formulas become much simpler in the case of free actions.

The present work is based on the author's Ph.D. thesis \cite{palmathesis} written at the University of Oslo. There are a few differences between the present work and \cite{palmathesis}, notably the greater generality of the types of actions involved. This improvement follows a suggestion of Dana Williams and John Quigg.

The author would like to thank his advisor Nadia Larsen for the very helpful discussions, suggestions and comments during the elaboration of this work. A word of appreciation goes also to John Quigg, Dana Williams and Erik B\'edos for some very helpful comments.\\

\section{Preliminaries}
\label{preliminaries section}

In this section we set up the conventions, notation, and background results which will be used throughout this work. We indicate the references where the reader can find more details, but we also provide proofs for those results which we could not find in the literature.\\

\begin{convention}
 The following convention for displayed equations will be used throughout this work: if a displayed formula starts with the equality sign, it should be read as a continuation of the previously displayed formula.

A typical example takes the following form:

\begin{eqnarray*}
 \mathrm{(expression\;\; 1)} & = & \mathrm{(expression\;\; 2)}\\
& = & \mathrm{(expression\;\; 3)}\,.
\end{eqnarray*}
\emph{By Theorem A and Lemma B it then follows that}
\begin{eqnarray*}
\qquad\qquad\quad & = & \mathrm{(expression\;\; 4)}\\
& = & \mathrm{(expression\;\; 5)}\,.
\end{eqnarray*}

\vskip0.4cm

Under our convention starting with the equality sign in the second array of equations simply means that $\mathrm{(expression\;\; 3)}$ is equal to $\mathrm{(expression\;\; 4)}$.\\
\end{convention}

\subsection{$^*$-Algebras and (pre-)$^*$-representations}

Let $\mathscr{V}$ be an inner product space over $\mathbb{C}$. Recall that a function $T: \mathscr{V} \to \mathscr{V}$ is said to be \emph{adjointable} if there exists a function $T^*: \mathscr{V} \to \mathscr{V}$ such that
\begin{align*}
 \langle T\xi \,,\, \eta \rangle = \langle \xi \,,\, T^* \eta \rangle\,,
\end{align*}
for all $\xi, \eta \in \mathscr{V}$. Recall also that every adjointable operator $T$ is necessarily linear and that $T^*$ is unique and adjointable with $T^{**} = T$. We will use the following notation:
\begin{itemize}
 \item $L(\mathscr{V})$\index{L(V)@$L(\mathscr{V})$} denotes the $^*$-algebra of all adjointable operators in $\mathscr{V}$
 \item $B(\mathscr{V})$\index{B(V)@$B(\mathscr{V})$} denotes the $^*$-algebra of all bounded adjointable operators in $\mathscr{V}$.
\end{itemize}
Of course, we always have $B(\mathscr{V}) \subseteq L(\mathscr{V})$, with both $^*$-algebras coinciding when $\mathscr{V}$  is a Hilbert space (see, for example, \cite[Proposition 9.1.11]{palmer}).

Following \cite[Def. 9.2.1]{palmer}, we define a \emph{pre-$^*$-representation} of a $^*$-algebra $A$ on an inner product space $\mathscr{V}$ to be a $^*$-homomorphism $\pi: A \to L(\mathscr{V})$ and a \emph{$^*$-representation} of $A$ on a Hilbert space $\mathscr{H}$ to be a $^*$-homomorphism $\pi: A \to B(\mathscr{H})$. As in \cite[Def. 4.2.1]{palmer11}, a pre-$^*$-representation $\pi:A \to L(\mathscr{V})$ is said to be \emph{normed} if $\pi(A) \subseteq B(\mathscr{V})$, i.e. if $\pi(a)$ is a bounded operator for all $a \in A$. \\

\begin{df}[\cite{palmer}, Def. 10.1.17]
 A $^*$-algebra $A$ is called a \emph{$BG^*$-algebra}  if all pre-$^*$-representations of $A$ are normed.\\
\end{df}

We now introduce our notion of an essential ideal. Our definition is not the usual one, but this choice of terminology will be justified in what follows.\\

\begin{df} Let $A$ be a $^*$-algebra. An ideal $I \subseteq A$ is said to be \emph{essential} if $aI \neq \{0\}$ for all $a \in A \setminus \{0\}$.\\
\end{df}

The usual definition of an essential ideal states that $I$ is essential if it has nonzero intersection with every other nonzero ideal. Our definition is stronger, but coincides with the usual one for the general class of semiprime $^*$-algebras. We recall from \cite[Definition 4.4.1]{palmer11} that a $^*$-algebra is said to be \emph{semiprime} if $aAa = \{0\}$ implies $a = 0$, where $a \in A$. The class of semiprime $^*$-algebras is quite large, containing all $^*$-algebras that have a faithful $^*$-representation on a Hilbert space (in particular, all $C^*$-algebras) and many other classes of $^*$-algebras (see \cite[Theorem 9.7.21]{palmer}).\\

\begin{prop}
\label{equiv semiprime}
 Let $A$ be an algebra and $I \subseteq A$ a nonzero ideal. We have
\begin{itemize}
 \item[i)] If $I$ is essential, then $I$ has a nonzero intersection with every other nonzero ideal of $A$.
 \item[ii)] The converse of $i)$ is true in case $A$ is semiprime.\\
\end{itemize}

\end{prop}

{\bf \emph{Proof:}} $i)$  Let $I$ be an essential ideal of $A$. Let $J \subseteq A$ be a nonzero ideal and $a \in J \setminus \{0\}$. Since $a$ is nonzero, then $a I \neq \{0\}$. Hence, $J \cdot I \neq \{0\}$, and since $J \cdot I \subseteq J \cap I$, we have $J \cap I \neq \{0\}$.

$ii)$ Suppose $A$ is semiprime. Suppose also that $I$ is not essential. Thus, there is $a \in A \setminus \{0\}$ such that $a I = \{0\}$. Let $J_a \subseteq A$ be the ideal generated by $a$. We have $J_a \cdot I = \{0\}$. Since $(J_a \cap I)^2 \subseteq J_a \cdot I$ we have $(J_a \cap I)^2 = \{0\}$.  Since $A$ is semiprime this implies that $J_a \cap I = \{0\}$ (see \cite[Theorem 4.4.3]{palmer11}). Hence, $I$ has zero intersection with a nonzero ideal. \qed\\

For $C^*$-algebras the focus is mostly on closed ideals. In this setting we still see that our definition is equivalent to the usual one (\cite[Definition 2.35]{morita equiv}):\\

\begin{prop}
 Let $A$ be a $C^*$-algebra and $I \subseteq A$ a closed ideal. The following are equivalent:
\begin{itemize}
 \item[i)] $I$ is essential.
 \item[ii)] $I$ has nonzero intersection with every other nonzero ideal of $A$.
 \item[iii)] $I$ has nonzero intersection with every other nonzero closed ideal of $A$.\\
\end{itemize}

\end{prop}

{\bf \emph{Proof:}} $i) \Longleftrightarrow ii)$ This was established in Proposition \ref{equiv semiprime}, since $C^*$-algebras are automatically semiprime.

$ii) \Longrightarrow iii)$ This is obvious.

$ii) \Longleftarrow iii)$ Let $J$ be a nonzero ideal of $A$ and $\overline{J}$ its closure. From $iii)$ we have $I \cap \overline{J} \neq \{0\}$. Since $I$ and $\overline{J}$ are both closed, and $A$ is a $C^*$-algebra, we have $I\cdot \overline{J} = I \cap \overline{J}$. Now, it is clear that $I \cdot J = \{0\}$ if and only if $I\cdot \overline{J} = \{0\}$. Hence, we necessarily have $I \cdot J \neq \{0\}$, which implies $I \cap J \neq \{0\}$. \qed\\

We now introduce the notion of an \emph{essential} $^*$-algebra. The class of essential $^*$-algebras seems to be the appropriate class of $^*$-algebras for which one can a define a multiplier algebra (as we shall see in Section \ref{algebraic multiplier algebras section}).\\

\begin{df}
 A $^*$-algebra $A$ is said to be \emph{essential} if $A$ is an essential ideal of itself, i.e. if $aA \neq \{0\}$ for all $a \in A \setminus \{0\}$.\\
\end{df}

Any unital $^*$-algebra is obviously essential. Also, it is easy to see that a semiprime $^*$-algebra is essential. The converse is false, so that essential $^*$-algebras form a more general class than that of semiprime $^*$-algebras:\\

\begin{ex}
Let $\mathbb{C}[X]$ be the polynomial algebra in one selfadjoint variable $X$. For any $n \geq 2$ the algebra $\mathbb{C}[X] / \langle X^n \rangle$ is essential, because it is unital, but it is not semiprime because $[X^{n-1}] \Big(\mathbb{C}[X] / \langle X^n \rangle \Big) [X^{n-1}] = \{0\}$.\\
\end{ex}

\subsection{$^*$-Algebraic multiplier algebras}
\label{algebraic multiplier algebras section}

Every $C^*$-algebra can be embedded in a unital $C^*$-algebra in a ``maximal'' way. These maximal unitizations of $C^*$-algebras enjoy a number of useful properties and certain concrete realizations of these algebras are commonly referred to as multiplier algebras. The reader is referred to \cite{morita equiv} for an account.

The definition of a multiplier algebra is quite standard in $C^*$-algebra theory, but this notion is in fact more general and applicable for more general types of rings and algebras. For example, in \cite[Section 1.1]{ara} it is explained how multiplier algebras can be defined for semiprime algebras.

In this section we are going to generalize this notion to the context of essential $^*$-algebras and derive their basic properties. We believe that essential $^*$-algebras are the appropriate class of $^*$-algebras for which one can define multiplier algebras, since the property $aA = \{0\} \Rightarrow a = 0$, which characterizes an essential $^*$-algebra, is constantly used in proofs.

Multiplier algebras are many times defined via the so-called double centralizers (see for example \cite{ara}), but since we are only interested in algebras with an involution a slightly simpler and more convenient approach can be given, analogue to the Hilbert $C^*$-module approach to $C^*$-multiplier algebras (presented in \cite[Section 2.3]{morita equiv}). This is the approach we follow.\\

\begin{df}
\label{max unit of *-alg}
Let $\mathcal{C}$ be a subclass of $^*$-algebras. A $^*$-algebra $A \in \mathcal{C}$ is said to have a \emph{maximal unitization in $\mathcal{C}$} if there exists a unital $^*$-algebra $B \in \mathcal{C}$ (called the \emph{maximal unitization} of $A$) and a $^*$-embedding $i: A \hookrightarrow B$ for which $i(A)$ is an essential ideal of $B$ and such that for every other $^*$-embedding $j$ of $A$ as an essential ideal of a unital $^*$-algebra $C \in \mathcal{C}$, there is a unique $^*$-homomorphism $\phi: C \to B$ such that
\begin{displaymath}
\xymatrix{ & B \\
A \ar[ur]^i \ar[r]_j & C \ar@{.>}[u]_{\phi}} 
\end{displaymath}
commutes.\\
\end{df}

\begin{lemma}
\label{phi is injective in diagram}
 In the above diagram the $^*$-homomorphism $\phi$ is always injective (even if $C$ was not unital).\\
\end{lemma}

{\bf \emph{Proof:}} We have that $j(A) \cap \mathrm{Ker}\, \phi = \{0\}$, because if $j(a) \in j(A) \cap \mathrm{Ker}\, \phi$, then $0 = \phi (j(a)) = i(a)$ and hence $a = 0$ and therefore $j(a)=0$. Hence, since $j(A)$ is an essential ideal of $C$, it follows from Proposition \ref{equiv semiprime} $i)$ that $\mathrm{Ker}\, \phi = \{0\}$. $\qed$\\

For $C^*$-algebras, one might expect to replace ``ideal'' by ``closed ideal'', in Definition \ref{max unit of *-alg}. This condition, however, follows automatically since $i(A)$ and $j(A)$ are automatically closed. Hence, this definition encompasses the usual definition of a maximal unitization for a $C^*$-algebra.\\

\begin{df}
 Let $A$ be a $^*$-algebra. By a \emph{right $A$-module} we mean a vector space $X$ together with a mapping $X \times A \to X$ satisfying the usual consistency conditions. An \emph{$A$-linear mapping} $T: X \to Y$ between $A$-modules is a linear mapping between the underlying vector spaces such that $T(xa) = T(x)a$, for all $x \in X$ and $a \in A$. We will often use the notation $Tx$, instead of $T(x)$.\\
\end{df}

Every $^*$-algebra $A$ is canonically a right $A$-module, with the action of right multiplication. This is the example we will use thoroughly in what follows.

Let $\langle \cdot , \cdot \rangle_A: A \times A \to A$ be the function
\begin{align*}
 \langle a, b \rangle_A := a^*b\,.
\end{align*}
The function $\langle \cdot , \cdot \rangle_A$ is an $A$-linear form, in the sense that the following properties are satisfied:
\begin{enumerate}
 \item[a)] $\langle a\,,\, \lambda_1 b_1 + \lambda_2 b_2 \rangle_A = \lambda_1 \langle a, b_1 \rangle_A + \lambda_2 \langle a, b_2 \rangle_A\,$,
 \item[b)] $\langle \lambda_1a_1 + \lambda_2 a_2\,,\, b \rangle_A = \overline{\lambda_1} \langle a_1, b \rangle_A + \overline{\lambda_2} \langle a_2, b \rangle_A\,$,
 \item[c)] $\langle a, b c \rangle_A = \langle a, b \rangle_A\, c\,$,
 \item[d)] $\langle ac, b \rangle_A = c^*\langle a, b \rangle_A\,$,
 \item[e)] $\langle a, b \rangle_A^* = \langle b, a \rangle_A\,$,
\end{enumerate}
for all $a, a_1, a_2, b, b_1, b_2 \in A$ and $\lambda_1, \lambda_2 \in \mathbb{C}$.\\

If the $^*$-algebra $A$ is essential we also have:
\begin{enumerate}
 \item[f)] If $\langle a\,,\, b \rangle_A = 0$ for all $b \in A$, then $a = 0\,$.\\
\end{enumerate}

\begin{df}
 Let $A$ be a $^*$-algebra. A function $T: A \to A$ is called \emph{adjointable} if there is a function $T^*: A \to A$ such that
\begin{align*}
 \langle T(a), b\rangle_A = \langle a, T^*(b) \rangle_A\,,
\end{align*}
for all $a, b \in A$.\\
\end{df}

\begin{prop}
 If $A$ is an essential $^*$-algebra, then every adjointable map $T: A \to A$ is $A$-linear. Moreover, the adjoint $T^*$ is unique and adjointable with $T^{**} = T$.\\
\end{prop}

{\bf \emph{Proof:}} Let $T$ be an adjointable map in $A$ and $x_1, x_2, y \in A$. We have
\begin{eqnarray*}
\langle T(\lambda_1x_1+\lambda_2x_2)\,,\, y \rangle_A & = & \langle \lambda_1x_1+ \lambda_2 x_2\,,\, T^*(y) \rangle_A\\
& = & \overline{\lambda_1}\,\langle x_1\,,\, T^*(y) \rangle_A + \overline{\lambda_2}\,\langle x_2\,,\, T^*(y) \rangle_A\\
& = & \overline{\lambda_1}\,\langle T(x_1)\,,\, y \rangle_A + \overline{\lambda_2}\,\langle T(x_2)\,,\, y \rangle_A\\
& = & \langle \lambda_1T(x_1)+\lambda_2T(x_2)\,,\, y \rangle_A\,.
\end{eqnarray*}
Hence, we have $\langle T(\lambda_1x_1+\lambda_2x_2)- \lambda_1T(x_1)+\lambda_2T(x_2)\,,\, y \rangle_A = 0$. We can then  conclude from f) that
\begin{align*}
 T(\lambda_1x_1+\lambda_2x_2)- \lambda_1T(x_1)+\lambda_2T(x_2) = 0\,,
\end{align*}
i.e. $T$ is a linear map.

Let us now check that $T$ is $A$-linear. For any $x, y, a \in A$ we have
\begin{eqnarray*}
\langle T(xa)\,,\, y \rangle_A & = & \langle xa\,,\, T^*(y) \rangle_A =  a^* \langle x\,,\, T^*(y) \rangle_A\\
& = & a^* \langle T(x)\,,\, y \rangle_A =  \langle T(x) a\,,\, y \rangle_A\,.
\end{eqnarray*}
Hence, we have $\langle T(xa) - T(x)a\,,\, y \rangle_A = 0$. We can then conclude from f) that $T(xa)- T(x)a= 0$, i.e. $T$ is $A$-linear.

Let us now prove the uniqueness of the adjoint $T^*$. Suppose there was a function $S: A \to A$ such that
\begin{align*}
 \langle x \,,\, T^*(y) \rangle_A = \langle x\,,\, S(y) \rangle_A\,.
\end{align*}
for all $x, y \in A$. Then, $\langle T^*(y) - S(y) \,,\, x  \rangle_A = 0$. We can then conclude from f) that $T^*(y) - S(y) = 0$, i.e. $T^*= S$.

 It remains to prove that $T^*$ is adjointable with $T^{**} = T$. This follows easily from the equality
\begin{align*}
 \langle T^* x\,,\, y \rangle_A  =  \langle y\,,\, T^*x \rangle_A^* = \langle Ty\,,\,x \rangle_A^* =  \langle x, Ty \rangle_A\,.
\end{align*}
\qed\\

\begin{df}
Let $A$ be an essential $^*$-algebra. The set of all adjointable maps on $A$ is called the \emph{multiplier algebra} of $A$ and is denoted by $M(A)$\index{M(A)@$M(A)$}.\\
\end{df}

The multiplier algebra is in fact a $^*$-algebra, and the proof of this fact is standard.\\

\begin{prop}
Let $A$ be an essential $^*$-algebra. The multiplier algebra of $A$ is a unital $^*$-algebra with the sum and multiplication given by pointwise sum and composition (respectively), and the involution given by the adjoint.\\
\end{prop}

\begin{prop}
\label{A is an ideal of M(A)}
Let $A$ be an essential $^*$-algebra. There is a $^*$-embedding $L: A \to M(A)$ of $A$ as an essential ideal of $M(A)$, given by
\begin{align*}
 a \mapsto L_a
\end{align*}
where $L_a: A \to A$ is the left multiplication by $a$, i.e. $L_a(b):=ab$.\\
\end{prop}

{\bf \emph{Proof:}} It is easy to see that, for every $a \in A$, $L_a$ is adjointable with adjoint $L_{a^*}$, thus the mapping $L$ is well-defined. Also clear is the fact that $L$ is a $^*$-homomorphism. Let us prove that it is injective: suppose $L_a = 0$ for some $a \in A$. Then, for all $b \in A$ we have $ab = L_a b = 0$ and since $A$ is essential this implies $a =0$. Thus, $L$ is injective.

It remains to prove that $L(A)$ is an essential ideal of $M(A)$. Let us begin by proving that it is an ideal. Let $T \in M(A)$. For every $a, b \in A$ we have
\begin{align*}
 TL_a (b)= T(ab) = T(a)b = L_{Ta} (b)\,,
\end{align*}
and also
\begin{eqnarray*}
 L_aT(b) & = & aT(b)   =  \langle a^*, T(b) \rangle\\
 & = & \langle T^*(a^*), b \rangle = (T^*(a^*))^*b\\
& = & L_{(T^*a^*)^*} (b)\,.
\end{eqnarray*}
Hence we have
\begin{align}
\label{properties of L in multiplier algebra}
 TL_a = L_{Ta} \qquad \text{and} \qquad L_aT = L_{(T^*a^*)^*}\,,
\end{align}
from which it follows easily that $L(A)$ is an ideal of $M(A)$.

 Let us now prove that this ideal is essential. Let $T \in M(A)$ be such that $TL(A) = \{0\}$. Then, in particular, $TL_a = 0$ for all $a \in A$, but as we have seen before $TL_a = L_{Ta}$, and since $L$ is injective we must have $Ta = 0$ for all $a \in A$, i.e $T = 0$. \qed\\

\begin{rem}
According to Proposition \ref{A is an ideal of M(A)}, an essential $^*$-algebra $A$ is canonically embedded in its multiplier algebra $M(A)$. We will often make no distinction of notation between $A$ and its embedded image in $M(A)$, i.e. we will often just write $a$ to denote an element of $A$ and to denote the element $L(a)$ of $M(A)$.
 No confusion will arise from this because the left equality in (\ref{properties of L in multiplier algebra}) simply means, in this notation, that $T \cdot a = T(a)$.\\
\end{rem}

\begin{thm}
 \label{property of M(A)}
Let $A$ be an essential $^*$-algebra and $L: A \to M(A)$ the canonical $^*$-embedding of $A$ in $M(A)$. If $j:A \to C$ is a $^*$-embedding of $A$ as an ideal of a $^*$-algebra $C$, then there exists a unique $^*$-homomorphism $\phi: C \to M(A)$ such that the following diagram commutes
\begin{displaymath}
\xymatrix{ & M(A) \\
A \ar[ur]^L \ar[r]_j & C \ar@{.>}[u]_{\phi}} 
\end{displaymath}
Moreover, if $j(A)$ is essential then $\phi$ is injective. \\
\end{thm}

{\bf \emph{Proof:}} For simplicity of notation let us assume, without any loss of generality, that $A$ itself is an ideal of a $^*$-algebra $C$, so that we avoid any reference to $j$ (or its inverse).  Let $\phi : C \to M(A)$ be the function defined by
\begin{align*}
 \phi(c): A \to A\\
\phi(c) a := ca\,,
\end{align*}
for every $c \in C$. It is a straightforward computation to check that $\phi(c) \in M(A)$ and that $\phi$ itself is a $^*$-homomorphism. It is also easy to see that $\phi(a) = L_a$, for every $a \in A$. Hence, $\phi \circ j = L$. Let us now prove the uniqueness of $\phi$ relatively to this property. Suppose $\widetilde{\phi}: C \to M(A)$ is another $^*$-homomorphism such that $\widetilde{\phi} \circ j = L$. Then, for all $c \in C$ and $a \in A$ we have
\begin{eqnarray*}
 \big(\widetilde{\phi}(c) - \phi(c) \big) L_a & = & \widetilde{\phi}(c)L_a - \phi(c)L_a\\
& = & \widetilde{\phi}(c)\widetilde{\phi}(a) - \phi(c)\phi(a)\\
& = & \widetilde{\phi}(ca) - \phi(ca)\\
& = & L_{ca} - L_{ca}\\
& = & 0\,.
\end{eqnarray*}
Since $L(A)$ is an essential ideal of $M(A)$ it follows that $\widetilde{\phi}(c) = \phi(c)$ for all $c \in C$, i.e. $\widetilde{\phi} = \phi$.

The last claim of the theorem, concerning injectivity of $\phi$, was proven in Lemma \ref{phi is injective in diagram}. \qed\\

\begin{cor}
\label{mult is max unit}
The multiplier algebra $M(A)$ is a maximal unitization of $A$ in the class of: essential $^*$-algebras, semiprime $^*$-algebras and $C^*$-algebras.\\
\end{cor}

{\bf \emph{Proof:}} By Theorem \ref{property of M(A)} we only need to check that if $A$ is an essential $^*$-algebra (respectively, semiprime $^*$-algebra or $C^*$-algebra), then the multiplier algebra has the same property.

Suppose $A$ is an essential $^*$-algebra. Let $T \in M(A)$ be such that $TM(A) = \{0\}$. Then, by the embedding of $A$ in $M(A)$ we have $Ta= 0$ for all $a \in A$, i.e. $T = 0$. Hence, $M(A)$ is also an essential $^*$-algebra.

Suppose $A$ is a semiprime $^*$-algebra. Let $T \in M(A)$ be such that $TM(A)T = \{0\}$. Then, we also have that $TL_a M(A)TL_a= \{0\}$ for any $a \in A$, and therefore $L_{T(a)}M(A) L_{T(a)} = \{0\}$. Thus, in particular, $L_{T(a)}L(A) L_{T(a)} = \{0\}$, and since $L$ is injective this means that $T(a) A T(a) = \{0\}$. Since $A$ is semiprime we conclude that $T(a) = 0$, and therefore $T = 0$. Hence, $M(A) $ is semiprime.

It is well-known that $M(A)$ is a $C^*$-algebra when $A$ is a $C^*$-algebra. \qed\\

An important feature of $C^*$-multiplier algebras is that a nondegenerate $^*$-representation of $A$ extends uniquely to $M(A)$. This result does not hold in general for essential $^*$-algebras. Nevertheless we can still extend a nondegenerate $^*$-representation of $A$ to a unique pre-$^*$-representation of $M(A)$:\\

\begin{thm}
\label{ext hom *-alg}
Let $A$ be an essential $^*$-algebra, $\pi: A \to B(\mathscr{H})$ a nondegenerate $^*$-representation of $A$ on a Hilbert space $\mathscr{H}$ and $\mathscr{V} \subseteq \mathscr{H}$ the dense subspace
\begin{align*}
 \mathscr{V}:=\pi(A)\mathscr{H} = \mathrm{span}\; \{\pi(a)\xi: a\in \a\,, \xi \in \mathscr{H}\}\,.
\end{align*}
Then there is a unique pre-$^*$-representation
\begin{align*}
 \widetilde{\pi}\index{pitilde@$\tilde{\pi}$}: M(A) \to L(\mathscr{V})
\end{align*}
such that $\widetilde{\pi} (a) = \pi(a)|_{\mathscr{V}}$ for every $a \in A$.\\
\end{thm}

{\bf \emph{Proof:}} We define the pre-$^*$-representation $\widetilde{\pi}: M(A) \to L(\mathscr{V})$ by
\begin{align*}
 \widetilde{\pi} (T)\;\big[\sum_{i=1}^n \pi(a_i) \xi_i \big]:= \sum_{i=1}^n \pi(Ta_i)\xi_i\,,
\end{align*}
for $n \in \mathbb{N}$, $a_1, \dots, a_n \in A$ and $\xi_1, \dots, \xi_n \in \mathscr{H}$. Let us first check that $\widetilde{\pi}$ is well-defined. Suppose $\sum_{i=1}^n \pi(a_i)\xi_i = \sum_{j=1}^m \pi(b_j) \eta_j$. Then, for every $z \in A$ we have
\begin{eqnarray*}
 \pi(z) \Big( \sum_{i=1}^n \pi(Ta_i)\xi_i - \sum_{j=1}^m \pi(Tb_i)\eta_j \Big) & = & \sum_{i=1}^n \pi(zTa_i)\xi_i - \sum_{j=1}^m \pi(zTb_i)\eta_j\\
& = & \pi(zT) \Big( \sum_{i=1}^n \pi(a_i)\xi_i - \sum_{j=1}^m \pi(b_i)\eta_j \Big)\\
& = & 0\,.
\end{eqnarray*}
Since the $^*$-representation $\pi$ is nondegenerate we necessarily have
\begin{align*}
 \sum_{i=1}^n \pi(Ta_i)\xi_i - \sum_{j=1}^m \pi(Tb_i)\eta_j = 0\,,
\end{align*}
which means that $\widetilde{\pi}(T)$ is well-defined.

Let us now check that $\widetilde{\pi}(T) \in L(\mathscr{V})$, i.e. that $\widetilde{\pi}(T)$ is indeed an adjointable operator in $\mathscr{V}$. We will in fact prove that $\widetilde{\pi}(T)^*=\widetilde{\pi}(T^*)$, which follows from the following equality
\begin{eqnarray*}
 \Big\langle \widetilde{\pi}(T) \sum_{i=1}^n \pi(a_i)\xi_i\,\,,\, \sum_{j=1}^m \pi(b_i)\eta_j \Big\rangle & = &  \sum_{i=1}^n \sum_{j=1}^m \langle \pi(Ta_i)\xi_i\,,\,\pi(b_j)\eta_j \rangle\\
& = & \sum_{i=1}^n \sum_{j=1}^m \langle \xi_i\,,\,\pi(a_i^*T^*)\pi(b_j)\eta_j \rangle\\
& = & \sum_{i=1}^n \sum_{j=1}^m \langle \xi_i\,,\,\pi(a_i^*T^*b_j)\eta_j \rangle\\
& = & \sum_{i=1}^n \sum_{j=1}^m \langle \pi(a_i)\xi_i\,,\,\pi(T^*b_j)\eta_j \rangle\\
& = & \Big\langle \sum_{i=1}^n \pi(a_i)\xi_i\,\,,\, \widetilde{\pi}(T^*) \sum_{j=1}^m \pi(b_i)\eta_j \Big\rangle\,.
\end{eqnarray*}
It is straightforward to see that $\widetilde{\pi}$ is linear, multiplicative and, as we have seen, $\widetilde{\pi}(T^*) = \widetilde{\pi}(T)^*$, hence $\widetilde{\pi}$ is a pre-$^*$-representation of $M(A)$ on $\mathscr{V}$.

It is also clear that, for any $a \in A$, $\widetilde{\pi}(a)$ is just $\pi(a)$ restricted to $\mathscr{V}$,  because of the equality
\begin{align*}
 \widetilde{\pi}(a) \sum_{i=1}^n \pi(a_i)\xi_i = \sum_{i=1}^n \pi(aa_i)\xi_i = \pi(a) \sum_{i=1}^n \pi(a_i)\xi_i\,.
\end{align*}
Let us now prove the uniqueness of $\widetilde{\pi}$. Suppose $\phi: M(A) \to L(\mathscr{V})$ is a pre-$^*$-representation such that $\phi(a) = \pi(a)|_{\mathscr{V}}$. Then, for every $z \in A$ and $v \in \mathscr{V}$ we have
\begin{eqnarray*}
 \pi(z) ( \phi(T)v -\widetilde{\pi}(T)v ) & = & \pi(z)\phi(T)v - \pi(z)\widetilde{\pi}(T)v\\
& = & \phi(z)\phi(T)v - \widetilde{\pi}(z)\widetilde{\pi}(T)v\\
& = & \phi(zT)v - \widetilde{\pi}(zT)v\\
& = & \pi(zT)v - \pi(zT)v\\
& = & 0\,.
\end{eqnarray*}
Since the $^*$-representation $\pi$ is nondegenerate, we necessarily have
\begin{align*}
 \phi(T)v -\widetilde{\pi}(T)v = 0\,,
\end{align*}
which means that $\phi(T) = \widetilde{\pi}(T)$, i.e. $\phi = \widetilde{\pi}$. \qed\\

\begin{rem}
  Theorem \ref{ext hom *-alg} can be interpreted in the following way: every nondegenerate $^*$-representation $\pi: A \to B(\mathscr{H})$ can be extended to $M(A)$ by possibly unbounded operators, defined on the dense subspace $\pi(A)\mathscr{H}$.\\
\end{rem}

\begin{df}
 Let $A$ be an essential $^*$-algebra. We will denote by $M_B(A)$\index{MB(A)@$M_B(A)$} the subset of $M(A)$ consisting of all the elements $T \in M(A)$ such that $\widetilde{\pi}(T) \in B(\mathscr{V})$ for all nondegenerate $^*$-representations $\pi:A \to B(\mathscr{H})$, where $\mathscr{V}:=\pi(A)\mathscr{H}$ and $\widetilde{\pi}$ is the unique pre-$^*$-representation extending $\pi$ as in Proposition \ref{ext hom *-alg}.\\
\end{df}

As stated in the next result, $M_B(A)$ is a $^*$-subalgebra of $M(A)$. The advantage of working with $M_B(A)$ over $M(A)$ is that nondegenerate $^*$-representations of $A$ always extend to $^*$-representations of $M_B(A)$. Easy examples of elements of $M_B(A)$ that might not belong to $A$ are the projections and unitaries of $M(A)$.\\

\begin{prop}
\label{ext reps to MB(A)}
 Let $A$ be an essential $^*$-algebra. The set $M_B(A)$ is a $^*$-subalgebra of $M(A)$ containing $A$. Moreover, if $\pi: A \to B(\mathscr{H})$ is a nondegenerate $^*$-representation of $A$, then there exists a unique $^*$-representation of $M_B(A)$ on $\mathscr{H}$ that extends $\pi$.\\
\end{prop}

{\bf \emph{Proof:}} Let $T, S \in M_B(A)$.  Let $\pi :A \to B(\mathscr{H})$ be any nondegenerate $^*$-representation of $A$ and $\widetilde{\pi}$ its extension to $L(\mathscr{V})$, in the sense of Theorem \ref{ext hom *-alg}, where  $\mathscr{V}:= \pi(A)\mathscr{H}$. By definition, $\widetilde{\pi}(T), \widetilde{\pi}(S) \in B(\mathscr{V})$, and therefore  $\widetilde{\pi}(T+S), \widetilde{\pi}(TS), \widetilde{\pi}(T^*) \in B(\mathscr{V})$, since $B(\mathscr{V})$ is a $^*$-algebra. Hence, $M_B(A)$ is a $^*$-subalgebra of $M(A)$. Moreover, $A \subseteq M_B(A)$ since $\widetilde{\pi}(a) = \pi(a)|_{\mathscr{V}} \in B(\mathscr{V})$.

Let us now prove the last claim of this proposition. Let $\pi: A \to B(\mathscr{H})$ be a nondegenerate $^*$-representation and $\widetilde{\pi}: M(A) \to L(\mathscr{V})$ its extension as in Theorem \ref{ext hom *-alg}. Then we obtain by restriction a pre-$^*$-representation $\widetilde{\pi}:M_B(A) \to L(\mathscr{V})$. By definition of $M_B(A)$ we actually have $\widetilde{\pi}(M_B(A)) \subseteq B(\mathscr{V})$. Hence $\widetilde{\pi}$ gives rise to a $^*$-representation $\widetilde{\pi}: M_B(A) \to B(\mathscr{H})$, since $\mathscr{V}$ is dense in $\mathscr{H}$.

 Let us now prove the uniqueness claim. Suppose $\varphi$ is another representation of $M_B(A)$ that extends $\pi$. For $T \in M_B(A)$, $a \in A$ and $\xi \in \mathscr{H}$ we have
\begin{eqnarray*}
 \varphi(T)\,\pi(a)\xi & = & \varphi(T)\, \varphi(a)\xi \;\; = \;\; \varphi(Ta)\xi\\
& = &  \pi(Ta)\xi \;\; =\;\; \widetilde{\pi}(T) \, \pi(a) \xi\,.
\end{eqnarray*}
By linearity and density it follows that $\varphi(T) = \widetilde{\pi}(T)$, i.e. $\varphi = \widetilde{\pi}$. \qed\\

The above result is a generalization of the well-known result for $C^*$-algebras which states that any nondegenerate $^*$-representation can be extended to the multiplier algebra (see for example \cite[Corollary 2.51]{morita equiv}), because $M(A) = M_B(A)$ for any $C^*$-algebra $A$.\\

\subsection{Hecke algebras}

We start by establishing some notation and conventions concerning left coset spaces and double coset spaces and we prove two resuls which will be useful later on.

Let $G$ be a group,  $B, C$ subgroups of $G$ and $e \in G$ the identity element. The double coset space $B \backslash G / C$ is the set
\begin{align}
 B \backslash G / C := \{BgC \subseteq G : g \in G \}\,.
\end{align}
It is easy to see that the sets of the form $BgC$ are either equal or disjoint, or in other words, we have an equivalence relation defined in $G$ whose equivalence classes are precisely the sets $BgC$.

The left coset space $G / C$ is the set
\begin{align}
  G / C := \{e\} \backslash G / C = \{gC \subseteq G : g \in G \}\,.
\end{align}

Given an element $g \in G$ and a double coset space $B \backslash G / C$  (which can in particular be a left coset space by taking $B = \{e\}$) we will denote by $[g]$ the double coset $BgC$. Thus, $[g]$ denotes the whole equivalence class for which $g\in G$ is a representative.

If $A$ is a subset of $G$ we define the double coset space $B \backslash A / C$ as the set of double cosets in  $B \backslash G / C$ which have a representative in $A$, i.e.
\begin{align}
 B \backslash A / C := \{BaC \subseteq G : a \in A \}\,.
\end{align}\\
 
\begin{prop}
\label{double coset spaces prop}
 Let $A,B$ and $C$ be subgroups of a group $G$. If $C \subseteq A$, then the following map is a bijective correspondence between the double coset spaces:
 \begin{align}
 \label{map BAC to BcapA AC}
  B\backslash A / C &\longrightarrow (B \cap A) \backslash A / C
  \end{align}
  \vskip-1cm
 \begin{align*}
 [a] & \mapsto [a]\,.\qquad
 \end{align*}
 Similarly, if $B \subseteq A$, then the following map is a bijective correspondence:
 \begin{align}
 \label{map BAC to BA AcapC}
  B\backslash A / C &\longrightarrow B \backslash A / (A \cap C)
  \end{align}
  \vskip-1cm
 \begin{align*}
 [a] & \mapsto [a]\,.\qquad\\
 \end{align*}
\end{prop}

{\bf \emph{Proof:}} We first need to show that the map (\ref{map BAC to BcapA AC}) is well defined, i.e. if $Ba_1 C = Ba_2 C$, for some $a_1, a_2 \in A$, then $(B \cap A)a_1 C = (B \cap A)a_2 C$. If $Ba_1 C = Ba_2 C$ then there exist $b \in B$ and $c \in C$ such that $a_1 = ba_2c$, from which it follows that $b = a_1c^{-1}a_2^{-1}$. Since $A$ is a subgroup and $C \subseteq A$, it follows readily that $b \in B \cap A$, and therefore $a_1 \in (B \cap A)a_2 C$, i.e. $(B \cap A)a_1 C = (B \cap A)a_2 C$.

The map defined in (\ref{map BAC to BcapA AC}) is clearly surjective. It is also injective because if $(B \cap A)a_1 C = (B \cap A)a_2 C$, then clearly $Ba_1 C = Ba_2 C$.

A completely analogous argument shows that map defined in (\ref{map BAC to BA AcapC}) is a bijection. \qed\\

Suppose a group $G$ acts (on the right) on a set $X$ and let $x \in  X$. We will henceforward denote by $\mathcal{S}_x$\index{s2 tabilizer@$\mathcal{S}_x$} the stabilizer of the point $x$, i.e.
\begin{align}
\label{stabilizer notation}
 \mathcal{S}_x := \{g \in G: xg = x \}\,.
\end{align}
Given a subset $Z \subseteq X$ and a subgroup $H \subseteq G$ we denote by $Z / H$ the set of $H$-orbits which have representatives in $Z$, i.e.
\begin{align*}
Z / H := \{zH: z \in Z \}\,.
\end{align*}
Suppose now that $H, K \subseteq G$ are subgroups and let $x \in X$ be a point. The following result establishes a correspondence between the set of $H$-orbits $(xK) / H$ and the double coset space $\mathcal{S}_x \backslash K / H$:\\

\begin{prop}
\label{xKH bijection with Sx K H}
 Let $G$ be a group which acts (on the right) on a set $X$. Let $x \in X$ be a point and $H, K \subseteq G$ be subgroups. We have a bijection
\begin{align*}
 (xK) / H \longrightarrow \mathcal{S}_x \backslash K / H\,,
\end{align*}
given by $xgH \mapsto \mathcal{S}_{x} g H$, where $g \in K$.\\
\end{prop}

{\bf \emph{Proof:}} Let us first prove that the map $xgH \mapsto \mathcal{S}_{x} g H$ is well defined, i.e. if $xg_1 H = xg_2 H $, then $ \mathcal{S}_{x} g_1 H =  \mathcal{S}_{x} g_2 H$. If $xg_1 H = xg_2 H $, then there exists $h \in h$ such that $xg_1 = xg_2 h$, which implies that $x = x g_2 h g_1^{-1}$, from which it follows that $g_2 h g_1^{-1} \in \mathcal{S}_x$. Thus we see that
\begin{eqnarray*}
\mathcal{S}_{x} g_1 H & = & \mathcal{S}_{x} g_2 h g_1^{-1} g_1 H \;\; = \;\; \mathcal{S}_{x} g_2 H\,.
\end{eqnarray*}
We conclude that the map is well-defined. The map is obviously surjective. It is also injective because if $ \mathcal{S}_{x} g_1 H =  \mathcal{S}_{x} g_2 H$, then there exists $r \in \mathcal{S}_x$ and $h \in H$ such that $g_1 = r g_2 h$, from which it follows that $xg_1H = xr g_2 hH = xg_2H$. \qed\\

We will mostly follow \cite{krieg} and \cite{schl} in what regards Hecke pairs and Hecke algebras and refer to these references for more details.

We start by establishing some notation which will be useful later on. Given a group $G$, a subgroup $\gm \subseteq G$ and $g \in G$, we will denote by $\gm^g$\index{G1 ammag@$\gm^g$} the subgroup
\begin{align}
\label{notation for intersection with conj}
 \gm^g := \gm \cap g\gm g^{-1}\,.
\end{align}
We now recall the definition of a Hecke pair:\\

\begin{df}
 Let $G$ be a group and $\gm$ a subgroup. The pair $(G , \gm)$\index{G2  gm@$(G, \gm)$} is called a \emph{Hecke pair} if every double coset $\gm g\gm$ is the union of finitely many right (and left) cosets. In this case, $\gm$ is also called a \emph{Hecke subgroup} of $G$.\\
\end{df}

Given a Hecke pair $(G, \gm)$ we will denote by $L$\index{L(g)@$L(g)$} and $R$\index{R(g)@$R(g)$}, respectively, the left and right coset counting functions, i.e.
\begin{align}
 L(g):= |\gm g \gm / \gm| = [\gm : \gm^g] < \infty\\
 R(g) :=|\gm \backslash \gm g \gm| = [\gm : \gm^{g^{-1}}] < \infty\,.
\end{align}
We recall that $L$ and $R$ are $\gm$-biinvariant functions which satisfy $L(g) = R(g^{-1})$ for all $g \in G$. Moreover, the function $\Delta: G \to \mathbb{Q^+}$\index{D1elta@$\Delta(g)$} given by
\begin{align}
\label{def modular function Hecke algebra}
 \Delta(g) := \frac{L(g)}{R(g)}\,,
\end{align}
is a group homomorphism, usually called the \emph{modular function} of $(G, \gm)$.\\

\begin{df}
 Given a Hecke pair $(G, \gm)$, the \emph{Hecke algebra} $\h(G, \gm)$\index{H(G,gm)@$\h(G, \gm)$} is the $^*$-algebra of finitely supported $\mathbb{C}$-valued functions on the double coset space $\gm \backslash G / \gm$ with the product and involution defined by
\begin{align}
 (f_1*f_2)(\gm g \gm) & := \sum_{h\gm \in G / \gm} f_1(\gm h \gm)f_2(\gm h^{-1}g\gm)\,,\\
f^*(\gm g\gm) & := \Delta(g^{-1}) \overline{f(\gm g^{-1} \gm)}\,.\\ \notag
\end{align}

\end{df}

Equivalently, we can define $\mathcal{H}(G, \gm)$ as the $^*$-algebra of finitely supported $\gm$-left invariant functions $f:G / \gm \to \mathbb{C}$ with the product and involution operations given by
\begin{align}
 (f_1*f_2)( g \gm) & := \sum_{h\gm \in G / \gm} f_1( h \gm)f_2( h^{-1}g\gm)\,,\\
f^*( g\gm) & := \Delta(g^{-1}) \overline{f( g^{-1} \gm)}\,.\\ \notag
\end{align}

\begin{rem}
Some authors, including Krieg \cite{krieg}, do not include the factor $\Delta$ in the involution. Here we adopt the convention of Kaliszewski, Landstad and Quigg \cite{schl} in doing so, as it gives rise to a more natural $L^1$-norm. We note, nevertheless, that there is no loss (or gain) in doing so, because these two different involutions give rise to $^*$-isomorphic Hecke algebras.\\
\end{rem}

The Hecke algebra has a natural basis, as a vector space, given by the characteristic functions of double cosets. We will henceforward identify a characteristic function of a double coset $1_{\gm g \gm}$ with the double coset $\gm g \gm$ itself. \\

\subsection{Fell bundles over discrete groupoids}
\label{Fell bundle section}

Let $X$ be a discrete groupoid. We will denote by $X^0$ the unit space of $X$ and by $\mathbf{s}$\index{source@$\mathbf{s}(x)$} and $\mathbf{r}$\index{range@$\mathbf{r}(x)$} the source and range functions $X \to X^0$, respectively.

We will essentially follow \cite{kum} when it comes to Fell bundles over groupoids. All the groupoids in this work are assumed to be discrete, so that the theory of Fell bundles admits a few simplifications. Basically a \emph{Fell bundle} over a discrete groupoid $X$ consists of:
\begin{itemize}
\item a space $\a$ together with a surjective map $p:\a \to X$, such that each \emph{fiber} $\a_x := p^{-1}(x)$ is a Banach space, for every $x \in X$,;
\item a multiplication operation between fibers over composable elements of the groupoid, which we suggestively write as $\a_x \cdot \a_y \subseteq \a_{xy}$;
 \item an involution $a \mapsto a^*$ which takes $\a_x$ into $\a_{x^{-1}}$.
 \end{itemize}
These operations and norms satisfy some consistency properties which are described in \cite[Section 2]{kum}. As it is well-known, each fiber over a unit element is naturally a $C^*$-algebra.\\

\begin{standing}
\label{assumption on Fell bundles}
 Given a Fell bundle $\a$ over a discrete groupoid $X$ we will always assume that the fibers over units are non-trivial, i.e. $\a_u \neq \{ 0\}$ for all $u \in X^0$.\\
\end{standing}

Assumption \ref{assumption on Fell bundles} is not very restrictive. In fact, removing from the groupoid $X$ all the units $u \in X^0$ for which $\a_u = \{0\}$ and also all the elements $x \in X$ such that $\mathbf{s}(x)$ or $\mathbf{r}(x)$ is $u$, we obtain a subgroupoid $Y$ for which the assumption holds (relatively to the restriction $\a |_Y$ of $\a$ to $Y$). Moreover, and this is the important fact, the algebras of finitely supported sections are canonically isomorphic, i.e. $C_c(\a|_Y) \cong C_c(\a)$.

The reason for us to follow Assumption \ref{assumption on Fell bundles} is because it will make our theory slightly simpler. Since we are interested mostly in algebras of sections, this assumption does not reduce the generality of the work in any way, as we observed in the previous paragraph.\\

\begin{df}
 Let $\a$ be a Fell bundle over a discrete groupoid $X$. An \emph{automorphism} of $\a$ is a bijective map $\beta:\a \to \a$ which preserves the bundle structure, i.e. such that
\begin{itemize}
 \item[i)] $\beta$ takes any fiber onto another fiber;
 \item[ii)] $\beta$ takes fibers over composable elements of $X$ to fibers over composable elements;
 \item[iii)] As a map between (two) fibers, $\beta$ is a linear map;
 \item[iv)] $\beta(a\cdot b) = \beta(a)\cdot \beta(b)$, whenever multiplication is defined;
 \item[v)] $\beta(a^*) = \beta(a)^*$.
\end{itemize}
The set of all automorphisms of $\a$ forms a group under composition and will be denoted by $Aut(\a)$.\\
\end{df}

It follows easily from $i)$ and $ii)$ above that every automorphism $\beta$ of $\a$ entails a groupoid automorphism $\beta_0$ of $X$ such that $\beta_0 \big(p(a) \big) = p(\beta(a))$. We also note that, by being a groupoid automorphism, $\beta_0$ takes units into units.\\

\begin{rem}
\label{automorphism gives isometry between fibers}
 The restricted map $\beta: \a_x \to \a_{\beta_0(x)}$ is an isometric linear map. Linearity was required in condition $iii)$, but the fact that the map is an isometry follows from the other axioms. To see this we note that
\begin{eqnarray*}
 \|\beta(a) \| = \| \beta(a)^*\beta(a) \|^{\frac{1}{2}} = \|\beta(a^*a)\|^{\frac{1}{2}}\,.
\end{eqnarray*}
Now $a^*a \in \a_{\mathbf{s}(x)}$ and $\mathbf{s}(x) \in X^0$. Thus, we also have $\beta_0(\mathbf{s}(x)) \in X^0$ and therefore both $\a_{\mathbf{s}(x)}$ and $\a_{\beta_0(\mathbf{s}(x) )}$ are $C^*$-algebras. It follows from $iii)$, $iv)$ and $v)$ that the restricted map $\beta:\a_{\mathbf{s}(x)} \to \a_{\beta_0(\mathbf{s}(x))}$ is a $C^*$-isomorphism and is therefore isometric. Hence we have
\begin{eqnarray*}
 \|\beta(a) \| =  \|\beta(a^*a)\|^{\frac{1}{2}} = \|a^*a\|^{\frac{1}{2}} = \|a\|\,,
\end{eqnarray*}
which shows that $\beta: \a_x \to \a_{\beta_0(x)}$ is an isometry.\\
\end{rem}

Given a Fell bundle $\a$ over a discrete groupoid $X$ we will denote by $C_c(\a)$ its corresponding $^*$-algebra of finitely supported sections. The following notation will be used throughout the rest of this work: for $x \in X$ and $a \in \a_x$ the symbol $a_x$ will always denote the element of $C_c(\a)$ such that
\begin{align}
 a_x (y) & := \begin{cases}
 a\,, \qquad\quad \text{if}\;\; y = x \\
 0, \qquad\quad \text{otherwise}\,.
\end{cases}
\end{align}\\

According to the notation above we can then write any $f \in C_c(\a)$ uniquely as a sum of the form
\begin{align}
\label{decomposition of f of CcA as a sum of ax}
 f = \sum_{x \in X} (f(x))_x\,.
\end{align}

We recall that the operations of multiplication and involution in $C_c(\a)$ are determined by
\begin{align*}
 a_x \cdot b_y & = \begin{cases}
 (ab)_{xy}\,, \qquad\quad \text{if}\;\; \mathbf{s}(x) = \mathbf{r}(y) \\
 0, \qquad\qquad\quad\; \text{otherwise}\,,
\end{cases}\\
 (a_x)^* & = (a^*)_{x^{-1}}\,,
\end{align*}
where $x, y \in X$ and $a \in \a_x$, $b \in \a_y$.\\

\section{Orbit space groupoids and Fell bundles}
\label{groupoids Fell bundles chapter}

In this section we present the basic set up which will enable us to define crossed products by Hecke pairs later in Section \ref{algebraic crossed products chapter}.

Our construction of a ($^*$-algebraic) crossed product $A \times^{alg} G / \gm$ of an algebra $A$ by a Hecke pair $(G, \gm)$ will make sense when $A$ is a certain algebra of sections of a Fell bundle over a discrete groupoid. In this section we show in detail what type of algebras $A$ are involved in the crossed product and how they are obtained.

\subsection{Group actions on Fell bundles}
\label{group actions on discrete groupoids section}

Throughout this section $G$ will denote a discrete group. One of our ingredients for defining crossed products by Hecke pairs consists of a group action on a Fell bundle over a groupoid (a concept we borrow from \cite[Section 6]{coacbundles}). Such actions always carry an associated action on the corresponding groupoid (by groupoid automorphisms). Since we are primarily interested in \emph{right} actions on groupoids, we start by recalling what they are:\\

\begin{df}
 Let $X$ be a groupoid.  A \emph{right action} of $G$ on $X$ is a mapping
\begin{align*}
 X \times G \to X\\
(x,g) \mapsto xg\,,
\end{align*}
which is a right action of $G$ on the underlying set of $X$, meaning that
\begin{itemize}
 \item[1)] $xe=x$, for all $x \in X$,
 \item[2)] $x(g_1g_2) = (xg_1)g_2$, for all $x \in X$, $g_1, g_2 \in G$,
\end{itemize}
which is compatible with the groupoid operations, meaning that
\begin{itemize}
 \item[3)] if $x$ and $y$ are composable  in $X$, then so are $xg$ and $yg$, for all $g \in G$, and moreover
\begin{align*}
 (xg)(yg) = (xy)g\,,
\end{align*}
 \item[4)] $(xg)^{-1}=x^{-1}g$, for all $x \in X$ and $g \in G$.
\end{itemize}
In other words, a right action of $G$ on $X$ is a right action on the set $X$ performed by groupoid automorphisms.\\
\end{df}

\begin{lemma}
\label{s(xg) = s(x)g}
 Let $X$ be a groupoid endowed with a right $G$-action. For every $x \in X$ and $g \in G$ we have
\begin{align*}
 \mathbf{s}(xg) = \mathbf{s}(x)g \qquad \text{and}\qquad \mathbf{r}(xg)=\mathbf{r}(x)g\,.
\end{align*}
In particular, $G$ restricts to an action on the unit space $X^0$.\\
\end{lemma}

{\bf \emph{Proof:}} It follows easily from the definition of a right $G$-action that
\begin{align*}
 \mathbf{s}(x)g = (x^{-1}x)g = (x^{-1}g)(xg) = (xg)^{-1}(xg) = \mathbf{s}(xg)\,,
\end{align*}
and similarly for the range function. \qed\\

\begin{rem}
 Given elements $x, y$ in a groupoid $X$ endowed with a right $G$-action and given $g \in G$, we will often drop the brackets in expressions like $(xg)y$ and simply use the notation $xgy$. No confusion arises from this since $G$ is only assumed to act on the right. On the other hand, we will never write an expression like $xyg$ without brackets, since it can be confusing on whether it means $x(yg)$ or $(xy)g$.\\
\end{rem}

\begin{df}{\cite[Section 6]{coacbundles}}
Let $G$ be a group and $\a$ a Fell bundle over a discrete groupoid $X$. An \emph{action} of $G$ on $\a$ consists of a homomorphism $\alpha:G \to Aut (\a)$.\\
\end{df}

As observed in Section \ref{Fell bundle section}, each automorphism of $\a$ carries with it an associated automorphism of the underlying groupoid $X$. Hence, an action of a group $G$ on $\a$ entails an action of $G$ on $X$ by groupoid automorphisms. Since we are interested only in right actions on groupoids, we just ensure that these associated actions are on the right simply by taking inverses. Moreover, even though we will typically denote by $\alpha$ the action of $G$ on $\a$, we will simply write $(x,g) \mapsto xg$ to denote its associated action on $X$ and it will be always assumed that this action comes from $\alpha$. To summarize what we have said so far: given an action $\alpha$ of $G$ on a Fell bundle $\a$ over a groupoid $X$, there is an associated right $G$-action $(x, g) \mapsto xg$ on $X$ such that
\begin{align}
 p(\alpha_g(a)) = p(a)g^{-1}\,.\\ \notag
\end{align}

\begin{rem}
Typically one would require the mapping $(a,g) \mapsto \alpha_g(a)$ to be continuous, but this is not necessary here since both $G$ and $X$ are discrete.\\
\end{rem}

\begin{prop}
\label{action on CcA}
 Let $\alpha$ be an action of a group $G$ on a Fell bundle $\a$ over a groupoid $X$. We have an associated action $\overline{\alpha}: G \to Aut(C_c(\a))$ of $G$ on $C_c(\a)$ given by
\begin{align*}
 \overline{\alpha}_g(f)\,(x) := \alpha_g(f(xg))\,,
\end{align*}
for $g \in G$, $f \in C_c(\a)$ and $x \in X$.\\
\end{prop}

{\bf \emph{Proof:}} Let us first prove that the action is well-defined, i.e. $\overline{\alpha}_g(f) \in C_c(\a)$. The fact that $\overline{\alpha}_g(f)$ is finitely supported is obvious, so the only thing one needs to check is that $\overline{\alpha}_g(f)$ is indeed a section of the bundle, i.e. $\alpha_g(f(xg)) \in \a_x$ for all $x \in X$, which is clear because $\alpha_g(\a_y) = \a_{yg^{-1}}$.

Let us now check that $\overline{\alpha}_g$ is indeed a $^*$-homomorphism for all $g \in G$. Linearity of $\overline{\alpha}_g$ is obvious. Let $f,f_1, f_2 \in C_c(\a)$. We have
\begin{eqnarray*}
 \overline{\alpha}_g (f_1 \cdot f_2)\,(x) & = & \alpha_g((f_1 \cdot f_2)\,(xg))\\  & = & \sum_{\substack{y,z \in X \\ yz=xg}} \alpha_g(f_1(y)f_2(z))\\
& = & \!\!\!\!\!\!\!\!\!\!\!\!\! \sum_{\substack{y,z \in X \\ \quad (yg^{-1})(zg^{-1})=x}} \!\!\!\!\!\!\!\! \alpha_g (f_1(y)) \alpha_g (f_2(z))\\ & = & \sum_{\substack{y,z \in X \\ yz=x}} \alpha_g (f_1(yg))\alpha_g (f_2(zg))\\
& = & \sum_{\substack{y,z \in X \\ yz=x}} \overline{\alpha}_g(f_1)(y) \overline{\alpha}_g(f_2)(z)\\ & = & \big(\overline{\alpha}_g(f_1) \cdot \overline{\alpha}_g(f_2)\big)\, (x)\,.
\end{eqnarray*}
Hence, $\overline{\alpha}_g(f_1 \cdot f_2) = \overline{\alpha}_g(f_1) \cdot \overline{\alpha}_g(f_2)$. Also,
\begin{eqnarray*}
 \overline{\alpha}_g (f^*) \,(x) & = & \alpha_g(f^* (xg)) = \alpha_g(f (x^{-1}g))^*\\
& = & \big( \overline{\alpha}_g(f)\,(x^{-1}) \big)^* = \big( \overline{\alpha}_g(f) \big)^* \, (x)\,.
\end{eqnarray*}
Hence, $\overline{\alpha}_g(f^*) = (\overline{\alpha}_g(f))^*$. The fact that $\overline{\alpha}_{g_1g_2} = \overline{\alpha}_{g_1} \circ \overline{\alpha}_{g_2}$ for every $g_1, g_2 \in G$ is also easily checked. \qed\\

\begin{df}
Let $\alpha$ be a group action of $G$ on a Fell bundle $\a$ over a groupoid $X$ and let $H$ be a subgroup of $G$. We will say that the $G$-action is \emph{$H$-good} if for any $x \in X$ and $h \in H$ we have
\begin{align}
 \mathbf{s}(x) h= \mathbf{s}(x) \; \Longrightarrow \; \alpha_{h^{-1}}(a) = a \qquad \forall a \in \a_x\,.
\end{align}
Also, a right $G$-action on a groupoid $X$ is said to be \emph{$H$-good} if for any $x \in X$ and $h \in H$ we have
\begin{align}
 \mathbf{s}(x) h= \mathbf{s}(x) \; \Longrightarrow \; xh = x\,.\\ \notag
\end{align}
\end{df}

It is clear from the definitions that if the action $\alpha$ of $G$ on $\a$ is $H$-good, then its associated right $G$-action on the underlying groupoid $X$ is also $H$-good. We will mostly use actions on Fell bundles. However, some of our results (namely Proposition \ref{X / H is a groupoid}) are about groupoids only, and this is the reason for defining $H$-good actions for groupoids as well.

We now give equivalent definitions of a $H$-good action. For that we recall from (\ref{stabilizer notation}) that given an action of $G$ on a set $X$ we denote by $\mathcal{S}_x$ the stabilizer of the point $x \in X$. We will also denote by $\mathcal{S}(\a_x)$ the set $\mathcal{S}(\a_x) := \{g \in G: \alpha_{g^{-1}}(a) = a \,, \forall a \in \a_x\}$.  \\

\begin{prop}
\label{equivalent def of H good action}
 Let $\alpha$ be an action of $G$ on a Fell bundle $\a$ over a groupoid $X$. The following statements are equivalent:
\begin{itemize}
 \item[i)] The action $\alpha$ is $H$-good.
 \item[ii)] For every $x \in X$ we have that $\mathcal{S}_{\mathbf{s}(x)} \cap H = \mathcal{S}(\a_x) \cap H$.
 \item[iii)] For any $x \in X$ we have
\begin{eqnarray}
\label{6 equality in good action}
 & & \mathcal{S}_{\mathbf{s}(x)} \cap H \;\; = \;\; \mathcal{S}_x \cap H \;\; = \;\; \mathcal{S}_{\mathbf{r}(x)} \cap H \;\; = \;\;\\
& = & \mathcal{S}(\a_{\mathbf{s}(x)}) \cap H \;\; = \;\; \mathcal{S}(\a_x) \cap H \;\; = \;\; \mathcal{S}(\a_{\mathbf{r}(x)}) \cap H\,. \notag
\end{eqnarray}

 \item[iv)] The stabilizers of the $H$-actions on $X$ and on the fibers of $\a$ are the same on composable pairs, i.e. if $x \in X$ and $y \in Y$ are composable, then
\begin{eqnarray*}
 & & \mathcal{S}_x \cap H = \mathcal{S}_y \cap H \;\; = \;\;\\
& = & \mathcal{S}(\a_x) \cap H \;\; = \;\; \mathcal{S}(\a_y) \cap H\,.
\end{eqnarray*}\\
\end{itemize}

\end{prop}

{\bf \emph{Proof:}} $i) \Longrightarrow ii)$ Since the action is $H$-good we have, by definition, that $\mathcal{S}_{\mathbf{s}(x)} \cap H \subseteq \mathcal{S}(\a_x) \cap H$. Also, if $h \in \mathcal{S}(\a_x) \cap H$, then we necessarily have $xh = x$, and therefore by Lemma \ref{s(xg) = s(x)g} we get $\mathbf{s}(x) = \mathbf{s}(xh) = \mathbf{s}(x)h$, from which we conclude that $h \in \mathcal{S}_{\mathbf{s}(x)} \cap H$. Hence we have $\mathcal{S}_{\mathbf{s}(x)} \cap H = \mathcal{S}(\a_x) \cap H$.

$ii) \Longrightarrow iii)$ Repeating a little bit of what we did above: if $h \in \mathcal{S}(\a_x) \cap H$, then we necessarily have that $xh = h$, and therefore $h \in \mathcal{S}_x \cap H$. Moreover, if $h \in \mathcal{S}_x \cap H$, then it follows that by Lemma \ref{s(xg) = s(x)g} that $h \in \mathcal{S}_{\mathbf{s}(x)} \cap H$. Thus, we have that
\begin{align*}
 \mathcal{S}(\a_x) \cap H = \mathcal{S}_x \cap H = \mathcal{S}_{\mathbf{s}(x)} \cap H\,.
\end{align*}

Since $\mathbf{s}(\mathbf{s}(x)) = \mathbf{s}(x)$, we also have, directly by our assumption of $ii)$, that $\mathcal{S}_{\mathbf{s}(x)} \cap H = \mathcal{S}(\a_{\mathbf{s}(x)}) \cap H$.

Since we have $(xg)^{-1} = x^{-1}g$, it follows easily that $\mathcal{S}_x = \mathcal{S}_{x^{-1}}$. Similarly, since $\alpha_{g}(a)^* = \alpha_g(a^*)$, it follows easily that $\mathcal{S}(\a_x) = \mathcal{S}(\a_{x^{-1}})$. Observing that $\mathbf{s}(x^{-1}) = \mathbf{r}(x)$, equality (\ref{6 equality in good action}) follows directly from what we proved above.

$iii) \Longrightarrow iv)$ Suppose $x \in X$ and $y \in X$ are composable. Then, $\mathbf{s}(x) = \mathbf{r}(y)$ and equality (\ref{6 equality in good action}) immediately yelds that
\begin{eqnarray*}
 & & \mathcal{S}_x \cap H = \mathcal{S}_y \cap H \;\; = \;\;\\
& = & \mathcal{S}(\a_x) \cap H \;\; = \;\; \mathcal{S}(\a_y) \cap H\,.
\end{eqnarray*}
$iv) \Longrightarrow i)$ Let $h \in H$ and $x \in X$ be such that $\mathbf{s}(x) h = \mathbf{s}(x)$. From $iv)$ it follows that $h \in \mathcal{S}(\a_x) \cap H$. This means that the action is $H$-good. \qed\\

It is easy to see that any $H$-good action is also $gHg^{-1}$-good for any conjugate $gHg^{-1}$, and also $K$-good for any subgroup $K \subseteq H$.

The following property will also be important for defining crossed products by Hecke pairs:\\

\begin{df}
 Let $X$ be a groupoid endowed with a right $G$-action and let $H$ be a subgroup of $G$. We will say that the action has the \emph{$H$-intersection property} if
\begin{align}
 uH \cap ugHg^{-1} = uH^g\,,
\end{align}
for every unit $u \in X^0$ and $g \in G$.

An action of $G$ on a Fell bundle $\a$ is said to have the \emph{$H$-intersection property} if its associated right $G$-action on the underlying groupoid has the $H$-intersection property.\\
\end{df}

We defer examples of $H$-good actions and actions with the $H$-intersection property for the next section. We now introduce one of the important ingredients for our definition of crossed products by Hecke pairs: the orbit space groupoid.

Let $G$ be a group, $H \subseteq G$ a subgroup and $X$ a groupoid endowed with a $H$-good right $G$-action. Then, the orbit space $X / H$ becomes a groupoid in a canonical way which we will now describe. For that, and throughout this text, we will use the following notation: given elements $x, y$ we define the set
\begin{align}
 H_{x,y}\index{Hxy@$H_{x,y}$} := \{h \in H : \mathbf{s}(x)h = \mathbf{r}(y)\}\,.
\end{align}
The groupoid structure on $X / H$ is described as follows:
\begin{itemize}
\item A pair $(xH, yH) \in (X / H)^2$ is composable if and only if $H_{x,y} \neq \emptyset$, or equivalently, $\mathbf{r}(y) \in \mathbf{s}(x) H$. This property is easily seen not to depend on the choice of representatives $x,y$ from the orbits $xH, yH$ respectively.
\item Given a composable pair $(xH, yH) \in (X / H)^2$, their product is
\begin{align}
  xH \,yH := x \widetilde{h} y H\,,
\end{align}
where $\widetilde{h}$ is any element of $H_{x,y}$. It will follow from the fact the action is $H$-good that $x\widetilde{h}$ does not depend on the representative $\widetilde{h}$ chosen from $H_{x,y}$. The result of the product $xH\, yH$ also does not depend on the choice of representatives $x, y$. We will prove this in the next result.
\item The inverse of the element $xH$ is simply the element $x^{-1}H$. It is also easy to see that this does not depend on the choice of representative $x$.\\
\end{itemize}

\begin{prop}
\label{X / H is a groupoid}
 Let $G$ be a group, $H \subseteq G$ a subgroup and $X$ a groupoid endowed with a $H$-good right $G$-action. The operations above give the orbit space $X / H$ the structure of a groupoid. Moreover, the unit space $(X / H)^0$ of this groupoid is $X^0 / H = \{uH: u \in X^0 \}$, where $X^0$ is the unit space of $X$, and the range and source functions satisfy
\begin{align*}
 \mathbf{s}(xH)=\mathbf{s}(x)H \qquad\qquad \text{and}\qquad\qquad \mathbf{r}(xH) = \mathbf{r}(x)H\,.\\
\end{align*}

\end{prop}

{\bf \emph{Proof:}} Let us first prove that the product is well-defined. Let $(xH, yH) \in (X / H)^2$ be a composable pair. The fact that $x\widetilde{h}$ does not depend on the representative $\widetilde{h}$ chosen from $H_{x,y}$ follows from the assumption that the action is $H$-good, since if $h_1, h_2 \in H_{x,y}$ then we have
\begin{align*}
 \mathbf{s}(x)h_1 = \mathbf{r}(y) = \mathbf{s}(x)h_2\,,
\end{align*}
and therefore $\mathbf{s}(x)h_1h_2^{-1} = \mathbf{s}(x)$, and because the action is $H$-good $xh_1h_2^{-1}=x$, i.e. $xh_1 = xh_2$.

Let us now prove that $X / H$ is a groupoid with the operations above. We check associativity first. Suppose $xH, yH, zH \in X / H$ are such that $(xH, yH)$ is composable and $(yH, zH)$ is composable. We want to prove that $(xHyH, zH)$ and $(xH, yHzH)$ are also composable and moreover $(xHyH)zH = xH(yHzH)$. We have by definition that $xHyH = x\widetilde{h_1}yH$ and $yHzH = y\widetilde{h_2}zH$, where $\widetilde{h_1}$ is any element of $H_{x,y}$ and $\widetilde{h_2}$ is any element of $H_{y,z}$. We now notice that
\begin{align*}
 H_{x\widetilde{h_1}y, z} = \{h \in H: \mathbf{s}(x\widetilde{h_1}y)h = \mathbf{r}(z)\} = \{h \in H: \mathbf{s}(y)h = \mathbf{r}(z)\} = H_{y, z}\,.
\end{align*}
Since $H_{y,z} \neq \emptyset$ it follows that $H_{x\widetilde{h_1}y, z} \neq \emptyset$, and therefore $(x\widetilde{h_1}yH, zH)$ is composable. Similarly,
\begin{eqnarray*}
 H_{x, y\widetilde{h_2}z} & = & \{h \in H: \mathbf{s}(x)h = \mathbf{r}(y\widetilde{h_2}z)\}\\
& = & \{h \in H: \mathbf{s}(x)h = \mathbf{r}(y)\widetilde{h_2}\}\\
& = & \{h \in H: \mathbf{s}(x)h\widetilde{h_2}^{-1} = \mathbf{r}(y)\}\\
& = & H_{x, y}\, \widetilde{h_2}\,.
\end{eqnarray*}
Hence, since $H_{x,y} \neq \emptyset$ it follows that $H_{x,y\widetilde{h_2}z} \neq \emptyset$, and therefore $(xH, y\widetilde{h_2}zH)$ is composable.

As we saw above $H_{x\widetilde{h_1}y, z} = H_{y, z}$, and since $\widetilde{h_2} \in H_{y, z}$, we can write
\begin{eqnarray*}
 (xHyH)zH & = & x\widetilde{h_1}yH zH \;\; = \;\; (x\widetilde{h_1}y)\widetilde{h_2}zH\\
& = & x\widetilde{h_1}\widetilde{h_2}y\widetilde{h_2}zH\,.
\end{eqnarray*}
Also seen above, we have that $H_{x, y\widetilde{h_2}z} = H_{x, y}\, \widetilde{h_2}$, so that $\widetilde{h_1}\widetilde{h_2} \in H_{x, y\widetilde{h_2}z}$. Hence, we conclude that
\begin{eqnarray*}
 (xHyH)zH & = & xH(yHzH)\,.
\end{eqnarray*}
We now check that for any element $xH \in X / H$ we have that $(xH, x^{-1}H)$ and $(x^{-1}H, xH)$ are composable pairs. We have that
\begin{eqnarray*}
 H_{x, x^{-1}} & = & \{h \in H: \mathbf{s}(x)h = \mathbf{r}(x^{-1})\} \;\; =\;\; \{h \in H: \mathbf{s}(x)h = \mathbf{s}(x)\}\,,
\end{eqnarray*}
and the identity element $e$ obviously belongs to the latter set. Hence we conclude that $H_{x, x^{-1}} \neq \emptyset$, and therefore $(xH, x^{-1}H)$ is composable. A similar observation shows that $(x^{-1}H, xH)$ is also composable.

To prove that $X / H$ is a groupoid it now remains to prove the inverse identities $xHyHy^{-1}H = xH$ and $y^{-1}HyHxH = xH$, in case $(xH, yH)$ is composable (for the first identity) and $(yH, xH)$ is composable (for the second identity). We first show that $yHy^{-1}H = \mathbf{r}(y)H$. We have that $yHy^{-1}H = y\widetilde{h}y^{-1}H$ for any element $\widetilde{h} \in H_{y, y^{-1}}$. Since, as we observed above, we always have $e \in H_{y, y^{-1}}$, it follows that we can take $\widetilde{h}$ as $e$. Thus, we get
\begin{eqnarray}
\label{r(yH)=r(y)H}
 yHy^{-1}H & = & yy^{-1}H \;\; =\;\; \mathbf{r}(y)H\,.
\end{eqnarray}
From this it follows that
\begin{eqnarray*}
 xHyHy^{-1}H & = & xH\mathbf{r}(y)H \;\; =\;\; x\widetilde{h_1}\mathbf{r}(y)H\,,
\end{eqnarray*}
where $\widetilde{h_1}$ is any element of $H_{x, \mathbf{r}(y)}$. By definition, $\widetilde{h_1}$ is such that $\mathbf{r}(y) = \mathbf{s}(x) \widetilde{h_1} = \mathbf{s}(x \widetilde{h_1})$. Hence we have that $x\widetilde{h_1}\mathbf{r}(y) = x\widetilde{h_1}$, and therefore
\begin{align*}
 xHyHy^{-1}H = x\widetilde{h_1}H = xH\,.
\end{align*}
The other identity $y^{-1}HyHxH = xH$ is proven in a similar fashion. Hence, we conclude that $X / H$ is a groupoid.

From equality (\ref{r(yH)=r(y)H}) it follows easily that the units of $X / H$ are precisely the elements of the form $uH$ where $u \in X^0$, so that we can write $(X / H)^0 = X^0 / H$. Also from (\ref{r(yH)=r(y)H}) it follows that the range function in $X / H$ satisfies:
\begin{align*}
 \mathbf{r}(xH) = \mathbf{r}(x)H\,.
\end{align*}
The analogous result for the source function is proven in a similar fashion. \qed\\

Let $\alpha$ be an action of $G$ on a Fell bundle $\a$ over a groupoid $X$. Assume that the action is $H$-good, where $H$ is a subgroup of $G$. We will now define a new Fell bundle $\a / H$\index{A/H@$\a / H$} over the groupoid $X / H$. First we set some notation. The set of $H$-orbits of the action $\alpha$ on $\a$ gives us a partition of $\a$ into equivalence classes. We will denote by $[a]$ the equivalence class of the element $a \in \a$, i.e.
\begin{align*}
[a] := \{\alpha_h(a)\}_{h \in H}\,.
\end{align*}
We define $\a / H$ as the set of all the $H$-orbits in $\a$, i.e.
\begin{align}
\label{orbit Fell bundle A H}
\a/ H := \{ [a] : a \in \a \}\,,
\end{align}
which, as we will now see, is a Fell bundle over $X / H$ in a natural way.\\

\begin{prop}
\label{orbit bundle is a Fell bundle}
 Let $\alpha$ be an action of a group $G$ on a Fell bundle $\a$ over a groupoid $X$ and $H \subseteq G$ be a subgroup for which the $G$-action is $H$-good. The set of $H$-orbits $\a / H$ forms a Fell bundle over the groupoid $X / H$ in the following way:
\begin{itemize}
 \item The associated projection $p_H: \a / H \to X / H$ is defined by $p_H([a]) := p(a)H$, where $p$ is the associated projection of the bundle $\a$.
 \item The vector space structure on each fiber $\big(\a / H \big)_{xH}$ is defined in the following way: if $a, b \in \a_x$ then $[a] + [b] := [a + b]$, and if $\lambda \in \mathbb{C}$ then $\lambda [a] := [\lambda a]$.
 \item The norm on $\a / H$ is defined by $\| [a] \| := \| a \|$.
 \item The multiplication maps $\big(\a / H \big)_{xH} \times \big(\a / H \big)_{yH} \to \big(\a / H \big)_{xH\cdot yH}$, for a composable pair $(xH, yH)$, are defined in the following way: if $a \in \a_x$ and $b \in \a_y$, then
\begin{align}
 [a][b] = [\alpha_{\widetilde{h}^{-1}}(a)b]\,,
\end{align}
where $\widetilde{h}$ is any element of $H_{x,y}$.
\item The involution map is defined by $[a]^*:= [a^*]$.\\
\end{itemize}
\end{prop}

\begin{lemma}
\label{representative of [a] is unique}
  Let $\alpha$ be an action of $G$ on a Fell bundle $\a$ over a groupoid $X$ and $H \subseteq G$ be a subgroup for which the $G$-action is $H$-good. Let $x \in X$ and $a \in \a_x$. Given any $y \in xH$ there exists a unique representative $b$ of $[a]$ such that $b \in \a_y$.\\
\end{lemma}

{\bf \emph{Proof:}} Given an element $y \in xH$ we have that $y = x h$ for some $h \in H$. The element $\alpha_{h^{-1}}(a)$ is then a representative of $[a]$ such that $\alpha_{h^{-1}}(a) \in \a_{xh} = \a_y$, thus existence is established.

The uniqueness claim follows from the fact the action is $H$-good. Suppose we have two representatives $b$ and $c$ of $[a]$ such that both $b$ and $c$ belong to $\a_y$. Being representatives of $[a]$ means that there are elements $h_1, h_2 \in H$ such that $b = \alpha_{h_1}(a)$ and $c = \alpha_{h_2}(a)$. Hence we have that
\begin{align*}
 \alpha_{h_2h_1^{-1}}(b) = c\,,
\end{align*}
and therefore $h_2h_1^{-1}$ takes $\a_y$ into $\a_y$. This means that $y h_1h_2^{-1} = y$ and therefore $\mathbf{s}(y) h_1h_2^{-1} = \mathbf{s}(y)$. Since the action is $H$-good it follows that $\alpha_{h_2h_1^{-1}}(b) = b$, and therefore $b = c$. \qed\\

{\bf \emph{Proof of Proposition \ref{orbit bundle is a Fell bundle}:}} First, it is clear that the vector space structure on each fiber $\big(\a / H \big)_{xH}$ is well-defined. By this we mean two things: first, given two elements $[a], [b] \in \big(\a / H \big)_{xH}$ there exist unique representatives $a,b$ such that $a,b \in \a_x$ for a given representative  $x$ of the orbit $xH$ (Lemma \ref{representative of [a] is unique}); second, the sum $[a + b]$ still lies in $\big(\a / H \big)_{xH}$ and does not depend on the choice of representatives $a$  and $b$ (provided only that $a$ and $b$ are in the same fiber).

The norm on $\a / H$ is also easily seen to be well-defined, i.e. independent of the choice of representative. This is true because any other representative of $[a]$ is of the form $\alpha_h(a)$ for some $h \in H$, and by Remark \ref{automorphism gives isometry between fibers} we know that $\alpha_h$ gives an isometry between fibers. It is also clear that each fiber $\big(\a / H \big)_{xH}$ is a Banach space under this norm.

The multiplication map is also easily seen to be well-defined: using the fact the $G$-action on $\a$ is $H$-good we know that $\alpha_{\widetilde{h}^{-1}}(a)b$ does not depend on the choice of element $\widetilde{h} \in H_{x,y}$. Moreover, $\alpha_{\widetilde{h}^{-1}}(a)b \in \a_{x\widetilde{h}y}$ and therefore $[\alpha_{\widetilde{h}^{-1}}(a)b] \in \a_{x\widetilde{h}yH}$. The fact that the multiplication map does not depend on the chosen representatives of the orbits $[a]$ and $[b]$ is also easily checked.

It follows from a routine computation that map $\big(\a / H \big)_{xH} \times \big(\a / H \big)_{yH} \to \big(\a / H \big)_{xH\cdot yH}$ is bilinear. Moreover, for $[a] \in \big(\a / H \big)_{xH}$ and $[b] \in \big(\a / H \big)_{yH}$, where we assume without loss of generality that $a \in \a_x$ and $b \in \a_y$, we have that
\begin{eqnarray*}
 \| [a][b]\| & = & \| \alpha_{\widetilde{h}^{-1}}(a)b \| \;\; \leq \;\; \| \alpha_{\widetilde{h}^{-1}}(a) \| \|b \|\\
& = & \| a \| \|b \| \;\; = \;\; \|[a] \| \|[b]\|\,.
\end{eqnarray*}
We will now check associativity of the multiplication maps. Let $(xH, yH)$ and $(yH, zH)$ be two composable pairs in $X / H$, and let $[a] \in (\a / H )_{xH}$, $[b] \in (\a / H )_{yH}$ and $[c] \in (\a / H )_{zH}$, where we assume without loss of generality that $a \in \a_x$, $b \in \a_y$ and $c \in \a_z$. By definition, we have $[a][b] = [\alpha_{\widetilde{h_1}^{-1}}(a)b]$, where $\widetilde{h_1}$ is any element of $H_{x,y}$. Thus, we have
\begin{eqnarray*}
 \big([a][b] \big)[c] & = & [\alpha_{\widetilde{h_1}^{-1}}(a)b][c] \;\; = \;\; [\alpha_{\widetilde{h_2}^{-1}}\big( \alpha_{\widetilde{h_1}^{-1}}(a)b \big) c]\\
 & = & [\alpha_{\widetilde{h_2}^{-1}\widetilde{h_1}^{-1} } (a) \alpha_{\widetilde{h_2}^{-1}}(b) c]\,,
\end{eqnarray*}
where $\widetilde{h_2}$ is any element of $H_{x\widetilde{h_1}y, z}$. One can  easily check (or see the proof of Proposition \ref{X / H is a groupoid} where this is done) that $H_{x\widetilde{h_1}y, z} = H_{y,z}$ and moreover that $\widetilde{h_1}\widetilde{h_2} \in H_{x,y\widetilde{h_2}z}$. From this observations it follows that
\begin{eqnarray*}
 \big([a][b] \big)[c] & = & [\alpha_{\widetilde{h_2}^{-1}\widetilde{h_1}^{-1} } (a) \alpha_{\widetilde{h_2}^{-1}}(b) c] \;\; = \;\; [a] [ \alpha_{\widetilde{h_2}^{-1}}(b) c]\\
& = & [a] \big([b][c])\,.
\end{eqnarray*}
Hence, the multiplication maps are associative.

The involution on $\a / H$ is also easily seen not to depend on choice of representative of the orbit, since the maps $\alpha_h$ preserve the involution of $\a$. Morevoer, it is easily checked that: if $[a] \in \big(\a / H \big)_{xH}$ then $[a]^* \in \big(\a / H \big)_{x^{-1}H}$, the associated map $\big(\a / H \big)_{xH} \to \big(\a / H \big)_{x^{-1}H}$ is conjugate linear, and $[a]^{**} = [a]$. Let us now check that $\big([a][b] \big)^* =[b]^*[a]^*$, whenever the multiplication is defined. Let us assume that $a \in \a_x$ and $b \in \a_y$ and that $(xH, yH)$ is composable. We have that
\begin{eqnarray*}
 \big([a][b] \big)^* & = & [\alpha_{\widetilde{h}^{-1}}(a)b]^* \;\; = \;\; [b^*\alpha_{\widetilde{h}^{-1}}(a^*)] \;\; = \;\; [\alpha_{\widetilde{h}}(b^*)a^*]\,,
\end{eqnarray*}
where $\widetilde{h}$ is any element of $H_{x,y}$. It is easily seen that $\widetilde{h}^{-1} \in H_{y^{-1}, x^{-1}}$, so that
\begin{eqnarray*}
 \big([a][b] \big)^* & = & [\alpha_{\widetilde{h}}(b^*)a^*] \;\; = \;\; [b^*][a^*] \;\; = \;\; [b]^*[a]^*\,.
\end{eqnarray*}
We also need to prove that $\|[a]^*[a] \| = \|[a]\|^2$. This is also easy because
\begin{align*}
 \|[a]^*[a] \| = \|[a^*][a] \| = \|[a^*a] \| = \|a^*a\| = \|a\|^2 = \|[a]\|^2\,.
\end{align*}
The last thing we need to check is that if $[a] \in (\a / H)_{xH}$, then $[a]^*[a]$ is a positive element of $(\a / H)_{\mathbf{s}(x)H}$ (seen as a $C^*$-algebra). We have that $[a]^*[a] = [a^*a]$.  We can assume without loss of generality that $a \in \a_x$, so that $a^*a \in \a_{\mathbf{s}(x)}$. Since $\a$ is a Fell bundle we have that $a^*a$ is a positive element of $\a_{\mathbf{s}(x)}$ (seen as a $C^*$-algebra). Hence, there exists an element $b \in \a_{\mathbf{s}(x)}$ such that $a^*a = b^*b$. Moreover, $[b] \in (\a / H)_{\mathbf{s}(x)}$ and it is now clear that
\begin{align*}
 [a]^*[a] = [a^*a] = [b^*b] = [b]^*[b]\,,
\end{align*}
i.e.  $[a]^*[a]$ is a positive element of $(\a / H)_{\mathbf{s}(x)}$. This finishes our proof that $\a / H$ is a Fell bundle. \qed\\

\begin{convention}
 For simplicity we will henceforward make the following convention. Given an orbit Fell bundle $\a / H$ as discribed in Proposition \ref{orbit bundle is a Fell bundle}, if we write that an element $[a]$ belongs to some fiber $(\a / H)_{xH}$, we will always assume that the representative $a$ of $[a]$ belongs to fiber over the representative $x$ of $xH$. In other words, if we write that $[a] \in (\a / H)_{xH}$, then we are implicitly assuming that $a \in \a_x$. This is possible and unambiguous by Lemma \ref{representative of [a] is unique}.

We apply this convention also for elements of $C_c(\a / H)$, meaning that a canonical element $[a]_{xH} \in C_c(\a / H)$ is always assumed to be written in a way that $a \in \a_x$.\\
\end{convention}

It is a strighforward fact that any function in $C_c(X / H)$ can also be seen as a complex-valued ($H$-invariant) function in $X$. This funcion in $X$ is in general no longer finitely supported, but it still makes sense as a function in $C(X)$, the vector space of all complex-valued functions in $X$. We will now see that something analogous can be said for the elements of $C_c(\a / H)$.

Given an element $f \in C_c(\a / H)$ we define a function $\iota(f) \in C(\a)$, where $C(\a)$ is the vector space of all sections of $\a$, by the following rule:
\begin{align}
 \iota(f)(x) := R_x(f(xH))\,,
\end{align}
where $R_x(f(xH))$ is the unique representative of $f(xH)$ such that $R_x(f(xH)) \in \a_x$, which is well-defined according to Lemma \ref{representative of [a] is unique}. It is then easy to see that the map $\iota$ is an injective linear map from $C_c(\a / H)$ to $C(\a)$.

For ease of reading we will henceforward drop the symbol $\iota$ and use the same notation both for elements of $C_c(\a / H)$ and for their correspondents in $C(\a)$. It will then be clear from context which one we are using.

Under this convention we can then write, for any $f \in C_c(\a / H)$ and $x \in X$, that $[f(x)] = f(xH)$. Moreover, the decomposition (\ref{decomposition of f of CcA as a sum of ax}) of $f \in C_c(\a / H)$ as a sum of elements of the form $[a]_{xH}$ can now be written as:
\begin{align}
\label{f written as a sum of [f(x)]}
 f = \sum_{xH \in X / H} \big(f(xH) \big)_{xH} = \sum_{xH \in X / H} [f(x)]_{xH}\,.\\ \notag
\end{align}

\subsection{Examples}

In this section we give some examples of $H$-good actions and actions satisfying the $H$-intersection property.
For the rest of the section we assume that $\a$ is a Fell bundle over a groupoid $X$ where a group $G$ acts and $H \subseteq G$ denotes a subgroup.

The first two examples (\ref{example action is free then it is good and int prop} and \ref{example action fixes every point then it is good and int prop}) show that $H$-good actions that satisfy the $H$-intersection property are present in actions that have completely opposite behaviours, such as free actions and actions that fix every point.\\

\begin{ex}
\label{example action is free then it is good and int prop}
 If the restricted action of $H$ on the unit space $X^0$ is free, then the action is $H$-good and satisfies the $H$-intersection property.\\
\end{ex}

\begin{ex}
\label{example action fixes every point then it is good and int prop}
 If the restricted action of $H$ on $\a$ fixes every point, then the action is $H$-good and satisfies the $H$-intersection property.\\
\end{ex}

The following example is one of the examples that motivated the development of this theory of crossed products by Hecke pairs. This example, and the study of the crossed products associated to it, seems to be valuable for obtaining a form of Katayama duality with respect to crossed products by ``coactions'' of discrete homogeneous spaces.\\

\begin{ex}
\label{katayama duality example}
  Suppose $X$ is the transformation groupoid $G \times G$. We recall that the multiplication and inversion operations on this groupoid are given by:
\begin{align*}
  (s,tr)(t,r)=(st,r)\qquad \text{and}\qquad (s,t)^{-1}=(s^{-1}, st)\,.
\end{align*}
Recall also that the source and range functions on $G \times G$ are defined by
\begin{align*}
 \mathbf{s}(s,t) = (e, t) \qquad \text{and}\qquad \mathbf{r}(s,t) = (e, st)\,.
\end{align*}

We observe that there is a natural right $G$-action on $G \times G$, given by
\begin{align}
\label{G action on G times G}
 (s,t)g := (s, tg)\,.
\end{align}

Let $\delta$ be a coaction of $G$ on a $C^*$-algebra $B$ and $\b$ the associated Fell bundle. Following \cite[Section 3]{ind}, we will denote by $\a:=\b \times G$ the corresponding Fell bundle over the groupoid $G \times G$. Elements of $\a$ have the form $(b_s, t)$, where $b_s \in \b_s$ and $s,t \in G$. Any such element lies in the fiber $\a_{(s,t)}$ over $(s,t)$.

It is easy to see that there is a canonical actin $\alpha$ of $G$ on $\a$, given by
\begin{align*}
 \alpha_g(b_s, t) := (b_s, tg^{-1})\,.
\end{align*}
This action of $G$ on $\a$ entails the natural right action of $G$ on $G \times G$, as described in (\ref{G action on G times G}). This $G$-action on $G \times G$ is free and therefore the action $\alpha$ is $H$-good and satisfies the $H$-intersection property with respect to any subgroup $H\subseteq G$.

 The orbit space groupoid $(G \times G) / H$ can be canonically identified with the  groupoid $G \times G / H$ of \cite{full}, whose operations are given by:
\begin{align*}
  (s,trH)(t,rH)=(st,rH)\qquad \text{and}\qquad (s,tH)^{-1}=(s^{-1}, stH)\,.
\end{align*}
 Moreover, the orbit Fell bundle $\a / H $ is canonically identified with the Fell bundle $\b \times G / H$ over $G \times G / H$ defined in \cite{full}, and in this way $C_c(\a / H)$ is canonically isomorphic with the Echterhoff-Quigg algebra $C_c(\b \times G / H)$, also from \cite{full}.\\
\end{ex}

\subsection{The algebra $M(C_c(\a))$}
\label{big algrebra MCcA section}

We will assume for the rest of this section that $G$ is a group, $H \subseteq G$ is a subgroup and $\a$ is a Fell bundle over a groupoid $X$ endowed with a $G$-action $\alpha$. We also assume that the action $\alpha$ is $H$-good. We recall that $\a / H$ stands for the orbit Fell bundle over the groupoid $X / H$, as defined in (\ref{orbit Fell bundle A H}).

For the purpose of defining crossed products by Hecke pairs it is convenient to have a ``large'' algebra which contains the algebras $C_c(\a / H)$ for different subgroups $H\subseteq G$. In this way we are allowed to multiply elements of $C_c(\a / H)$ and $C_c(\a / K)$, for different subgroups $H, K \subseteq G$, in a meaningful way. This large algebra will be the multiplier algebra $M(C_c(\a))$. This section is thus devoted to show how algebras such as $C_c(\a / H)$ and $C_c(X^0 / H)$ embed in $M(C_c(\a))$ in a canonical way.

 Our first result shows that there is a natural inclusion $C_c(\a / H) \subseteq M(C_c(\a))$.\\

\begin{thm}
\label{embeding of Cc(a H) in M(Cc(a))}
There is an embedding $\iota$ of $C_c(\a /H)$ into $M(C_c(\a))$ determined by the following rule: for any $x, y \in X$, $a \in \a_x$ and $b \in \a_y$ we have
\begin{align}
\label{axH as a multiplier expression}
 \iota([a]_{xH}) b_y & := \begin{cases}
 (\alpha_{\widetilde{h}^{-1}}(a)b)_{x\widetilde{h}y}\,, \qquad\; \text{if}\;\; H_{x,y} \neq \emptyset \\
 0, \qquad\qquad\qquad\qquad \text{otherwise,}\;\;
\end{cases}
\end{align}
where $\widetilde{h}$ is any element of $H_{x, y}$.\\
\end{thm}

\begin{rem}
 The above result allows us to see $C_c(\a / H)$ as a $^*$-subalgebra of $M(C_c(\a))$. We shall henceforward drop the symbol $\iota$ and make no distinction of notation between an element of $C_c(\a / H)$ and its correspondent multiplier in $M(C_c(\a))$.\\
\end{rem}

{\bf \emph{Proof of Theorem \ref{embeding of Cc(a H) in M(Cc(a))}:}} Let us first show that expression (\ref{axH as a multiplier expression}) does indeed define an element of $M(C_c(\a))$. For this it is enough to check that $\langle \iota([a]_{xH}) b_y, c_z\rangle = \langle b_y, \iota([a]^*_{x^{-1}H})\, c_z \rangle$, for all $b \in \a_y$ and $c \in \a_z$, with $y, z \in X$. For $\iota([a]_{xH})b_y$ to be non-zero, we must necessarily have $H_{x,y} \neq \emptyset$, and in this case $\iota([a]_{xH})b_y = (\alpha_{\widetilde{h}^{-1}}(a)b)_{x\widetilde{h}y}$, where $\widetilde{h} \in H_{x,y}$. Now,
\begin{eqnarray*}
 \langle \iota([a]_{xH}) b_y, c_z\rangle & = & \langle (\alpha_{\widetilde{h}^{-1}}(a)b)_{x\widetilde{h}y} , c_z \rangle  = (b^*\alpha_{\widetilde{h}^{-1}}(a)^*)_{y^{-1}(x^{-1}\widetilde{h})} c_z\\
& = & b^*_{y^{-1}} \alpha_{\widetilde{h}^{-1}}(a)^*_{x^{-1}\widetilde{h}}\, c_z
\end{eqnarray*}
For $\alpha_{\widetilde{h}^{-1}}(a)^*_{x^{-1}\widetilde{h}}\, c_z$ to be non-zero we must necessarily have $\mathbf{r}(z) = \mathbf{s}(x^{-1})\widetilde{h}$, i.e. $\widetilde{h} \in H_{x^{-1}, z}$. So, to summarize, for $\langle [a]_{xH} b_y, c_z\rangle$ to be non-zero we must have $H_{x,y} \cap H_{ x^{-1}, z} \neq \emptyset$ and in this case we obtain
\begin{align*}
 \langle \iota([a]_{xH}) b_y, c_z\rangle = b^*_{y^{-1}}  \alpha_{\widetilde{h}^{-1}}(a)^*_{x^{-1}\widetilde{h}}\, c_z\,,
\end{align*}
where $\widetilde{h}$ is any element of $H_{x,y} \cap H_{ x^{-1}, z}$. A similar computation for $\langle b_y, \iota([a]^*_{x^{-1}H})\, c_z \rangle$ yelds the exact same result.

Recall from (\ref{decomposition of f of CcA as a sum of ax}) that any $f \in C_c(\a / H)$ can be written as
\begin{align*}
 f = \sum_{xH \in X / H} \big( f(xH) \big)_{xH}\,.
\end{align*}
From this we are able to define a multiplier $\iota(f) \in M(C_c(\a))$, simply by extending expression (\ref{axH as a multiplier expression}) by linearity.

 We want to show that $\iota$ is an injective $^*$-homomorphism. First, we claim that given $[a]_{xH}, [b]_{yH} \in C_c(\a / H)$ we have
\begin{align*}
 \iota([a]_{xH}) \iota([b]_{yH}) = \iota([a]_{xH} [b]_{yH})\,.
\end{align*}
 This amounts to proving that
\begin{align*}
 \iota([a]_{xH}) \iota( [b]_{yH}) & = \begin{cases}
 \iota([\alpha_{\widetilde{h}^{-1}}(a)b]_{x\widetilde{h}yH})\,, \qquad\;\; \text{if}\;\; H_{x,y} \neq \emptyset \\
 0, \qquad\qquad\qquad\qquad\qquad \text{otherwise}\;\; 
\end{cases}
\end{align*}
with $\widetilde{h}$ being any element of $H_{x,y}$. To see this, let $c_z \in \a_z$, with $z \in X$. We have
\begin{eqnarray*}
 \iota([a]_{xH})\iota( [b]_{yH}) c_z \!\!\!\! & = & \!\!\! \begin{cases}
 \iota([a]_{xH})(\alpha_{h_0^{-1}}(b)c)_{yh_0z}\,, \qquad\; \text{if}\;\; H_{y,z} \neq \emptyset \\
 0, \qquad\qquad\qquad\qquad\qquad\quad\;\; \text{otherwise}\;\;
\end{cases}\\
& = & \!\!\! \begin{cases}
 (\alpha_{h_1^{-1}}(a)\alpha_{h_0^{-1}}(b)c)_{xh_1yh_0z}\,, \; \text{if}\;\; H_{y,z} \neq \emptyset\;\, \text{and}\;\, H_{x, yh_0z} \neq \emptyset \\
 0, \qquad\qquad\qquad\qquad\qquad\quad \text{otherwise}\;\;
\end{cases}
\end{eqnarray*}
with $h_0 \in H_{y,z}$ and $h_1 \in H_{x, yh_0z}$. But $H_{x, yh_0z} = H_{x, y h_0} = H_{x, y} \, h_0$, hence the above can be written as
\begin{eqnarray*}
& = & \begin{cases}
 (\alpha_{h_0^{-1}\widetilde{h}^{-1}}(a)\alpha_{h_0^{-1}}(b)c)_{x\widetilde{h}h_0yh_0z}\,, \qquad\; \text{if}\;\; H_{y,z} \neq \emptyset\;\; \text{and}\;\; H_{x, y} \neq \emptyset \\
 0, \qquad\qquad\qquad\qquad\qquad\qquad\qquad\;\; \text{otherwise}\;\;
\end{cases}\\
& = & \begin{cases}
 \big(\alpha_{h_0^{-1}}(\alpha_{\widetilde{h}^{-1}}(a)b)c \big)_{(x\widetilde{h}y)h_0z}\,, \qquad\; \text{if}\;\; H_{y,z} \neq \emptyset\;\; \text{and}\;\; H_{x, y} \neq \emptyset \\
 0, \qquad\qquad\qquad\qquad\qquad\qquad\quad\; \text{otherwise}\;\;
\end{cases}
\end{eqnarray*}
where $\widetilde{h} \in H_{x,y}$. Also, $H_{y,z} = H_{x\widetilde{h}y, z}$. Thus, we obtain
\begin{eqnarray*}
& = & \begin{cases}
 \big(\alpha_{h_0^{-1}}(\alpha_{\widetilde{h}^{-1}}(a)b)c \big)_{(x\widetilde{h}y)h_0z}\,, \qquad\; \text{if}\;\; H_{x\widetilde{h}y,z} \neq \emptyset\;\; \text{and}\;\; H_{x, y} \neq \emptyset \\
 0, \qquad\qquad\qquad\qquad\qquad\qquad\quad \text{otherwise}\;\;
\end{cases}\\
& = & \begin{cases}
 \iota([\alpha_{\widetilde{h}^{-1}}(a)b]_{x\widetilde{h}yH})\, c_z\,, \qquad\;\;\; \text{if}\;\;  H_{x, y} \neq \emptyset \\
 0, \qquad\qquad\qquad\qquad\qquad\quad \text{otherwise}\;\;
\end{cases}
\end{eqnarray*}

Since $\iota$ is linear and multiplicative on the elements of the form $[a]_{xH}$, it is necessarily an homomorphism. 
Now the fact that $\iota([a]_{xH})^* = \iota(([a]_{xH})^*) = \iota( [a]^*_{x^{-1}H})$ follows directly from the computations in the beginning of this proof. Hence, $\iota$ is a $^*$-homomorphism.

Let us now prove injectivity of $\iota$. Suppose $f \in C_c(\a / H)$ is such that $\iota(f) = 0$. Decomposing $f$ as a sum of elements of the form $[a]_{xH}$, following (\ref{f written as a sum of [f(x)]}), we get
\begin{align*}
 0 = \iota(f) = \sum_{xH \in X/H} \iota \Big( \big( f(xH) \big)_{xH} \Big) = \sum_{xH \in X/H} \iota \big( [ f(x)]_{xH} \big) \,.
\end{align*}
For any $y \in X$ we then have
\begin{eqnarray*}
 0 & = & \sum_{xH \in X/H} \iota  \big( [f(x)]_{xH} \big) (f(y)^*)_{y^{-1}}\\
& = & \sum_{\substack{xH \in X/H \\ \mathbf{s}(y) \in \mathbf{s}(x)H}} \iota \big( [f(x)]_{xH} \big) (f(y)^*)_{y^{-1}}\\
& = & \sum_{\substack{xH \in X/H \\ \mathbf{s}(y) \in \mathbf{s}(x)H}}  \big( \alpha_{\widetilde{h_x}^{-1}}(f(x)) f(y)^* \big)_{x \widetilde{h_{x}}y^{-1}}\,,
\end{eqnarray*}
where $\widetilde{h_{x}}$ is any element of $H_{x, y^{-1}}$. Now the elements $x \widetilde{h_{x}} y^{-1}$ in the sum above are all different, because if we had $x_1\widetilde{h_{x_1}} y^{-1} = x_2 \widetilde{h_{x_2}} y^{-1}$, then we would have $x_1\widetilde{h_{x_1}} = x_2 \widetilde{h_{x_2}}$ and therefore $x_1 H = x_2 H$. Therefore each of the summands in the above sum is zero, and in particular we must have
\begin{eqnarray*}
 0 & = & \big( \alpha_{\widetilde{h_y}^{-1}}(f(y)) f(y)^* \big)_{y \widetilde{h_{y}}y^{-1}}\\
& = & \big( f(y) f(y)^* \big)_{\mathbf{r}(y)}\,,
\end{eqnarray*}
and therefore $f(y) f(y)^* = 0$. Hence we get $f(y) = 0$, and since this is true for any $y\in X$, we have $f = 0$, i.e. $\iota$ is injective. \qed\\

\begin{prop}
\label{Cb(X0) embedds in MCc(A)}
 There is an embedding $\iota$ of $C_b(X^0)$ into $M(C_c(\a))$ defined by
\begin{align}
\label{expression for f in Cb(X0) as a multiplier}
\iota(f) \, b_y := f(\mathbf{r}(y))b_y\,.
\end{align}
for every $f \in C_b(X^0)$, $y \in X$ and $b \in \a_y$.\\
\end{prop}

\begin{rem}
 The above result allows us to see $C_b(X^0)$ as a $^*$-subalgebra of $M(C_c(\a))$. We shall henceforward drop the symbol $\iota$ and make no distinction of notation between an element of $C_b(X^0)$ and its correspondent multiplier in $M(C_c(\a))$.\\
\end{rem}

{\bf \emph{Proof of Proposition \ref{Cb(X0) embedds in MCc(A)} :}} It is easy to see that $\langle \iota(f) b_y\,,\, c_z \rangle = \langle b_y\,,\, \iota(f^*)c_z \rangle$ for any $y,z \in X$, $b \in \a_y$ and $c \in \a_z$,  so that the expression (\ref{expression for f in Cb(X0) as a multiplier}) does define an element of $M(C_c(\a))$.

Hence we get a linear map $\iota: C_b(X^0) \to M(C_c(\a))$. Given two elements $f_1, f_2 \in C_b(X^0)$, we have that
\begin{align*}
 \iota(f_1) \iota(f_2) b_y = f_1(\mathbf{r}(y))f_2(\mathbf{r}(y)) b_y = \iota(f_1f_2)b_y
\end{align*}
for any $y \in X$ and$b \in \a_y$, so that $\iota$ is a $^*$-homomorphism. Hence, we only need to prove that $\iota$ is injective. This is not difficult to see: given $f \in C_b(X^0)$ such that $\iota(f) = 0$ we have, for any unit $u \in X^0$ and $b \in \a_u$, that
\begin{eqnarray*}
 0 & = & \iota(f) b_u\;\; = \;\; f(u)b_u\,.
\end{eqnarray*}
Hence, $f(u) = 0$ because each fiber $\a_u$ is non-zero by our assumption on Fell bundles (see Assumption \ref{assumption on Fell bundles}). Since this is true for any $u \in X^0$ we get $f = 0$, i.e. $\iota$ is injective. \qed\\

 Recall, from Lemma \ref{s(xg) = s(x)g}, that the action of $G$ on $X$ restricts to an action of $G$ on the set $X^0$. Thus it makes sense to talk about the commutative $^*$-algebra
\begin{align*}
 C_c(X^0 / H) \subseteq C_b(X^0)\,.
\end{align*}
Since there is a canonical embedding, given by Proposition \ref{Cb(X0) embedds in MCc(A)}, of $C_b(X^0)$ into $M(C_c(\a))$, we have in particular an embedding of $C_c(X^0 / H)$ into $M(C_c(\a))$ which identifies an element $f \in C_c(X^0 / H)$ with the multiplier $f \in M(C_c(\a))$ given by:
\begin{align*}
 fb_y := f(\mathbf{r}(y)H)b_y\,.
\end{align*}
Moreover Proposition \ref{Cb(X0) embedds in MCc(A)} applied to the groupoid $X / H$ and the Fell bundle $\a / H$ shows that there is a canonical embedding of $C_b(X^0 / H)$ into $M(C_c(\a/ H))$, which identifies an element $f \in C_b(X^0 / H)$ with the multiplier $f \in M(C_c(\a / H))$ given by
\begin{align}
\label{expression for embedding of CcXH into CcAH}
 f[b]_{yH} := f(\mathbf{r}(y)H)[b]_{yH}\,.
\end{align}
Since both $C_c(X^0 / H)$ and $C_c(\a / H)$ are canonically embedded in $M(C_c(\a))$, it is convenient to understand what happens (inside $M(C_c(\a))$) when one multiplies an element of $C_c(X^0 / H)$ by an element $C_c(\a/H)$. Perhaps unsurprisingly, this product is given exactly by expression (\ref{expression for embedding of CcXH into CcAH}), which models the action of $C_c(X^0/H)$ on $C_c(\a / H)$ as multipliers of the latter algebra. In other words, it makes no difference to view $C_c(X^0 / H)$ inside $M(C_c(\a / H))$ or inside $M(C_c(\a))$ when it comes to multiplication by elements of $C_c(\a / H)$.

We will now show how the multiplication of elements of $C_c(\a / H)$ by elements of $C_c(X^0)$ is determined (inside $M(C_c(\a))$). Before we proceed we will first introduce some notation that will be used throughout this work: Given a set $A \subset X^0$ we will denote by $1_A \in C_b(X^0)$\index{1 1A@$1_A$} the characteristic function of $A$. In case $A$ is a singleton $\{u\}$ we will simply write $1_u$\index{1 1u@$1_u$}.\\

\begin{prop}
\label{multiplication of CcAgm by CcX0 formula}
 Inside $M(C_c(\a))$ we have that, for $x \in X$, $a \in \a_x$ and $u \in X^0$,
\begin{align*}
 [a]_{xH} 1_u = \begin{cases}
 \alpha_{\widetilde{h}^{-1}}(a)_{x\widetilde{h}}\,, \qquad\; \text{if}\;\; H_{x,u} \neq \emptyset \\
 0, \qquad\qquad\qquad \text{otherwise,}\;\;
\end{cases}
\end{align*}
where $\widetilde{h}$ is any element of $H_{x,u}$.\\
\end{prop}

{\bf \emph{Proof:}} Let $y \in X$ and $b \in \a_y$. For the product $[a]_{xH} 1_u\,b_y$ to be non-zero we must necessarily have $u = \mathbf{r}(y)$ (from (\ref{expression for f in Cb(X0) as a multiplier})), and in this case we obtain
\begin{align*}
 [a]_{xH} 1_u\,b_y = [a]_{xH}b_y = (\alpha_{\widetilde{h}^{-1}}(a)b)_{x\widetilde{h}y} = \alpha_{\widetilde{h}^{-1}}(a)_{x\widetilde{h}}b_y\,,
\end{align*}
where $\widetilde{h}$ is any element of $H_{x,y}$. Since $u = \mathbf{r}(y)$, we have $H_{x,y} = H_{x,u}$, and this concludes the proof. \qed\\

It will be of particular importance to know how to multiply, inside $M(C_c(\a))$, elements of $C_c(\a / H)$ with elements of $C_c(\a / K)$ when $K \subseteq H$ is an arbitrary subgroup. It turns out that the algebra $C_c(\a / K)$ is preserved by multiplication by elements of $C_c(\a / H)$, as we show in the next result:\\

\begin{prop}
\label{multiplication of CcAH with CcAK inside MCcA}
 Let $K \subseteq H$ be any subgroup. We have that
\begin{align}
\label{prod of axH with byK inside MCcA}
 [a]_{xH} [b]_{yK} & = \begin{cases}
 [\alpha_{\widetilde{h}^{-1}}(a)b]_{x\widetilde{h}yK}\,, \qquad\; \text{if}\;\; H_{x,y} \neq \emptyset \\
 0, \qquad\qquad\qquad\qquad \text{otherwise,}\;\;
\end{cases}
\end{align}
where $x, y \in X$, $a \in \a_x$ and $b \in \a_y$. In particular $C_c(\a / K)$ is invariant under multiplication by elements of $C_c(\a / H)$.\\
\end{prop}

{\bf \emph{Proof:}} First we observe that since the action is assumed to be $H$-good, it is automatically $K$-good, so that we can form the groupoid $X / K$ and the Fell bundle $\a / K$.

Let $z \in X$ and $c \in \a_z$. We have that
\begin{eqnarray*}
 [a]_{xH} [b]_{yK} c_z \!& = & \! \begin{cases}
 [a]_{xH}(\alpha_{\widetilde{k}^{-1}}(b)c)_{y\widetilde{k}z}\,, \qquad\; \text{if}\;\; K_{y,z} \neq \emptyset \\
 0, \qquad\qquad\qquad\qquad\qquad \text{otherwise,}
\end{cases}\\
& = & \!\begin{cases}
 (\alpha_{\widetilde{h}^{-1}}(a)\alpha_{\widetilde{k}^{-1}}(b)c)_{x\widetilde{h}y\widetilde{k}z}\,, \qquad\; \text{if}\;\; H_{x, y\widetilde{k}z}\neq \emptyset\;\;\text{and}\;\;K_{y,z} \neq \emptyset \\
 0, \qquad\qquad\qquad\qquad\qquad\qquad \text{otherwise,}
\end{cases}\\
& = &\! \begin{cases}
 (\alpha_{\widetilde{h}^{-1}}(a)\alpha_{\widetilde{k}^{-1}}(b)c)_{\big(x\widetilde{h}\widetilde{k}^{-1}y\big)\widetilde{k}z}\,, \quad \text{if}\;\; H_{x, y\widetilde{k}z}\neq \emptyset\;\text{and}\;K_{y,z} \neq \emptyset \\
 0, \qquad\qquad\qquad\qquad\qquad\qquad\quad \text{otherwise,}
\end{cases}\\
\end{eqnarray*}
where $\widetilde{k}$ is any element of $K_{y,z}$ and $\widetilde{h}$ is any element of $H_{x, y\widetilde{k}z}$. Now, since  $H_{x, y\widetilde{k}z} = H_{x, y\widetilde{k}} = H_{x, y}\widetilde{k}$, it follows that $\widetilde{h}\widetilde{k}^{-1} \in H_{x, y}$, and moreover since $K_{y,z} = K_{x\widetilde{h}\widetilde{k}^{-1}y, z}$, we conclude that
\begin{eqnarray*}
& = & \begin{cases}
 (\alpha_{\widetilde{k}^{-1}}(\alpha_{\widetilde{k}\widetilde{h}^{-1}}(a)b)c)_{\big(x\widetilde{h}\widetilde{k}^{-1}y\big)\widetilde{k}z}\,, \qquad\; \text{if}\;\; H_{x, y}\neq \emptyset\;\;\text{and}\;\;K_{x\widetilde{h}\widetilde{k}^{-1}y,z} \neq \emptyset \\
 0, \qquad\qquad\qquad\qquad\qquad\qquad\qquad\quad \text{otherwise,}
\end{cases}\\
& = & \begin{cases}
 [\alpha_{\widetilde{k}\widetilde{h}^{-1}}(a)b]_{x\widetilde{h}\widetilde{k}^{-1}yK}c_{z}\,, \qquad\; \text{if}\;\; H_{x, y}\neq \emptyset\\
 0, \qquad\qquad\qquad\qquad\qquad\quad \text{otherwise.}
\end{cases}
\end{eqnarray*}
Thus (\ref{prod of axH with byK inside MCcA}) follows immediately (the element $\widetilde{h}$ in (\ref{prod of axH with byK inside MCcA}) is simply the element denoted by $\widetilde{h}\widetilde{k}^{-1}$ above). \qed\\

In case the subgroup $K$ has finite index in $H$ we can strengthen Proposition \ref{multiplication of CcAH with CcAK inside MCcA} in the following way:\\

\begin{prop}
\label{axH and sum axhK inside MCcA proposition}
 Let $K \subseteq H$ be a subgroup such that $[H:K] < \infty$. Inside $M(C_c(\a))$ we have that
\begin{align}
\label{axH and sum axhK in statement 111}
 [a]_{xH} \; = \; \sum_{[h] \in \mathcal{S}_x \backslash H / K} [\alpha_{h^{-1}}(a)]_{x h K}\,,
\end{align}
for any $x \in X$ and $a \in \a_x$. In particular, inside $M(C_c(\a))$ we have that $C_c(\a / H)$ is a $^*$-subalgebra of $C_c(\a / K)$.\\
\end{prop}

{\bf \emph{Proof:}} First we notice that since $[H:K] < \infty$ we have that the right hand side of (\ref{axH and sum axhK in statement 111}) is a finite sum and therefore does indeed define an element of $C_c(\a / K)$. To prove this result it suffices to show that
\begin{align}
\label{axH and sum axhK in proof 111}
[a]_{xH} b_y = \sum_{[h] \in \mathcal{S}_x \backslash H / K} [\alpha_{h^{-1}}(a)]_{xh K} b_y\,,
\end{align}
for all $y \in X$ and $b \in \a_y$. First we notice that both the right and left hand sides of (\ref{axH and sum axhK in proof 111}) are zero unless $\mathbf{r}(y) \in \mathbf{s}(x) H$. In case $\mathbf{r}(y) \in \mathbf{s}(x) H$ we have
\begin{align*}
[a]_{xH} b_y = (\alpha_{\widetilde{h}^{-1}}(a)b)_{x\widetilde{h}y}\,,
\end{align*}
where $\widetilde{h}$ is any element of $H_{x, y}$.

Recall from Proposition \ref{xKH bijection with Sx K H} that there is a bijective correspondence between the set of $K$-orbits $(xH)/ K$ and the double coset space $\mathcal{S}_x \backslash H / K$. It is clear that $[a]_{x \widetilde{h} K}b_y = (\alpha_{\widetilde{h}^{-1}}(a)b)_{x \widetilde{h} y}$. Moreover, for all the classes $[h] \neq [\widetilde{h}]$ in $\mathcal{S}_x \backslash H / K$ we have $\mathbf{r}(y) \notin \mathbf{s}(x) h K$, because $\mathbf{r}(y) \in \mathbf{s}(x) \widetilde{h} K$. Hence, for all the classes $[h] \neq [\widetilde{h}]$ in $\mathcal{S}_x \backslash H / K$ we have $[\alpha_{h^{-1}}(a)]_{xhK}b_y = 0$. We conclude that
\begin{align*}
\sum_{[h] \in \mathcal{S}_x \backslash H / K} [\alpha_{h^{-1}}(a)]_{xh K} b_y = [\alpha_{\widetilde{h}^{-1}}(a)]_{x \widetilde{h} K}b_y =  (\alpha_{\widetilde{h}^{-1}}(a)b)_{x\widetilde{h}y}\,,
\end{align*}
and equality (\ref{axH and sum axhK in proof 111}) is proven. \qed\\

\begin{rem}
 In Proposition \ref{axH and sum axhK inside MCcA proposition} the fact that $[H:K] < \infty$ was only used to ensure that the sum on the right hand side of (\ref{axH and sum axhK in statement 111}) was finite. One could more generally just require that the sets $\mathcal{S}_x \backslash H / K$ are finite for all $x \in X$, but this generality will not be used here.\\
\end{rem}

As we saw in Proposition \ref{action on CcA} we have an action $\overline{\alpha}$ of $G$ on $C_c(\a)$. This action can be extended in a unique way to an action on $M(C_c(\a))$, which we will still denote by $\overline{\alpha}$, by the following formula:
\begin{align}
\overline{\alpha}_g(T)f := \overline{\alpha}_g\big(T \overline{\alpha}_{g^{-1}}(f)  \big)\,,
\end{align}
where $g \in G$, $T \in M(C_c(\a))$ and $f \in C_c(\a)$. We will now show what this action on $M(C_c(\a))$ does to the algebras  $C_b(X^0)$, $C_c(\a / H)$ and $C_c(X^0 / H)$.\\

\begin{prop}
\label{a / H --- a / gHg-1}
 The extension of the action $\overline{\alpha}$ to $M(C_c(\a))$, also denoted by $\overline{\alpha}$, satisfies the following properties:
\begin{itemize}
 \item[(i)] The restriction of $\overline{\alpha}$ to $C_b(X^0)$ is precisely the action that comes from the $G$-action on $X^0$.
 \item[(ii)] For any $g \in G$ the automorphism $\overline{\alpha}_g$ takes $C_c(X^0 / H)$ to $C_c(X^0 / gHg^{-1})$, by
\begin{align}
 \overline{\alpha}_g(1_{xH}) = 1_{(xg^{-1})(gHg^{-1})}\,.
\end{align}
\item[(iii)] For any $g \in G$ the automorphism $\overline{\alpha}_g$ takes $C_c(\a / H)$ to $C_c(\a / gHg^{-1})$, by
\begin{align}
 \overline{\alpha}_g([a]_{xH}) = [\alpha_{g}(a)]_{(xg^{-1})(gHg^{-1})}\,.
\end{align}
\item[(iv)] Both $C_c(\a / H)$ and $C_c(X^0 / H)$ are contained in $M(C_c(\a))^H$, the algebra of $H$-fixed points.\\
\end{itemize}
\end{prop}

{\bf \emph{Proof:}} $(i)$ Let $y \in X$, $b \in \a_y$ and $f \in C_b(X^0)$. For any $g \in G$ let us denote by $f_g \in C_b(X^0)$ the function defined by $f_g(x) = f(xg)$. By definition of the extension of $\alpha$ to $M(C_c(\a))$, we have
\begin{eqnarray*}
 \overline{\alpha}_g(f) \, b_y & = & \overline{\alpha}_g(f \cdot \overline{\alpha}_g^{-1}(b_y)) \;\;=\;\; \overline{\alpha}_g(f \cdot \alpha_{g^{-1}}(b)_{yg})\\
& = & \overline{\alpha}_g(f(\mathbf{r}(yg))\alpha_{g^{-1}}(b)_{yg}) \;\; =\;\; \overline{\alpha}_g(f(\mathbf{r}(y)g)\alpha_{g^{-1}}(b)_{yg})\\
& = & f(\mathbf{r}(y)g)b_{y} \;\; = \;\; f_g(\mathbf{r}(y)) b_y\\
& = & f_g \cdot b_y\,.
\end{eqnarray*}
 Hence we conclude that $\overline{\alpha}_g(f) = f_g$ and therefore the action $\overline{\alpha}$ on $C_b(X^0)$ is just the action that comes from the  $G$-action on $X^0$.

$(ii)$ This follows directly from $(i)$.

$(iii)$ Let $y \in X$ and $b \in \a_y$. By definition of the extension of $\overline{\alpha}$ to $M(C_c(\a))$, we have
\begin{align*}
 \overline{\alpha}_g([a]_{xH}) \, b_y = \overline{\alpha}_g([a]_{xH} \overline{\alpha}_g^{-1}(b_y)) = \overline{\alpha}_g([a]_{xH} \alpha_{g^{-1}}(b)_{yg})\,.
\end{align*}
Also, we can see that
\begin{eqnarray*}
 \overline{\alpha}_g([a]_{xH}\alpha_{g^{-1}}(b)_{yg}) & = & \begin{cases}
 \overline{\alpha}_g\big((\alpha_{\widetilde{h}^{-1}}(a) \alpha_{g^{-1}}(b))_{x\widetilde{h}(yg)}\big)\,, \qquad\; \text{if}\;\; H_{x,yg} \neq \emptyset \\
 0, \qquad\qquad\qquad\qquad\qquad\qquad\qquad \text{otherwise}\;\;
\end{cases}\\
& = & \begin{cases}
 (\alpha_{g\widetilde{h}^{-1}}(a)b)_{x\widetilde{h}g^{-1}y}\,, \quad\qquad\; \text{if}\;\; H_{x,yg} \neq \emptyset \\
 0, \qquad\qquad\qquad\qquad\qquad\;\; \text{otherwise}\;\;
\end{cases}\\
& = & \begin{cases}
 (\alpha_{g\widetilde{h}^{-1}}(a)b)_{xg^{-1}g\widetilde{h}g^{-1}y}\,, \quad\;\; \text{if}\;\; H_{x,yg} \neq \emptyset \\
 0, \qquad\qquad\qquad\qquad\qquad\;\; \text{otherwise}\;\;
\end{cases}
\end{eqnarray*}
where $\widetilde{h} \in H_{x, yg}$. Now an easy computation shows that we have
\begin{align*}
 H_{x,yg} = g^{-1} \big(gHg^{-1} \big)_{xg^{-1}, y}\,g\,,
\end{align*}
and thereby we obtain, for $t \in \big(gHg^{-1} \big)_{xg^ {-1},y}$,
\begin{eqnarray*}
 \overline{\alpha}_g([a]_{xH}) \, b_y & = & \begin{cases}
 (\alpha_{gg^{-1}t^{-1}g}(a)b)_{xg^{-1}ty}\,, \qquad \text{if}\;\; \big(gHg^{-1} \big)_{xg^ {-1},y} \neq \emptyset \\
 0, \qquad\qquad\qquad\qquad\qquad\;\; \text{otherwise}\;\;
\end{cases}\\
& = & \begin{cases}
 (\alpha_{t^{-1}g}(a)b)_{xg^{-1}ty}\,, \qquad \text{if}\;\; \big(gHg^{-1} \big)_{xg^ {-1},y} \neq \emptyset \\
 0, \qquad\qquad\qquad\qquad\;\; \text{otherwise}\;\;
\end{cases}\\
& = & [\alpha_{g}(a)]_{(xg^{-1})(gHg^{-1})}\,b_y\,.
\end{eqnarray*}

$(iv)$ This follows directly from $(ii)$ and $(iii)$. \qed\\

It is important to know how to multiply an element of $C_c(\a / H)$ with an element of $C_c(X^0/g Hg^{-1})$ inside $M(C_c(\a))$. This is easy if we are under the assumption that $G$-action satisfies the $H$-intersection property. We recall from (\ref{notation for intersection with conj}) that $H^g$ stands for the subgroup $H \cap g H g^{-1}$.\\

\begin{prop}
\label{product a_xH chi_s(x)gHg-1}
  If the $G$-action moreover satisfies the $H$-intersection property, then for every $x \in X$ and $g \in G$ the following equality holds in $M(C_c(\a))$:
\begin{align*}
 [a]_{xH}\,1_{\mathbf{s}(x)gHg^{-1}} = [a]_{x H^g}\,.\\
\end{align*}
\end{prop}

{\bf \emph{Proof:}} For any $y \in X$ and $b \in \a_y$ we have
\begin{eqnarray*}
 & & [a]_{xH}\, 1_{\mathbf{s}(x)gHg^{-1}}\, b_y \;\;=\;\;\\
 & = & \begin{cases}
 [a]_{xH} b_y\,, \qquad\;\; \text{if}\;\; \mathbf{r}(y) \in \mathbf{s}(x)gHg^{-1} \\
 0, \qquad\qquad\quad \text{otherwise}\;\; 
\end{cases}\\
& = & \begin{cases}
 (\alpha_{\widetilde{h}^{-1}}(a)b)_{x\widetilde{h}y}\,, \qquad\; \text{if}\;\; \mathbf{r}(y) \in \mathbf{s}(x)gHg^{-1} \;\; \text{and}\;\; \mathbf{r}(y) \in \mathbf{s}(x)H \\
 0, \qquad\qquad\qquad\qquad \text{otherwise}\;\; \,
\end{cases}\\
& = & \begin{cases}
 (\alpha_{\widetilde{h}^{-1}}(a)b)_{x\widetilde{h}y}\,, \qquad\; \text{if}\;\; \mathbf{r}(y) \in \mathbf{s}(x)H  \cap \mathbf{s}(x)gHg^{-1} \\
 0, \qquad\qquad\qquad\qquad \text{otherwise}\;\; \,
\end{cases}
\end{eqnarray*}
where $\widetilde{h} \in H_{x,y}$. Now, by the $H$-intersection property, we obtain
\begin{align*}
 = \begin{cases}
 (\alpha_{\widetilde{h}^{-1}}(a)b)_{x\widetilde{h}y}\,, \qquad\; \text{if}\;\; \mathbf{r}(y) \in \mathbf{s}(x)H^g \\
 0, \qquad\qquad\qquad\qquad \text{otherwise}\;\; \,.
\end{cases}
\end{align*}
Of course, we have $(H^g)_{x,y} \subseteq H_{x,y}$, and hence we can choose $\widetilde{h}$ as an element of $(H^g)_{x,y}$, thereby obtaining
\begin{align*}
 = [a]_{xH^g}\,b_y\,, \qquad\qquad\qquad\qquad\qquad
\end{align*}
which finishes the proof. \qed\\

\section{$^*$-Algebraic crossed product by a Hecke pair}
\label{algebraic crossed products chapter}

In this section we introduce our notion of a ($^*$-algebraic) crossed product by a Hecke pair and we explore its basic properties and its representation theory. Throughout the rest of this work we impose the following standing assumption, based on the tools developed in Section \ref{group actions on discrete groupoids section}.\\

\begin{standing}
\label{standing assumption}
 We assume from now on that $(G,\gm)$ is a Hecke pair, $\a$ is a Fell bundle over a groupoid $X$ endowed with a $\gm$-good right $G$-action $\alpha$ satisfying the $\gm$-intersection property.\\
\end{standing}

\subsection{Definition of the crossed product and basic properties}

In this section we aim at defining the ($^*$-algebraic) crossed product of $C_c(\a / \gm)$ by the Hecke pair $(G,\gm)$. For that we are going to define some sort of a bundle over $G / \gm$, where the fiber over each $g \gm$ is precisely $C_c(\a / \gm^g)$. Recall that we denote by $\overline{\alpha}$ the associated action of $G$ on $C_c(\a)$ and also its extension to $M(C_c(\a))$.\\

\begin{df}
 Let $B(\a, G, \gm)$\index{BAGgm@$B(\a, G, \gm)$} be the vector space of finitely supported functions $f : G / \gm \to M(C_c(\a))$ satisfying the following compatibility condition
\begin{align}
\label{compatibility cond}
 f(\gamma g \gm) = \overline{\alpha}_{\gamma}(f(g \gm))\,,
\end{align}
for all $\gamma \in \gm$ and $g \gm \in G / \gm$.\\
\end{df}

\begin{lemma}
\label{fibers are gm^g invariant}
 For every $f \in B(\a, G, \gm)$ and $g\gm \in G / \gm$ we have
\begin{align*}
 f(g\gm) \in M(C_c(\a))^{\gm^g}\,.\\
\end{align*}
\end{lemma}

{\bf \emph{Proof :}} This follows directly from the compatibility condition (\ref{compatibility cond}), since for every $\gamma \in \gm^g$ we have $\overline{\alpha}_{\gamma}(f(g\gm)) = f(\gamma g \gm) = f(g \gm)$. \qed\\

\begin{df}
\label{def of crossed product by hecke pair}
 The vector subspace of $B(\a, G, \gm)$ consisting of the functions $f:G / \gm \to  M(C_c(\a))$ satisfying the compatibility condition (\ref{compatibility cond}) and the property
\begin{align}
\label{f(g gm) in Cc (A Gmg)}
 f(g\gm) \in C_c(\a / \gm^g)\,,
\end{align}
will be denoted by $C_c(\a / \gm) \times_{\alpha}^{alg} G / \gm $\index{crossed product (algebraic)@$C_c(\a / \gm) \times_{\alpha}^{alg} G / \gm $} and will be called the \emph{$^*$-algebraic crossed product} of $C_c(\a / \gm)$ by the Hecke pair $(G, \gm)$.\\
\end{df}

 It is relevant to point out that the definitions of the spaces $B(\a, G , \gm)$ and $C_c(\a / \gm) \times_{\alpha}^{alg} G / \gm $ seem more suitable for Hecke pairs $(G, \gm)$, as in general a function in $B(\a, G, \gm)$ could only have support on those elements $g\gm \in G / \gm$ such that $|\gm g\gm / \gm| < \infty$.

We now define a product and an involution in $B(\a, G, \gm)$ by:
\begin{align}
\label{conv prod}
 (f_1 * f_2)(g \gm) & := \sum_{[h] \in G/ \gm} f_1(h \gm)\, \overline{\alpha}_{h}(f_2(h^{-1} g\gm))\,,\\
\label{invol def} (f^*)\,(g \gm) & := \Delta(g^{-1})\overline{\alpha}_{g}(f(g^{-1} \gm))^*\,.
\end{align}\\

\begin{prop}
\label{B is an algebra}
 $B(\a, G, \gm)$ becomes a unital $^*$-algebra under the product and involution defined above, whose identity element is the function $f$ such that $f(\gm) = 1$ and is zero in the remaining points of $G / \gm$.\\
\end{prop}

{\bf \emph{Proof:}} First, we claim that the expression for the product defined above is well-defined in $B(\a, G, \gm)$, i.e. for $f_1, f_2 \in B(\a, G, \gm)$ the expression
\begin{align*}
 (f_1*f_2)(g\gm) := \sum_{[h] \in G/ \gm} f_1(h \gm)\, \overline{\alpha}_{h}(f_2(h^{-1} g\gm))
\end{align*}
is independent from the choice of the representatives $[h]$ and also that it has finitely many summands. Independence from the choice of the representatives $[h] \in G / \gm$ follows directly from the compatibility condition (\ref{compatibility cond}) and the fact that the sum is finite follows simply from the fact that $f_1$ has finite support.

Now we claim that $f_1 * f_2$ has also finite support, for $f_1, f_2 \in B(\a, G, \gm)$. Let $S_1, S_2 \subseteq G / \gm$ be the supports of the functions $f_1$ and $f_2$ respectively. We will regard $S_1$ and $S_2$ as subsets of $G$ (being finite unions of left cosets). It is easy to check that the function $G \times G \to M(C_c(\a))$
\begin{align*}
 (h,g) \mapsto f_1(h\gm) \overline{\alpha}_h(f_2(h^{-1}g\gm))
\end{align*}
has support contained in $S_1 \times (S_1 \cdot S_2)$. Since $(G, \gm)$ is a Hecke pair, the product $S_1 \cdot S_2$ is also a finite union of left cosets. Hence, $f_1 * f_2$ has finite support.

We also notice that $f_1* f_2$ satisfies the compatibility condition (\ref{compatibility cond}), thus defining an element of $B(\a, G, \gm)$, since for any $ \gamma \in \gm$ we have
\begin{eqnarray*}
 (f_1*f_2) ( \gamma g\gm) & = & \sum_{[h] \in G/ \gm} f_1(h \gm)\, \overline{\alpha}_{h}(f_2(h^{-1} \gamma g\gm))\\
 & = & \sum_{[h] \in G/ \gm} f_1(\gamma h \gm)\, \overline{\alpha}_{\gamma h}(f_2(h^{-1}g\gm))\\
 & = & \sum_{[h] \in G/ \gm} \overline{\alpha}_{\gamma}(f_1( h \gm))\, \overline{\alpha}_{\gamma} \circ \overline{\alpha}_{h}(f_2(h^{-1}g\gm))\\
 & = & \overline{\alpha}_{\gamma} \big( (f_1 * f_2)(g\gm) \big)\,.
\end{eqnarray*}

In a similar way we can see that the expression that defines the involution is well-defined in $B(\a, G, \gm)$. There are now a few things that need to be checked before we can say that $B(\a, G, \gm)$ is a $^*$-algebra, namely that the product is associative and the involution is indeed an involution relatively to this product (the fact that the product is distributive and the properties concerning multiplication by scalars are obvious). The proofs of these facts are essentially just a mimic of the corresponding proofs for ``classical'' crossed products by groups. Thus, we can say that $B(\a, G,\gm)$ is $^*$-algebra under this product and involution. \qed\\

\begin{thm}
 $C_c(\a / \gm) \times_{\alpha}^{alg} G / \gm $ is a $^*$-ideal of $B(\a, G, \gm)$. In particular it is a $^*$-algebra for the above operations.\\

\end{thm}

{\bf \emph{Proof:}} It is easy to see that the space $C_c(\a / \gm) \times_{\alpha}^{alg} G / \gm $ is invariant for  the involution, i.e.
\begin{align*}
 f \in C_c(\a / \gm) \times_{\alpha}^{alg} G / \gm\;  \Longrightarrow\; f^* \in C_c(\a / \gm) \times_{\alpha}^{alg} G / \gm\,.
\end{align*}

Thus, to prove that $C_c(\a / \gm) \times_{\alpha}^{alg} G / \gm $ is a (two-sided) $^*$-ideal of $B(\a, G, \gm)$ it is enough to prove that is a right ideal, i.e. if $f_1 \in B(\a, G, \gm)$ and $f_2 \in C_c(\a / \gm) \times_{\alpha}^{alg} G / \gm$ then $f_1*f_2 \in C_c(\a / \gm) \times_{\alpha}^{alg} G / \gm$, because any right $^*$-ideal is automatically two-sided. Hence, all we need to prove is that $(f_1 * f_2)(g\gm) \in C_c(\a / \gm^g)$, for every $f_1 \in B(\a, G, \gm)$ and $f_2 \in C_c(\a / \gm) \times_{\alpha}^{alg} G / \gm $. The proof of this fact will follow the following steps:
\begin{itemize}
 \item[1)]  Prove that: given a subgroup $H \subseteq G$, $f \in C_c(\a/ H)$ and a unit $u \in X^0$, we have $f \cdot 1_u \in C_c(\a)$.
 \item[2)] Let $T:= (f_1*f_2)(g\gm)= \sum_{[h] \in G / \gm} f_1(h\gm) \overline{\alpha}_h(f_2(h^{-1}g\gm))$. Use $1)$ to show that $T\cdot 1_u \in C_c(\a)$ for any unit $u \in X^0$.
 \item[3)] Fix a unit $u \in X^0$. By $2)$ we have $T\, 1_u = \sum_{i} (a_i)_{x_i}$, where the elements $x_i \in X$ are such that $\mathbf{s}(x_i) = u$. Show that $T\, 1_{u \gm^g} = \sum_i [a_i]_{x_i \gm^g}$, and conclude that $T\, 1_{u \gm^g} \in C_c(\a / \gm^g)$.
 \item[4)] Prove that there exists a finite set of units $\{u_1, \dots, u_n\} \subseteq X^0$ such that $T = \sum_{i=1}^n T\, 1_{u_i\gm^g}$. Conclude that $ T \in C_c(\a / \gm^g)$.\\
\end{itemize}

\begin{itemize}
 \item Proof of 1) : $\quad$ This follows immediately from Proposition \ref{multiplication of CcAgm by CcX0 formula}.
 \item Proof of 2) : $\quad$ We know that $f_2(h^{-1}g \gm) \in C_c(\a / \gm^{h^{-1}g})$. Thus, from Proposition \ref{a / H --- a / gHg-1}, we conclude that $\overline{\alpha}_h(f_2(h^{-1}g)) \in C_c\big(\a/ h \gm h^{-1} \cap g \gm g^{-1}\big)$. Now, using $1)$, we see that $\overline{\alpha}_h(f_2(h^{-1}g))\, 1_{u} \in C_c(\a )$ and consequently $f_1(h\gm)\overline{\alpha}_h(f_2(h^{-1}g))\,1_u \in C_c(\a)$. Hence, $T \,1_u \in C_c(\a)$.
\item Proof of 3) : $\quad$ For any $\gamma \in \gm^g$ we have, using Lemma \ref{fibers are gm^g invariant},
\begin{eqnarray*}
 T\,1_{u\gamma} & = & \overline{\alpha}_{\gamma^{-1}}(T)\, 1_{u \gamma} = \overline{\alpha}_{\gamma^{-1}} \big(T \overline{\alpha}_{\gamma}(1_{u \gamma})\big)\\
& = & \overline{\alpha}_{\gamma^{-1}} \big( T \, 1_u \big) = \sum_{i} \alpha_{\gamma^{-1}}(a_i)_{x_i \gamma}\,.
\end{eqnarray*}
Let $y \in X$ and $b \in \a_y$. We have
\begin{eqnarray*}
 T \, 1_{u \gm^g}\, b_y & =  & \begin{cases}
 T b_y\,, \qquad\; \text{if}\;\; \mathbf{r}(y) \in u\gm^g \\
 0, \qquad\quad\;\; \text{otherwise}\;\; .
\end{cases}
\end{eqnarray*}
Assume now that $\mathbf{r}(y) \in u \gm^g$ and let $\widetilde{\gamma} \in \gm^g$ be such that $\mathbf{r}(y) = u \widetilde{\gamma}$. We then have
\begin{align*}
 T b_y \;\; = \;\; T \, 1_{u \widetilde{\gamma}}\, b_y \;\; = \;\; \sum_{i} \alpha_{\widetilde{\gamma}^{-1}}(a_i)_{x_i \widetilde{\gamma}}\, b_y\,.
\end{align*}
Since $\mathbf{s}(x_i) = u$, we have $\mathbf{s}(x_i \widetilde{\gamma}) = u \widetilde{\gamma} = \mathbf{r}(y)$. Hence,
\begin{align*}
 Tb_y \;\; = \;\; \sum_{i} (\alpha_{\widetilde{\gamma}^{-1}}(a_i) b)_{x_i \widetilde{\gamma} y}\,.
\end{align*}
We conclude that
\begin{eqnarray*}
T \, 1_{u \gm^g}\, b_y & = & \begin{cases}
 \sum_{i} (\alpha_{\widetilde{\gamma}^{-1}}(a_i) b)_{x_i \widetilde{\gamma} y}\,, \qquad\; \text{if}\;\; \mathbf{r}(y) \in u\gm^g \\
 0, \qquad\qquad\qquad\qquad\qquad\; \text{otherwise}\;\; .
\end{cases}\\
& = & \sum_{i} [a_i]_{x_i \gm^g} \, b_y\,.
\end{eqnarray*}
Thus, $T\, 1_{u \gm^g} = \sum_i [a_i]_{x_i\gm^g} \in C_c(\a / \gm^g)$.
\item Proof of 4) : $\quad$ For easiness of reading of this last part of the proof we introduce the following definition: given $F \in M(C_c(\a))$ we define the \emph{support} of $F$ to be the set
$\{u \in X^0 : F \,1_u \neq 0\}$. Notice in particular that the support of an element $[a]_{xH}$, with $a \neq 0$, is the set $\mathbf{s}(x)H$.

 Since $\overline{\alpha}_h(f_2(h^{-1}g\gm)) \in C_c(\a / h \gm h^{-1} \cap g \gm g^{-1})$, there exists a finite number of units $v_1, \dots, v_k \in X^0$ such that $\overline{\alpha}_h(f_2(h^{-1}g\gm))$ has support in 
\begin{align*}
 \bigcup_{i = 1}^k v_i \big(h \gm h^{-1} \cap g \gm g^{-1} \big)\; \subseteq \; \bigcup_{i=1}^k v_i  g \gm g^{-1}\,.
\end{align*}
Hence, there is a finite number of units $w_1, \dots, w_l \in X^0$ such that $T$ has support contained in
\begin{align*}
 \bigcup_{i=1}^l w_i g \gm  g^{-1}\,.
\end{align*}
Therefore, $T$ has support contained in
\begin{align*}
 \bigcup_{i=1}^l \bigcup_{j=1}^m w_i \theta_j  \gm^g\,,
\end{align*}
where $\theta_1, \dots, \theta_m$ are representatives of the classes of $ g \gm g^{-1} / \gm^g$ (being a finite number because $(G,\gm)$ is a Hecke pair). Thus, we have proven that there is a finite number of units  $u_1, \dots, u_n \in X^0$ such that $T$ has support inside $\bigcup_{i=1}^n u_i \gm^g$. Moreover, we can suppose we have chosen the units $u_1, \dots, u_n$ such that the corresponding orbits $u_i \gm^g$ are mutually disjoint. It is now easy to see that we have $T = \sum_{i=1}^n T \, 1_{u_i \gm^g}$. Indeed, given $y \in X$ and $b \in \a_y$, if $\mathbf{r}(y) \notin \bigcup_{i=1}^n u_i \gm ^g$, then
\begin{align*}
 T b_y = T \, 1_{\mathbf{r}(y)}\, b_y= 0 = \sum_{i=1}^n T \, 1_{u_i \gm^g} \, b_y \,,
\end{align*}
 and if $\mathbf{r}(y) \in \bigcup_{i=1}^n u_i \gm ^g$, then $\mathbf{r}(y)$ belongs to precisely one of the orbits, say $u_{i_0} \gm^g$, and we have
\begin{align*}
 \sum_{i=1}^n T \, 1_{u_i \gm^g} \, b_y = T 1_{u_{i_0} \gm^g} \,b_y = T b_y\,.
\end{align*}
Hence, we must have $T = \sum_{i=1}^n T\, 1_{u_i \gm^g}$, and by $3)$ we conclude that $T \in C_c(\a / \gm^g)$. \qed\\
\end{itemize}

As it is well-known, when working with crossed products $A \times G$ by discrete groups, one always has an embedded copy of $A$ inside the crossed product. Something analogous happens in the case of crossed products by Hecke pairs, where $C_c(\a / \gm)$ is canonically embedded in $C_c(\a / \gm) \times_{\alpha}^{alg} G / \gm$, as is stated in the next result (whose proof amounts to routine verification).\\

\begin{prop}
\label{embedding of C_c(A  gm)}
 There is a natural embedding of the $^*$-algebra $C_c(\a / \gm)$ in $C_c(\a / \gm) \times_{\alpha}^{alg} G / \gm$, which identifies an element $f \in C_c(\a / \gm)$ with the function $\iota(f) \in C_c(\a / \gm) \times_{\alpha}^{alg} G / \gm$ such that
\begin{align*}
 \iota(f)(\gm) = f \qquad \text{and} \qquad \iota(f)\;\, \text{is zero elsewhere}\,.\\
\end{align*}
\end{prop}

\begin{rem}
 The above result says that we can identify $C_c(\a / \gm)$ with the functions of $C_c(\a / \gm) \times_{\alpha}^{alg} G / \gm$ with support in $\gm$. We shall, henceforward, make no distinctions in notation between an element of $C_c(\a /\gm)$ and its correspondent in $C_c(\a / \gm) \times_{\alpha}^{alg} G / \gm$.\\
\end{rem}

\begin{thm}
\label{cross is essential ideal of B}
 $C_c(\a / \gm) \times_{\alpha}^{alg} G / \gm$ is an essential $^*$-ideal of $B(\a, G, \gm)$. In particular, $C_c(\a / \gm) \times_{\alpha}^{alg} G / \gm$ is an essential $^*$-algebra. Moreover, there are natural embeddings
\begin{align*}
 C_c(\a / \gm) \times_{\alpha}^{alg} G / \gm \hookrightarrow B(\a, G, \gm) \hookrightarrow M(C_c(\a / \gm) \times_{\alpha}^{alg} G / \gm )\,,
\end{align*}
that make the following diagram commute
\begin{displaymath}
\xymatrix{ & M(C_c(\a / \gm) \times_{\alpha}^{alg} G / \gm ) \\
C_c(\a / \gm) \times_{\alpha}^{alg} G / \gm \ar[ur]^L \ar[r] & B(\a, G, \gm)\,. \ar[u]}
\end{displaymath}\\

\end{thm}

{\bf \emph{Proof:}} We have already proven that $C_c(\a / \gm) \times_{\alpha}^{alg} G / \gm$ is a $^*$-ideal of $B(\a, G, \gm)$, thus we only need to check that this ideal is in fact essential. Suppose $f \in B(\a, G, \gm)$ is such that $f * \big(C_c(\a / \gm) \times_{\alpha}^{alg} G / \gm \big) = \{0\}$. Then, in particular, using Proposition \ref{embedding of C_c(A  gm)}, we must have $f *  \big(C_c(\a / \gm) \big) = \{0\}$. Let $g \in G$ and take $[a]_{x\gm} \in C_c(\a / \gm)$, we then have
\begin{align*}
0 = \Big(f* [a]_{x\gm} \Big)\; (g\gm) = f(g\gm) \overline{\alpha}_g([a]_{x\gm}) = f(g \gm) [\alpha_g(a)]_{xg^{-1}g \gm g^{-1}}\,.
\end{align*}
Thus, multiplying by $1_{\mathbf{s}(x)g^{-1}} \in M(C_c(\a))$ we get
\begin{align*}
 0 = f(g \gm) [\alpha_g(a)]_{xg^{-1}g \gm g^{-1}} 1_{\mathbf{s}(x)g^{-1}} = f(g\gm)\alpha_g(a)_{xg^{-1}} = f(g\gm)\overline{\alpha}_g(a_{x})\,.
\end{align*}
Since this true for all $a \in \a_x$ and $x \in X$ and given that $\alpha$ takes fibers of $\a$ bijectively into fibers of $\a$, we must have $f(g\gm)b_y = 0$ for all $b \in \a_y$ and $y \in X$. Hence, we must have $f(g\gm) = 0$.  Thus, $f = 0$ and we conclude that $C_c(\a / \gm) \times_{\alpha}^{alg} G / \gm$ is an essential $^*$-ideal of $B(\a, G, \gm)$.

Since $C_c(\a / \gm) \times_{\alpha}^{alg} G / \gm$ is a $^*$-subalgebra of $B(\a, G, \gm)$, we immediately conclude that $C_c(\a / \gm) \times_{\alpha}^{alg} G / \gm$ is an essential $^*$-algebra.

The embedding of $B(\a, G, \gm)$ in $M(C_c(\a / \gm) \times_{\alpha}^{alg} G / \gm)$ then follows from the universal property of multiplier algebras, Theorem \ref{property of M(A)}. \qed\\

In the theory of crossed products $A \times G$ by groups, one always has an embedded copy of the group algebra $\mathbb{C}(G)$ inside the multiplier algebra $M(A \times G)$. Something analogous happens in the case of crossed products by Hecke pairs, where the Hecke algebra $\h(G, \gm)$ is canonically embedded in the multiplier algebra $M(C_c(\a / \gm) \times_{\alpha}^{alg} G / \gm)$, as is stated in the next result (whose proof amounts to routine verification).\\

\begin{prop}
\label{embedding of the Hecke alg}
 The Hecke $^*$-algebra $\h(G, \gm)$ embeds in $B(\a, G, \gm)$ in the following way: an element $f \in \h(G, \gm)$ is identified with the element $\widetilde{f} \in B(\a, G, \gm)$ given by
\begin{align*}
 \widetilde{f}(g\gm) := f( g \gm)\mathbf{1}\,,
\end{align*}
where $\mathbf{1}$ is the unit of $M(C_c(\a))$.\\
\end{prop}

The next result does not typically play an essential role in the case of crossed products by groups, but will be extremely important for us in case of crossed products by Hecke pairs. The proof is also just routine verification.\\

\begin{prop}
\label{C_c(U gm) embeds in B}
 The algebra $C_c(X^0 / \gm)$ embeds in $B(\a, G, \gm)$ in the following way: an element $f \in C_c(X^0 / \gm)$ is identified with the function $\iota(f) \in B(\a, G, \gm)$ given by
\begin{align*}
 \iota(f)(\gm) = f \qquad \text{and} \qquad \iota(f)\;\, \text{is zero elsewhere}\,.\\
\end{align*}
\end{prop}

\begin{rem}
Propositions \ref{embedding of the Hecke alg} and \ref{C_c(U gm) embeds in B} allow us to view both the Hecke $^*$-algebra $\h(G, \gm)$ and $C_c(X^0 / \gm)$ as  $^*$-subalgebras of $B(\a, G, \gm)$. We shall henceforward make no distinctions in notation between an element of $\h (G, \gm)$ or $C_c(X^0 /\gm)$ and its correspondent in $B(\a, G, \gm)$.\\
\end{rem}

The purpose of the following diagram is to illustrate, in a more condensed form, all the canonical embeddings we have been considering so far:

\begin{displaymath}
\xymatrix{ C_c(\a / \gm) \ar[r] & C_c(\a / \gm) \times_{\alpha}^{alg} G / \gm \ar[dr] \\
  \h (G,\gm) \ar[rr] & & B(\a, G , \gm) \ar[r] & M(C_c(\a / \gm) \times_{\alpha}^{alg} G / \gm))\\
  C_c(X^0 / \gm) \ar[urr] }
\end{displaymath}\\

\begin{rem}
\label{reasons for considering BAGgm}
The reason for considering the algebra $B(\a, G, \gm)$ is two-fold. On one side $B(\a, G, \gm)$ made it easier to make sure the convolution product (\ref{conv prod}) was well-defined in $C_c(\a / \gm) \times_{\alpha}^{alg} G / \gm$. On the other (perhaps more important) side, the fact that both $\h(G, \gm)$ and $C_c(X^0 / \gm)$ are canonically embedded in $B(\a, G, \gm)$ allows us to treat the elements of $\h(G, \gm)$ and $C_c(X^0 / \gm)$ both as multipliers in $M(C_c(\a / \gm) \times_{\alpha}^{alg} G / \gm))$, but also allows us to operate these elements with the convolution product and involution expressions (\ref{conv prod}) and (\ref{invol def}), as these are defined in $B(\a, G, \gm)$.\\
\end{rem}

As it is well-known in the theory of crossed products by discrete groups,  a ($^*$-algebraic) crossed product $A \times G$ is spanned by elements of the form $a * g$, where $a \in A$ and $g \in G$ (here $g$ is seen as an element of the group algebra $\mathbb{C}(G) \subseteq M(A \times G)$). We will now explore something analogous in the case of crossed products by Hecke pairs. It turns out that $C_c(\a / \gm) \times_{\alpha}^{alg} G / \gm$ is spanned by elements of the form $[a]_{x\gm} * \gm g \gm * 1_{\mathbf{s}(x) g \gm}$, where $x \in X$, $a \in \a_x$ and $g\gm \in G / \gm$, as we show in the next result.\\

\begin{thm}
\label{prop decomp of f in crossed prod}
 For any $f \in C_c(\a / \gm) \times_{\alpha}^{alg} G / \gm$ we have
\begin{align}
\label{decomp of f in crossed prod}
 f = \sum_{[g] \in \gm \backslash G / \gm}\; \sum_{x\gm^g \in X / \gm^g} \Big[f(g\gm) (x) \Big]_{x\gm} * \gm g \gm * 1_{\mathbf{s}(x) g \gm}\,.
\end{align}
In particular, $C_c(\a / \gm) \times_{\alpha}^{alg} G / \gm$ is spanned by elements of the form
\begin{align*}
 [a]_{x\gm} * \gm g \gm * 1_{\mathbf{s}(x) g \gm}\,,
\end{align*}
with $x \in X$, $a \in \a_x$ and $g\gm \in G/\gm$.\\
\end{thm}

The following lemma is needed in order to prove the above result:\\

\begin{lemma}
\label{a_xH HgH chi_s(x)gHg}
 Let $x \in X$, $a \in \a_x$ and $g\gm \in G / \gm$. We have
\begin{align*}
 [a]_{x\gm} * \gm g \gm * 1_{\mathbf{s}(x) g \gm} \;\, (h \gm) = \begin{cases}
 [\alpha_{\gamma}(a)]_{x \gamma^{-1} \gm^{\gamma g}}\,, \qquad \text{if}\;\; h\gm = \gamma g\gm, \;\;\text{with}\;\; \gamma \in \gm \\
 0, \qquad\qquad\qquad\quad\;\; \text{otherwise}\;\;.
\end{cases}
\end{align*}
In particular,
\begin{align*}
 [a]_{x\gm} * \gm g \gm * 1_{\mathbf{s}(x) g \gm} \;\, (g \gm) =  [a]_{x \gm^{ g}}\,.
\end{align*}\\
\end{lemma}

{\bf \emph{Proof:}} An easy computation yields
\begin{align*}
 [a]_{x\gm} * \gm g \gm * 1_{\mathbf{s}(x) g \gm} \;\, (h \gm)\, = \, [a]_{x \gm}\cdot\gm g \gm (h \gm) \cdot \overline{\alpha}_h(1_{\mathbf{s}(x)g \gm})\,,
\end{align*}
from which we conclude that $[a]_{x\gm} * \gm g \gm * 1_{\mathbf{s}(x) g \gm}$ is supported in the double coset $\gm g \gm$. Now, evaluating at the point $g\gm \in G / \gm$ we get
\begin{eqnarray*}
 [a]_{x\gm} * \gm g \gm * 1_{\mathbf{s}(x) g \gm} \;\, (g \gm) & = & [a]_{x \gm}\cdot\gm g \gm (g \gm) \cdot \overline{\alpha}_g(1_{\mathbf{s}(x)g \gm})\\
& = & [a]_{x \gm} \cdot \overline{\alpha}_g(1_{\mathbf{s}(x)g \gm})\\
& = & [a]_{x \gm} \cdot 1_{\mathbf{s}(x)g\gm g^{-1}}\\
& = & [a]_{x \gm^g}\,,
\end{eqnarray*}
where the last equality comes from Proposition \ref{product a_xH chi_s(x)gHg-1}. From the compatibility condition (\ref{compatibility cond}) and Proposition \ref{a / H --- a / gHg-1} it then follows that, for $\gamma \in \gm$,
\begin{eqnarray*}
 [a]_{x\gm} * \gm g \gm * 1_{\mathbf{s}(x) g \gm} \;\, (\gamma g \gm) & = & \overline{\alpha}_{\gamma}([a]_{x \gm^g})\\
& = & [\alpha_{\gamma}(a)]_{x \gamma^{-1} \gm^{\gamma g}}\,.
\end{eqnarray*}\qed\\

{\bf \emph{Proof of Theorem \ref{prop decomp of f in crossed prod}:}} Let us first prove that the expression on the right hand side of (\ref{decomp of f in crossed prod}) is well-defined. It is easy to see that for every $g  \in G$, the expression
\begin{align*}
 \sum_{x\gm^g \in X / \gm^g} \Big[f(g\gm) (x) \Big]_{x\gm} * \gm g \gm * 1_{\mathbf{s}(x) g \gm}
\end{align*}
does not depend on the choice of the representative $x$ of $x\gm^g$. Now, let us see that it also does not depend on the choice of the representative $g$ in $\gm  g  \gm$. Let $\gamma g \theta$, with $\gamma, \theta \in \gm$, be any other representative. We have
\begin{eqnarray*}
&  & \sum_{x\gm^{\gamma g \theta} \in X / \gm^{\gamma g \theta}} \Big[f(\gamma g \theta\gm) (x) \Big]_{x\gm} * \gm \gamma g \theta \gm * 1_{\mathbf{s}(x) \gamma g \theta \gm} =\\
& = & \sum_{x\gm^{\gamma g} \in X / \gm^{\gamma g}} \Big[f(\gamma g\gm) (x) \Big]_{x\gm} * \gm g\gm * 1_{\mathbf{s}(x) \gamma g\gm}\\
& = & \sum_{x\gm^{\gamma g} \in X / \gm^{\gamma g}} \Big[\overline{\alpha}_{\gamma}( f(g\gm)) (x) \Big]_{x\gm} * \gm g\gm * 1_{\mathbf{s}(x) \gamma g\gm}\\
& = & \sum_{x\gm^{\gamma g} \in X / \gm^{\gamma g}} \Big[\alpha_{\gamma}( f(g\gm) (x \gamma)) \Big]_{x\gm} * \gm g\gm * 1_{\mathbf{s}(x) \gamma g\gm}
\end{eqnarray*}
We notice that there is a well-defined bijective correspondence $X / \gm^g \to X / \gm^{\gamma g}$ given by $x \gm^g \mapsto x\gamma^{-1} \gm^{\gamma g}$. Thus, we get
\begin{eqnarray*}
& = & \sum_{x\gm^g \in X / \gm^{g}} \Big[\alpha_{\gamma}(f(g\gm) (x )) \Big]_{x \gamma^{-1}\gm} * \gm g\gm * 1_{\mathbf{s}(x \gamma^{-1}) \gamma g\gm}\\
& = & \sum_{x\gm^g \in X / \gm^{g}} \Big[f(g\gm) (x ) \Big]_{x \gm} * \gm g\gm * 1_{\mathbf{s}(x) g\gm}\,.
\end{eqnarray*}
Hence, the expression in (\ref{decomp of f in crossed prod}) is well-defined. Let us now prove the decomposition in question. For any $t \gm \in G / \gm$ we have
\begin{eqnarray*}
&  & \sum_{[g] \in \gm \backslash G / \gm}\; \sum_{x\gm^g \in X / \gm^{g}} \Big[f(g\gm) (x ) \Big]_{x \gm} * \gm g\gm * 1_{\mathbf{s}(x) g\gm}\; (t\gm) =\\
& = & \sum_{x\gm^t \in X / \gm^t} \Big[f(t\gm) (x) \Big]_{x \gm} * \gm t\gm * 1_{\mathbf{s}(x) t\gm}\; (t\gm)\,.
\end{eqnarray*}
By Lemma \ref{a_xH HgH chi_s(x)gHg} it follows that
\begin{eqnarray*}
& = & \sum_{x\gm^t \in X / \gm^t} \Big[f(t\gm) (x) \Big]_{x \gm^t} \\
& = & f(t\gm)\,,
\end{eqnarray*}
and this finishes the proof. \qed\\

In the following result we collect some useful equalities concerning products in $C_c(\a / \gm) \times_{\alpha}^{alg} G / \gm$, which will be useful later on. One should observe the similarities between the equalities (\ref{a x gm times gm g gm expression}) and (\ref{gm g gm times a x gm expression}) and the equalities obtained by an Huef, Kaliszewski and Raeburn in \cite[Lemma 1.3 (i) and (ii)]{cov} if in their setting one was allowed to somehow ``drop'' the representations. The similarity is more than a coincidence and will be addressed in the sequel of this article.\\

\begin{prop}
\label{invol in crossed prod}
 In $C_c(\a / \gm) \times_{\alpha}^{alg} G / \gm$ the following equalities hold:
\begin{align}
\label{invol for span elements expression}
 \big( [a]_{x\gm} * \gm g \gm * 1_{\mathbf{s}(x) g \gm} \big)^* & = \Delta(g)\; [\alpha_{g^{-1}}(a^*)]_{x^{-1}g\gm} * \gm g^{-1} \gm * 1_{\mathbf{s}(x^{-1})\gm}\,,
\end{align}
\begin{align}
\label{alt prod expression}
 1_{\mathbf{r}(x) \gm} *  \gm g \gm *  [\alpha_{g^{-1}}(a)]_{xg\gm} & = [a]_{x \gm} * \gm g \gm * 1_{\mathbf{s}(x) g \gm}\,,
\end{align}
\begin{align}
\label{a x gm times gm g gm expression}
[a]_{x \gm}*\gm g \gm  & = \sum_{[\gamma] \in \mathcal{S}_{x} \backslash \gm / \gm^g} [a]_{x \gm} * \gm g  \gm* 1_{\mathbf{s}(x) \gamma g \gm}   \,.
\end{align}
\begin{align}
\label{gm g gm times a x gm expression}
\gm g \gm * [a]_{x \gm}  = \sum_{[\gamma] \in \mathcal{S}_{x} \backslash \gm / \gm^{g^{-1}}} 1_{\mathbf{r}(x) \gamma g^{-1} \gm} * \gm g  \gm * [a]_{x \gm}\,.
\end{align}
In particular, from (\ref{alt prod expression}) we see that $C_c(\a / \gm) \times_{\alpha}^{alg} G / \gm$ is also spanned by all elements of the form $1_{\mathbf{r}(x) \gm} *  \gm g \gm *  [a]_{xg\gm}$, with $g \in G$, $x \in X$ and $a \in \a_x$.\\
\end{prop}

{\bf \emph{Proof:}} Let us first prove equality (\ref{invol for span elements expression}). First we notice that
\begin{align*}
 \big( [a]_{x\gm} * \gm g \gm * 1_{\mathbf{s}(x) g \gm} \big)^* = \Delta(g)\; 1_{\mathbf{s}(x) g \gm} *  \gm g^{-1} \gm *  [a^*]_{x^{-1}\gm}\,,
\end{align*}
which means that $\big( [a]_{x\gm} * \gm g \gm * 1_{\mathbf{s}(x) g \gm} \big)^*$ has support in the double coset $\gm g^{-1}\gm$. Now evaluating this element on $g^{-1} \gm$ we get,
\begin{eqnarray*}
 & & \big( [a]_{x\gm} * \gm g \gm * 1_{\mathbf{s}(x) g \gm} \big)^* \,(g^{-1}\gm) \;\; =\;\;\\
 & = & \Delta(g)\; \overline{\alpha}_{g^{-1}} \big( ( [a]_{x\gm} * \gm g \gm * 1_{\mathbf{s}(x) g \gm} ) \,(g\gm) \big)^* \\
& = & \Delta(g)\; \overline{\alpha}_{g^{-1}}( [a]_{x \gm^g} )^*\\
& = & \Delta(g)\; [\alpha_{g^{-1}}(a^*)]_{x^{-1}g \gm^{g^{-1}}}\\
& = & \Delta(g)\; \big( [\alpha_{g^{-1}}(a^*)]_{x^{-1}g \gm} * \gm g^{-1} \gm * 1_{\mathbf{s}(x^{-1}) \gm} \big) \,(g^{-1} \gm)\,.
\end{eqnarray*}

 Let us now prove equality (\ref{alt prod expression}). We have
\begin{eqnarray*}
 1_{\mathbf{r}(x) \gm} *  \gm g \gm *  [\alpha_{g^{-1}}(a)]_{xg\gm} \!\! & = & \!\! \Delta(g)\;\big([\alpha_{g^{-1}}(a^*)]_{x^{-1}g\gm} *  \gm g^{-1} \gm * 1_{\mathbf{r}(x) \gm} \big)^*\\
& = & \!\! \Delta(g)\; \big([\alpha_{g^{-1}}(a^*)]_{x^{-1}g\gm} *  \gm g^{-1} \gm * 1_{\mathbf{s}(x^{-1}g)g^{-1} \gm} \big)^*\,,
\end{eqnarray*}
which together with (\ref{invol for span elements expression}) yields
\begin{eqnarray*}
& = & \Delta(g) \Delta(g^{-1})\; [a]_{x\gm} *  \gm g \gm * 1_{\mathbf{s}(xg) \gm}\\
& = & [a]_{x\gm} *  \gm g \gm * 1_{\mathbf{s}(xg) \gm}\,.
\end{eqnarray*}

Let us now prove (\ref{a x gm times gm g gm expression}). An easy computation yields
\begin{eqnarray*}
[a]_{x \gm} * \gm g \gm\; (h \gm) & = & [a]_{x\gm} \cdot \gm g \gm (h \gm)\,,
\end{eqnarray*}
from which we conclude that $[a]_{x \gm} * \gm g \gm$ has support in $\gm g \gm$. Evaluating this element on the point $ g\gm$ we get
\begin{eqnarray*}
[a]_{x \gm} * \gm g \gm \;(g \gm) & = & [a]_{x\gm} \cdot \gm g \gm (g \gm) \;\; = \;\; [a]_{x\gm}\,.
\end{eqnarray*}
From Proposition \ref{axH and sum axhK inside MCcA proposition} one always has the following decomposition
\begin{align*}
 [a]_{x \gm} \;=\; \sum_{[\gamma] \in \mathcal{S}_x \backslash \gm / \gm^g} [\alpha_{\gamma^{-1}}(a)]_{x \gamma \gm^g}\,.
\end{align*}
Together with Lemma \ref{a_xH HgH chi_s(x)gHg} we get
\begin{eqnarray*}
 [a]_{x \gm} * \gm g \gm\; (g \gm) & = & [a]_{x \gm}\\
 & = & \sum_{[\gamma] \in \mathcal{S}_x \backslash \gm / \gm^g} [\alpha_{\gamma^{-1}}(a)]_{x \gamma \gm^g}\\
& = & \sum_{[\gamma] \in \mathcal{S}_x \backslash \gm / \gm^g} [\alpha_{\gamma^{-1}}(a)]_{x \gamma \gm} *\gm g \gm * 1_{\mathbf{s}(x) \gamma g \gm}\;(g \gm)\\
& = & \sum_{[\gamma] \in \mathcal{S}_x \backslash \gm / \gm^g} [a]_{x \gm} *\gm g \gm * 1_{\mathbf{s}(x) \gamma g \gm}\;(g \gm)\,,
\end{eqnarray*}
and equality (\ref{a x gm times gm g gm expression}) is proven.

Equality (\ref{gm g gm times a x gm expression}) follows easily from (\ref{a x gm times gm g gm expression}) by taking the involution and using the fact that $\mathcal{S}_x = \mathcal{S}_{x^{-1}}$.

The last claim of this proposition  follows simply from (\ref{alt prod expression}) and Proposition \ref{prop decomp of f in crossed prod}. \qed\\

In the theory of crossed products $A \times G$ by discrete groups one has a ``covariance relation'' of the form $g* a *g^{-1} = \alpha_g(a)$. This relation is essential in the passage from covariant representations of the system $(A, G, \alpha)$ to representations of the crossed product. More generally, the following relation holds in $A \times G$:
\begin{align*}
 g * a *h \;=\; \alpha_{g}(a) * gh\,.
\end{align*}
We will now explore how this generalizes to the setting of crossed products by Hecke pairs. What we are aiming for is a description of how products of the form $\gm g \gm * [a]_{x\gm} * \gm s \gm$ can be expressed by the canonical spanning set of elements of the form $[b]_{y\gm}* \gm h \gm * 1_{\mathbf{s}(x)h\gm}$ (according to Theorem \ref{prop decomp of f in crossed prod}). This will be achieved in Corollary \ref{cor prod in corossed prod} below and will play an important role in the representation theory of crossed products by Hecke pairs, particularly in the definition of covariant representations. One should observe the similarities between the expressions we obtain both in Theorem \ref{product formula} and Corollary \ref{cor prod in corossed prod} and the expression provided by an Huef, Kaliszewski and Raeburn in \cite[Definition 1.1]{cov} (if one ``forgets'' the representations in their setting). Once again, this is more than a coincidence as we will see in in the sequel to this article. In fact, an Huef, Kaliszewski and Raeburn's definition served as a guiding line for our results below and for the definition of a covariant representation (Definition \ref{covariant pre rep def}) which we shall present in the next section.

Before we establish the results we are aiming for we need to establish some notation, which will be used throughout this work. For $w, v \in G$ and a unit $y \in X^0$ we define the sets
\begin{align}
\mathfrak{n}_{w, v}^y\index{n(frak)wvy@$\mathfrak{n}_{w, v}^y$} & := \big\{[r] \in \gm w \gm / \gm: r^{-1}w v\gm \subseteq \gm v\gm\; \text{and} \;\; y  w^{-1} \in y \gm r^{-1}\big\} \,,\\
\mathfrak{d}_{w, v}^y\index{d2 (frac)wvy@$\mathfrak{d}_{w, v}^y$} & := \big\{[r] \in \gm w \gm / \gm: r^{-1} wv\gm \subseteq \gm v\gm\; \text{and} \;\; y  w^{-1} \in y \gm r^{-1} \gm^{wv}\big\}  \,.
\end{align}
and the numbers
\begin{align}
\label{number n}
n_{w, v}^y\index{n2 wvy@$n_{w, v}^y$} & := \#\; \mathfrak{n}_{w, v}^y\,,\\
\label{number d}
d_{w, v}^y\index{d2 wvy@$d_{w, v}^y$} & := \#\;  \mathfrak{d}_{w, v}^y \,,\\
\label{number N}
N_{w,v}^y\index{N3 wvy@$N_{w, v}^y$} & := \frac{n_{w, v}^y}{d_{w, v}^y}\,.
\end{align}
We will also denote by $E_{u,v}^y$\index{Euvy@$E_{u,v}^y$} the double coset space
\begin{align}
\label{double coset space space Euv y}
E_{u, v}^y := \mathcal{S}_y \backslash \gm / ( u \gm u^{-1} \cap v\gm v^{-1} )\,.
\end{align}\\

\begin{thm}
\label{product formula}
 Let $g,s \in G$ and $y \in X^0$. We have that
\begin{eqnarray*}
 \gm g \gm * 1_{y\gm} * \gm s \gm \!\! & = & \!\! \sum_{\substack{[w] \in \gm g \gm / \gm \\ [v] \in \gm s \gm / \gm}}\; \sum_{[\gamma] \in E_{w^{-1},v}^y} \frac{N_{w, v}^{y \gamma}}{L(wv)}\; \big(1_{y \gamma w^{-1} \gm} * \gm w v \gm * 1_{y \gamma v \gm}\big)\\
& = & \!\! \sum_{[v] \in \gm s \gm / \gm}\; \sum_{[\gamma] \in E_{g^{-1},v}^y} \frac{L(g)N_{g, v}^{y\gamma}}{L(gv)}\; \big(1_{y \gamma g^{-1} \gm} * \gm g v \gm * 1_{y \gamma v \gm}\big)\\
& = & \!\! \sum_{\substack{[u] \in \gm g^{-1} \gm / \gm \\ [v] \in \gm s \gm / \gm}}\; \sum_{[\gamma] \in E_{u,v}^y} \frac{\Delta(g)N_{u^{-1}, v}^{y\gamma}}{L(u^{-1}v)} \; \big(1_{y \gamma u \gm} * \gm u^{-1} v \gm * 1_{y \gamma v \gm}\big)\,.\\
\end{eqnarray*}
\end{thm}

In order to prove the above result we will need the following lemma, which gives some properties of the numbers $n_{w, v}^y$ and $d_{w, v}^y$.\\

\begin{lemma}
\label{lemma numbers n and d}
 Let $w, v, \in G$, $\theta \in \gm$ and $y \in X^0$. The numbers $n_{w, v}^y$ and $d_{w,v}^y$ satisfy the following properties:
\begin{itemize}
 \item[i)] $n_{w, v \theta}^y = n_{w, v}^y$ $\qquad\qquad\qquad\qquad\qquad\;\;$ $i')\;$ $d_{w, v \theta}^y = d_{w, v}^y$
 \item[ii)] $n_{\theta w, v}^y = n_{w, v}^y$ $\qquad\qquad\qquad\qquad\qquad\;$ $ii')\;$ $d_{\theta w, v}^y = d_{w, v}^y$
 \item[iii)] $n_{w, \theta^{-1}v}^{y \theta} = n_{w \theta^{-1}, v}^y$ $\qquad\qquad\qquad\qquad$ $iii')\;$ $d_{w, \theta^{-1}v}^{y \theta} = d_{w \theta^{-1}, v}^y$
\end{itemize}
More generally, if $\widetilde{w}, \widetilde{v} \in G$ and $\widetilde{y} \in X^0$ are such that $\gm \widetilde{w} \gm = \gm w \gm$, $\gm \widetilde{v} \gm = \gm v \gm$, $\widetilde{y}\gm = y \gm$, $\widetilde{w}\widetilde{v} \gm = wv\gm$ and $\widetilde{y} \widetilde{w}^{-1} \gm^{wv} = yw^{-1} \gm^{wv}$, then
\begin{itemize}
 \item[iv)] $n_{w, v}^y = n_{\widetilde{w}, \widetilde{v}}^{\widetilde{y}}$ $\qquad\qquad\qquad\qquad\qquad\;\;$ $iv')$ $d_{w, v}^y = d_{\widetilde{w}, \widetilde{v}}^{\widetilde{y}}$\\
\end{itemize}

\end{lemma}

{\bf \emph{Proof:}} Assertions $i)$ and $i')$ are obvious.

Assertion $ii)$ follows from the observation that $[r] \mapsto [\theta^{-1} r]$ establishes a bijection between the sets $\mathfrak{n}_{w, v}^y$ and $\mathfrak{n}_{\theta w, v}^y$.

Assertion $ii')$ is proven in a similar fashion as assertion $ii)$.

To prove assertion $iv)$, let $\theta \in \gm^{wv}$ be such that $\widetilde{y}\widetilde{w}^{-1} = yw^{-1}\theta$. We have
\begin{eqnarray*}
\mathfrak{n}_{\widetilde{w}, \widetilde{v}}^{\widetilde{y}} & = & \big\{[r] \in \gm \widetilde{w} \gm / \gm: r^{-1}\widetilde{w}  \widetilde{v}\gm \subseteq \gm \widetilde{v}\gm\; \text{and} \;\; \widetilde{y} \widetilde{w}^{-1} \in \widetilde{y} \gm r^{-1}\big\}\\
& = & \big\{[r] \in \gm w \gm / \gm: r^{-1}w  v\gm \subseteq \gm v\gm\; \text{and} \;\; y w^{-1} \theta \in y \gm r^{-1}\big\}\,.
\end{eqnarray*}
Since $\theta \in \gm^{wv}$ we have $\theta wv \gm = wv \gm$, so that
\begin{eqnarray*}
& = & \big\{[r] \in \gm \theta^{-1} w \gm / \gm: r^{-1}\theta^{-1}w  v\gm \subseteq \gm v\gm\; \text{and} \;\; y w^{-1} \theta \in y \gm r^{-1}\big\}\\
& = & \mathfrak{n}_{\theta^{-1}w, v}^{y}\,.
\end{eqnarray*}
Now, from assertion $ii)$, it follows that $n_{\widetilde{w}, \widetilde{v}}^{\widetilde{y}} = n_{\theta^{-1}w, v}^{y} = n_{w, v}^{y}$.

As for assertion $iv')$, taking $\theta \in \gm^{wv}$ again as such that $\widetilde{y}\widetilde{w}^{-1} = yw^{-1}\theta$, we notice that
\begin{eqnarray*}
\mathfrak{d}_{\widetilde{w}, \widetilde{v}}^{\widetilde{y}} & = & \big\{[r] \in \gm \widetilde{w} \gm / \gm: r^{-1}\widetilde{w}  \widetilde{v}\gm \subseteq \gm \widetilde{v}\gm\; \text{and} \;\; \widetilde{y} \widetilde{w}^{-1} \in \widetilde{y} \gm r^{-1} \gm^{\widetilde{w }\widetilde{v}}\big\}\\
& = & \big\{[r] \in \gm w \gm / \gm: r^{-1}w  v\gm \subseteq \gm v\gm\; \text{and} \;\; y w^{-1} \theta \in y \gm r^{-1}\gm^{wv}\big\}\\
& = & \big\{[r] \in \gm w \gm / \gm: r^{-1}w  v\gm \subseteq \gm v\gm\; \text{and} \;\; y w^{-1} \in y \gm r^{-1}\gm^{wv}\big\}\\
& = & \mathfrak{d}_{w, v}^{y}\,.
\end{eqnarray*}

Assertions $iii)$ and $iii')$ are a direct consequence of $iv)$ and $iv')$. \qed\\

{\bf \emph{Proof of Theorem \ref{product formula}:}} We have
\begin{eqnarray*}
\gm g \gm * 1_{y\gm} * \gm s \gm \;(t \gm) & = & \sum_{[w] \in G / \gm} \gm g \gm (w \gm)\, \overline{\alpha}_w \big((1_{y\gm} * \gm s \gm) \; (w^{-1} t \gm) \big)\\
& = & \sum_{[w] \in \gm g \gm / \gm} \overline{\alpha}_w \big((1_{y\gm} * \gm s \gm) \; (w^{-1} t \gm) \big)\\
& = & \sum_{[w] \in \gm g \gm / \gm} \overline{\alpha}_w \big(1_{y\gm} \cdot \gm s \gm (w^{-1} t \gm) \big)\\
& = & \sum_{\substack{[w] \in \gm g \gm / \gm \\  w^{-1} t \gm \subseteq \gm s\gm}} \overline{\alpha}_w (1_{y\gm} )\\
& = & \sum_{\substack{[w] \in \gm g \gm / \gm \\  w^{-1} t \gm \subseteq \gm s\gm}} 1_{y\gm w^{-1}}
\end{eqnarray*}

We now claim that

\begin{align}
\label{equality in proof of product}
 \sum_{\substack{[w] \in \gm g \gm / \gm \\  w^{-1} t \gm \subseteq \gm s\gm}} 1_{y\gm w^{-1}}\; = \; \sum_{\substack{[w] \in \gm g \gm / \gm \\  w^{-1} t \gm \subseteq \gm s\gm}} \sum_{[\gamma] \in E_{w^{-1},w^{-1}t}^y} N_{w,w^{-1}t}^{y \gamma}\, 1_{y\gamma w^{-1}\gm^t}\,.
\end{align}

To see this, we will evaluate both the right and left expressions above on all points $x \in X^0$ and see that we obtain the same value. First, we note that if $x \in X^0$ is not of the form $y \theta \widetilde{w}^{-1}$, for some $\theta \in \gm$ and $\widetilde{w} \in \gm g \gm$ such that $\widetilde{w}^{-1} t \gm \subseteq \gm s\gm$, then both expressions are zero. Suppose now that $x = y \theta \widetilde{w}^{-1}$ for some $\widetilde{w} \in \gm g \gm$ such that $\widetilde{w}^{-1} t \gm \subseteq \gm s\gm$. Evaluating the left expression we get
\begin{align*}
  \sum_{\substack{[w] \in \gm g \gm / \gm \\  w^{-1} t \gm \subseteq \gm s\gm}} 1_{y\gm w^{-1}} (y \theta \widetilde{w}^{-1}) = \sum_{\substack{[w] \in \gm \widetilde{w} \gm / \gm \\  w^{-1} \widetilde{w} \widetilde{w}^{-1} t \gm \subseteq \gm \widetilde{w}^{-1} t\gm}} 1_{y\gm w^{-1}} (y \theta \widetilde{w}^{-1}) = n_{\widetilde{w}, \widetilde{w}^{-1}t}^{y\theta}\,.
\end{align*}
As for the right expression, first we observe that if $y \theta \widetilde{w}^{-1} \in y \gamma w^{-1} \gm^t$, then by Lemma \ref{lemma numbers n and d} $iv)$ and $iv')$ we have $N_{\widetilde{w}, \widetilde{w}^{-1}t }^{y \theta} = N_{w, w^{-1}t}^{y \gamma}$. Thus, evaluating the right expression we get
\begin{eqnarray*}
& & \sum_{\substack{[w] \in \gm g \gm / \gm \\  w^{-1} t \gm \subseteq \gm s\gm}} \sum_{[\gamma] \in E_{w^{-1},w^{-1}t}^y} N_{w,w^{-1}t}^{y\gamma}\, 1_{y\gamma w^{-1}\gm^t} \,(y \theta \widetilde{w}^{-1})  = \\
&  = & \sum_{\substack{[w] \in \gm g \gm / \gm \\  w^{-1} t \gm \subseteq \gm s\gm}} \sum_{[\gamma] \in E_{w^{-1},w^{-1}t}^y} N_{\widetilde{w},\widetilde{w}^{-1}t}^{y \theta}\, 1_{y\gamma w^{-1}\gm^t} \,(y \theta \widetilde{w}^{-1})\\
&  = & N_{\widetilde{w},\widetilde{w}^{-1}t}^{y \theta} \sum_{\substack{[w] \in \gm g \gm / \gm \\  w^{-1} t \gm \subseteq \gm s\gm}} \sum_{[\gamma] \in E_{w^{-1},w^{-1}t}^y} \, 1_{y\gamma w^{-1}\gm^t} \,(y \theta \widetilde{w}^{-1})
\end{eqnarray*}
Using Proposition \ref{xKH bijection with Sx K H} we notice that
\begin{eqnarray*}
 \sum_{[\gamma] \in E_{w^{-1},w^{-1}t}^y} \, 1_{y\gamma w^{-1}\gm^t} & = & \sum_{[\gamma] \in E_{w^{-1},w^{-1}t}^y} \, 1_{y\gamma (w^{-1}\gm w \cap w^{-1}t \gm t^{-1} w) w^{-1}}\\ 
& = & 1_{y\gm w^{-1}\gm^t} \,,
\end{eqnarray*}
from which we obtain that,
\begin{eqnarray*}
& & N_{\widetilde{w},\widetilde{w}^{-1}t}^{y \theta} \sum_{\substack{[w] \in \gm g \gm / \gm \\  w^{-1} t \gm \subseteq \gm s\gm}} \sum_{[\gamma] \in E_{w^{-1},w^{-1}t}^y} \, 1_{y\gamma w^{-1}\gm^t} \,(y \theta \widetilde{w}^{-1}) \;\;  = \\
& = &  N_{\widetilde{w},\widetilde{w}^{-1}t}^{y \theta} \sum_{\substack{[w] \in \gm g \gm / \gm \\  w^{-1} t \gm \subseteq \gm s\gm}} 1_{y\gm w^{-1}\gm^t} \,(y \theta \widetilde{w}^{-1})\\
& = &  N_{\widetilde{w},\widetilde{w}^{-1}t}^{y \theta} \sum_{\substack{[w] \in \gm \widetilde{w} \gm / \gm \\  w^{-1} \widetilde{w} \widetilde{w}^{-1} t \gm \subseteq \gm \widetilde{w}^{-1}t\gm}} 1_{y\gm w^{-1}\gm^t} \,(y \theta \widetilde{w}^{-1})\\
& = &  N_{\widetilde{w},\widetilde{w}^{-1}t}^{y \theta}\; d_{\widetilde{w}, \widetilde{w}^{-1}t}^{y \theta}\\
& = &  n_{\widetilde{w}, \widetilde{w}^{-1}t}^{y \theta}\,.
\end{eqnarray*}

So, equality (\ref{equality in proof of product}) is established.

Now, by Proposition \ref{a_xH HgH chi_s(x)gHg}, we see that
\begin{align*}
 \sum_{\substack{[w] \in \gm g \gm / \gm \\  w^{-1} t \gm \subseteq \gm s\gm}} \sum_{[\gamma] \in E_{w^{-1},w^{-1}t}^y} N_{w,w^{-1}t}^{y\gamma}\, 1_{y\gamma w^{-1}\gm^t} = \qquad\qquad\qquad\qquad\\
 = \sum_{\substack{[w] \in \gm g \gm / \gm \\  w^{-1} t \gm \subseteq \gm s\gm}} \sum_{[\gamma] \in E_{w^{-1},w^{-1}t}^y}  N_{w,w^{-1}t}^{y\gamma}\; \big(1_{y \gamma w^{-1} \gm} * \gm t \gm * 1_{y \gamma w^{-1}t \gm}\big) (t\gm)
\end{align*}
Now, using the fact that condition $w^{-1} t \gm \subseteq \gm s \gm$ means that there exists a (necessarily unique) element $[v] \in \gm s \gm / \gm$ such that $w^{-1}t \gm = v \gm$, or equivalently, $t \gm = wv\gm$, we obtain
\begin{eqnarray*}
 = \sum_{\substack{[w] \in \gm g \gm / \gm \\ [v] \in \gm s \gm / \gm \\  wv\gm = t \gm}} \sum_{[\gamma] \in E_{w^{-1},w^{-1}t}^y}  N_{w,w^{-1}t}^{y\gamma}\; \big(1_{y \gamma w^{-1} \gm} * \gm t \gm * 1_{y \gamma w^{-1}t \gm}\big) (t\gm)\\
 = \sum_{\substack{[w] \in \gm g \gm / \gm \\ [v] \in \gm s \gm / \gm \\  wv\gm = t \gm}} \sum_{[\gamma] \in E_{w^{-1},v}^y}  N_{w,v}^{y\gamma}\; \big(1_{y \gamma w^{-1} \gm} * \gm wv \gm * 1_{y \gamma v \gm}\big) (t\gm)\,.
\end{eqnarray*}

We now claim that
\begin{align*}
\sum_{\substack{[w] \in \gm g \gm / \gm \\ [v] \in \gm s \gm / \gm \\  wv\gm = t \gm}} \sum_{[\gamma] \in E_{w^{-1},v}^y}  N_{w,v}^{y\gamma}\; \big(1_{y \gamma w^{-1} \gm} * \gm wv \gm * 1_{y \gamma v \gm}\big) (t\gm) =\\
= \sum_{\substack{[w] \in \gm g \gm / \gm \\ [v] \in \gm s \gm / \gm}}\; \sum_{[\gamma] \in E_{w^{-1},v}^y} \frac{N_{w,v}^{y\gamma}}{L(wv)}\; \big(1_{y \gamma w^{-1} \gm} * \gm w v \gm * 1_{y \gamma v \gm}\big) (t\gm)
\end{align*}
To prove this we note that, given any $[w] \in \gm g \gm / \gm$ and $[v] \in \gm s \gm / \gm$, the element $\big(1_{y \gamma w^{-1} \gm} * \gm w v \gm * 1_{y \gamma v \gm}\big) (t\gm)$ is nonzero if and only if $\gm t \gm = \gm wv\gm$, so that we can write
\begin{eqnarray*}
& &\!\!\!\! \sum_{\substack{[w] \in \gm g \gm / \gm \\ [v] \in \gm s \gm / \gm}}\; \sum_{[\gamma] \in E_{w^{-1},v}^y} \frac{N_{w,v}^{y\gamma}}{L(wv)}\; \big(1_{y \gamma w^{-1} \gm} * \gm w v \gm * 1_{y \gamma v \gm}\big) (t\gm) =\\
& = &\!\!\!\! \sum_{\substack{[w] \in \gm g \gm / \gm \\ [v] \in \gm s \gm / \gm \\ wv \gm \subseteq \gm t \gm}}\; \sum_{[\gamma] \in E_{w^{-1},v}^y} \frac{N_{w,v}^{y\gamma}}{L(wv)}\; \big(1_{y \gamma w^{-1} \gm} * \gm w v \gm * 1_{y \gamma v \gm}\big) (t\gm) \\
& = &\!\!\!\! \sum_{[\theta] \in \gm / \gm^t} \sum_{\substack{[w] \in \gm g \gm / \gm \\ [v] \in \gm s \gm / \gm \\ wv \gm = \theta t \gm}}\; \sum_{[\gamma] \in E_{w^{-1},v}^y} \frac{N_{w,v}^{y\gamma}}{L(wv)}\; \big(1_{y \gamma w^{-1} \gm} * \gm w v \gm * 1_{y \gamma v \gm}\big) (t\gm) \\
& = &\!\!\!\! \sum_{[\theta] \in \gm / \gm^t} \sum_{\substack{[w] \in \gm g \gm / \gm \\ [v] \in \gm s \gm / \gm \\ \theta wv \gm = \theta t \gm}}\; \sum_{[\gamma] \in E_{ w^{-1} \theta^{-1},v}^y} \!\! \frac{N_{\theta w, v}^{y \gamma}}{L( \theta wv)}\; \big(1_{y \gamma w^{-1} \theta^{-1} \gm} * \gm \theta w v \gm * 1_{y \gamma v \gm}\big) (t\gm) \\
& = &\!\!\!\! \sum_{[\theta] \in \gm / \gm^t} \sum_{\substack{[w] \in \gm g \gm / \gm \\ [v] \in \gm s \gm / \gm \\ wv \gm = t \gm}}\; \sum_{[\gamma] \in E_{w^{-1},v}^y} \frac{N_{\theta w, v}^{y \gamma}}{L(wv)}\; \big(1_{y \gamma w^{-1} \gm} * \gm w v \gm * 1_{y \gamma v \gm}\big) (t\gm)
\end{eqnarray*}
By Lemma \ref{lemma numbers n and d} $ii)$ and $ii')$ we know that $N_{\theta w, v}^{y \gamma} = N_{ w, v}^{y \gamma}$, hence

\begin{eqnarray*}
& = & \sum_{[\theta] \in \gm / \gm^t} \sum_{\substack{[w] \in \gm g \gm / \gm \\ [v] \in \gm s \gm / \gm \\ wv \gm = t \gm}}\; \sum_{[\gamma] \in E_{w^{-1},v}^y} \frac{N_{ w, v}^{y \gamma}}{L(wv)}\; \big(1_{y \gamma w^{-1} \gm} * \gm w v \gm * 1_{y \gamma v \gm}\big) (t\gm)\\
& = & L(t) \sum_{\substack{[w] \in \gm g \gm / \gm \\ [v] \in \gm s \gm / \gm \\ wv \gm = t \gm}}\; \sum_{[\gamma] \in E_{w^{-1},v}^y} \frac{N_{ w, v}^{y \gamma}}{L(wv)}\; \big(1_{y \gamma w^{-1} \gm} * \gm w v \gm * 1_{y \gamma v \gm}\big) (t\gm)\\
& = & \sum_{\substack{[w] \in \gm g \gm / \gm \\ [v] \in \gm s \gm / \gm \\ wv \gm = t \gm}}\; \sum_{[\gamma] \in E_{w^{-1},v}^y} N_{ w, v}^{y \gamma}\; \big(1_{y \gamma w^{-1} \gm} * \gm w v \gm * 1_{y \gamma v \gm}\big) (t\gm)\,.
\end{eqnarray*}

Hence, we have proven that
\begin{align*}
 \gm g \gm * 1_{y\gm} * \gm s \gm & = & \sum_{\substack{[w] \in \gm g \gm / \gm \\ [v] \in \gm s \gm / \gm}}\; \sum_{[\gamma] \in E_{w^{-1},v}^y} \frac{N_{w, v}^{y \gamma}}{L(wv)}\; \big(1_{y \gamma w^{-1} \gm} * \gm w v \gm * 1_{y \gamma v \gm}\big)\,.
\end{align*}
Also,
\begin{eqnarray*}
 &  & \sum_{\substack{[w] \in \gm g \gm / \gm \\ [v] \in \gm s \gm / \gm}}\; \sum_{[\gamma] \in E_{w^{-1},v}^y} \frac{N_{w, v}^{y \gamma}}{L(wv)}\; \big(1_{y \gamma w^{-1} \gm} * \gm w v \gm * 1_{y \gamma v \gm}\big)\\
 & = & \sum_{\substack{[\theta] \in \gm  / \gm^g \\ [v] \in \gm s \gm / \gm}}\; \sum_{[\gamma] \in E_{ g^{-1} \theta^{-1},v}^y} \frac{N_{\theta g, v}^{y \gamma}}{L(\theta gv)}\; \big(1_{y \gamma g^{-1} \theta^{-1} \gm} * \gm \theta g v \gm * 1_{y \gamma v \gm}\big)\\
& = & \sum_{\substack{[\theta] \in \gm / \gm^g \\ [v] \in \gm s \gm / \gm}}\; \sum_{[\gamma] \in E_{g^{-1},v}^y} \frac{N_{g, v}^{y \gamma}}{L(gv)}\; \big(1_{y \gamma g^{-1} \gm} * \gm g v \gm * 1_{y \gamma v \gm}\big)\\
& = & \sum_{[v] \in \gm s \gm / \gm}\; \sum_{[\gamma] \in E_{g^{-1},v}^y} \frac{L(g)N_{g, v}^{y \gamma}}{L(gv)}\; \big(1_{y \gamma g^{-1} \gm} * \gm g v \gm * 1_{y \gamma v \gm}\big)\,.
\end{eqnarray*}
Moreover, we also have
\begin{eqnarray*}
& & \sum_{[v] \in \gm s \gm / \gm}\; \sum_{[\gamma] \in E_{g^{-1},v}^y} \frac{L(g)N_{g, v}^{y \gamma}}{L(gv)}\; \big(1_{y \gamma g^{-1} \gm} * \gm g v \gm * 1_{y \gamma v \gm}\big)\\
& = & L(g^{-1}) \sum_{[v] \in \gm s \gm / \gm}\; \sum_{[\gamma] \in E_{g^{-1},v}^y} \frac{\Delta(g)N_{g, v}^{y \gamma}}{L(gv)}\; \big(1_{y \gamma g^{-1} \gm} * \gm g v \gm * 1_{y \gamma v \gm}\big)\\
& = & \sum_{\substack{[\theta] \in \gm / \gm^{g^{-1}}\\ [v] \in \gm s \gm / \gm}}\; \sum_{[\gamma] \in E_{g^{-1},v}^y} \frac{\Delta(g)N_{g, v}^{y \gamma}}{L(gv)}\; \big(1_{y \gamma g^{-1} \gm} * \gm g v \gm * 1_{y \gamma v \gm}\big)\\
& = & \sum_{\substack{[\theta] \in \gm / \gm^{g^{-1}}\\ [v] \in \gm s \gm / \gm}}\; \sum_{[\gamma] \in E_{ g^{-1},\theta^{-1}v}^y} \frac{\Delta(g)N_{g ,  \theta^{-1}v}^{y \gamma}}{L(g\theta^{-1} v)}\; \big(1_{y \gamma g^{-1} \gm} * \gm g \theta^{-1} v \gm * 1_{y \gamma  \theta^{-1} v \gm}\big)\,,
\end{eqnarray*}
but since there is a well-defined bijection $ E_{\theta g^{-1}, v}^y \rightarrow E_{g^{-1}, \theta^{-1} v}^y$ given by
 $[\gamma] \mapsto [\gamma \theta]$, we obtain
\begin{eqnarray*}
& = & \sum_{\substack{[\theta] \in \gm / \gm^{g^{-1}}\\ [v] \in \gm s \gm / \gm}}\; \sum_{[\gamma] \in E_{ \theta g^{-1},v}^y} \frac{\Delta(g)N_{g , \theta^{-1} v}^{y \gamma \theta}}{L(g\theta^{-1} v)}\; \big(1_{y \gamma \theta g^{-1} \gm} * \gm g \theta^{-1} v \gm * 1_{y \gamma \theta  \theta^{-1} v \gm}\big)\\
\end{eqnarray*}
and from Lemma \ref{lemma numbers n and d} we get $N_{g , \theta^{-1} v}^{y \gamma \theta} = N_{g \theta^{-1} ,  v}^{y \gamma}$, thus
\begin{eqnarray*}
& = & \sum_{\substack{[\theta] \in \gm / \gm^{g^{-1}}\\ [v] \in \gm s \gm / \gm}}\; \sum_{[\gamma] \in E_{ \theta g^{-1},v}^y} \frac{\Delta(g)N_{g \theta^{-1} ,  v}^{y \gamma}}{L(g\theta^{-1} v)}\; \big(1_{y \gamma \theta g^{-1} \gm} * \gm g \theta^{-1} v \gm * 1_{y \gamma  v \gm}\big)\\
& = & \sum_{\substack{[u] \in \gm g^{-1} \gm / \gm\\ [v] \in \gm s \gm / \gm}}\; \sum_{[\gamma] \in E_{ u,v}^y} \frac{\Delta(g)N_{u^{-1},v}^{y \gamma}}{L(u^{-1} v)}\; \big(1_{y \gamma u \gm} * \gm u^{-1} v \gm * 1_{y \gamma  v \gm}\big)\,.
\end{eqnarray*}\qed\\

\begin{cor}
\label{cor prod in corossed prod}
Similarly, for $a \in \a_x$ with $x \in X$, we have
\begin{eqnarray*}
 & & \gm g \gm * [a]_{x\gm} * \gm s \gm \;\; = \;\;\\
 & = & \sum_{\substack{[w] \in \gm g \gm / \gm \\ [v] \in \gm s \gm / \gm}}\; \sum_{[\gamma] \in E_{w^{-1},v}^{\mathbf{s}(x)}} \frac{N_{w, v}^{\mathbf{s}(x) \gamma}}{L(wv)}\; \big([\alpha_{w\gamma^{-1}}(a)]_{x \gamma w^{-1} \gm} * \gm w v \gm * 1_{\mathbf{s}(x) \gamma v \gm}\big)\\
& = & \sum_{[v] \in \gm s \gm / \gm}\; \sum_{[\gamma] \in E_{g^{-1},v}^{\mathbf{s}(x)}} \frac{L(g)N_{g, v}^{\mathbf{s}(x)\gamma}}{L(gv)}\; \big([\alpha_{g\gamma^{-1}}(a)]_{x \gamma g^{-1} \gm} * \gm g v \gm * 1_{\mathbf{s}(x) \gamma v \gm}\big)\\
& = & \sum_{\substack{[u] \in \gm g^{-1} \gm / \gm \\ [v] \in \gm s \gm / \gm}}\; \sum_{[\gamma] \in E_{u,v}^{\mathbf{s}(x)}} \frac{\Delta(g)N_{u^{-1}, v}^{\mathbf{s}(x)\gamma}}{L(u^{-1}v)} \; \big([\alpha_{u^{-1}\gamma^{-1}}(a)]_{x \gamma u \gm} * \gm u^{-1} v \gm * 1_{\mathbf{s}(x) \gamma v \gm}\big)\,.
\end{eqnarray*}\\
\end{cor}

{\bf \emph{Proof:}} According to equality (\ref{gm g gm times a x gm expression}) in Proposition \ref{invol in crossed prod} we have
\begin{eqnarray*}
& & \gm g \gm * [a]_{x\gm} * \gm s \gm \;\; = \;\;\\
 & = & \sum_{[\theta] \in \mathcal{S}_x \backslash \gm / \gm^{g^{-1}}} 1_{\mathbf{r}(x) \theta g^{-1} \gm} * \gm g \gm * [a]_{x \gm} * \gm s \gm\\
 & = & \sum_{[\theta] \in \mathcal{S}_x \backslash \gm / \gm^{g^{-1}}} 1_{\mathbf{r}(x) \theta g^{-1} \gm} * \gm g \gm * [\alpha_{g^{-1}}(\alpha_{g\theta^{-1}}(a))]_{x \theta g^{-1}g\gm} * \gm s \gm
\end{eqnarray*}
and by (\ref{alt prod expression}) in the same proposition we get
\begin{eqnarray*}
 & = & \sum_{[\theta] \in \mathcal{S}_x \backslash \gm / \gm^{g^{-1}}} [\alpha_{g \theta^{-1}}(a)]_{x \theta g^{-1} \gm} * \gm g \gm * 1_{\mathbf{s}(x) \gm} * \gm s \gm\,,
\end{eqnarray*}
and by Theorem \ref{product formula} we obtain
\begin{eqnarray*}
 & = & \sum_{\substack{[\theta] \in \mathcal{S}_x \backslash \gm / \gm^{g^{-1}}\\ [w] \in \gm g \gm / \gm \\ [v] \in \gm s \gm / \gm \\ [\gamma] \in E_{w^{-1},v}^{\mathbf{s}(x)}}} \frac{N_{w, v}^{\mathbf{s}(x) \gamma}}{L(wv)}\; [\alpha_{g \theta^{-1}}(a)]_{x \theta g^{-1} \gm} * 1_{\mathbf{s}(x) \gamma w^{-1} \gm} * \gm w v \gm * 1_{\mathbf{s}(x) \gamma v \gm}\,.
\end{eqnarray*}
For each fixed $w$, $v$ and $\gamma$ all the summands in the expression 
\begin{align*}
 \sum_{[\theta] \in \mathcal{S}_x \backslash \gm / \gm^{g^{-1}}} \frac{N_{w, v}^{\mathbf{s}(x) \gamma}}{L(wv)}\; [\alpha_{g \theta^{-1}}(a)]_{x \theta g^{-1} \gm} * 1_{\mathbf{s}(x) \gamma w^{-1} \gm} * \gm w v \gm * 1_{\mathbf{s}(x) \gamma v \gm}\,,
\end{align*}
are zero except precisely for one summand and we have
\begin{eqnarray*}
&  & \sum_{[\theta] \in \mathcal{S}_x \backslash \gm / \gm^{g^{-1}}} \frac{N_{w, v}^{\mathbf{s}(x) \gamma}}{L(wv)}\; [\alpha_{g \theta^{-1}}(a)]_{x \theta g^{-1} \gm} * 1_{\mathbf{s}(x) \gamma w^{-1} \gm} * \gm w v \gm * 1_{\mathbf{s}(x) \gamma v \gm}\\
& = & \frac{N_{w, v}^{\mathbf{s}(x) \gamma}}{L(wv)}\; [\alpha_{w \gamma^{-1}}(a)]_{x \gamma w^{-1} \gm} * \gm w v \gm * 1_{\mathbf{s}(x) \gamma v \gm}\,.
\end{eqnarray*}
Hence we obtain
\begin{eqnarray*}
& & \gm g \gm * a_{x\gm} * \gm s \gm\;\; = \;\;\\
 & = & \sum_{\substack{[w] \in \gm g \gm / \gm \\ [v] \in \gm s \gm / \gm}}\; \sum_{[\gamma] \in E_{w^{-1},v}^{\mathbf{s}(x)}} \frac{N_{w, v}^{\mathbf{s}(x) \gamma}}{L(wv)}\; [\alpha_{w \gamma^{-1}}(a)]_{x \gamma w^{-1} \gm} * \gm w v \gm * 1_{\mathbf{s}(x) \gamma v \gm}\,.
\end{eqnarray*}
The remaining equalities in the statement of this corollary are proven in a similar fashion. \qed\\

\subsection{Basic Examples}

\begin{ex}
We will now show that when $\gm$ is a normal subgroup of $G$ our notion of a crossed product by the Hecke pair $(G,\gm)$ is precisely the usual crossed product by the group $G / \gm$. Normality of the subgroup $\gm$ implies that the $G$-action $\overline{\alpha}$ on $M(C_c(\a))$ gives rise to an action of $G / \gm$ on $C_c(\a / \gm)$. Moreover, we have $\gm^g = \gm$ for all $g \in G$, and it follows easily from the definitions that $C_c(\a / \gm) \times^{alg} G / \gm$ is nothing but the usual crossed product by the action of the group $G / \gm$.

It is also interesting to observe that any usual crossed product $C_c(\b) \times^{alg} G / \gm$ coming from an action of the group $G / \gm$ on a Fell bundle $\b$ over a groupoid $Y$ is actually a crossed by the Hecke pair $(G, \gm)$ in our sense. To see this we note that the action of $G / \gm$ on $\b$ lifts to an action of $G$ on $\b$. In this lifted action the subgroup $\gm$ acts trivially, so that the action is $\gm$-good. Moreover, since $\gm$ is normal in $G$, the $\gm$-intersection property is also trivially satisfied. It is clear that $Y / \gm $ is just $Y$ and $\b / \gm$ coincides with $\b$. Thus,  forming the crossed product by the Hecke pair $(G, \gm)$ will give nothing but the usual crossed product by $G / \gm$, i.e. $C_c(\b / \gm) \times^{alg} G / \gm \cong C_c(\b) \times^{alg} G / \gm$.\\
\end{ex}

\begin{ex}
\label{Hecke algebra as an algebraic crossed product example}
 We will now explain how the Hecke algebra $\h(G, \gm)$ is an example of a crossed product by a Hecke pair, namely $\h(G, \gm) \cong \mathbb{C} \times^{alg} G / \gm$, just like group algebras are examples of crossed products by groups.

We start with a groupoid $X$ consisting of only one element, i.e. $X = \{*\}$, and we take $\mathbb{C}$ as the Fell $\a$ bundle over $X$, i.e. $\a_* = \mathbb{C}$. We take also the trivial $G$-action $\alpha$ on $\a$. Since the $G$-action fixes every element of $\a$, it is indeed $\gm$-good and in this case we have $X / \gm = X = \{*\}$. For the orbit bundle we have that $\a / \gm = \a$, and moreover
\begin{align*}
 C_c(\a / \gm) \cong C_c(X / \gm) \cong C_c(X) \cong \mathbb{C}\,.
\end{align*}
Hence, we are in the conditions of the Standing Assumption \ref{standing assumption} and we can form the crossed product $C_c(\a / \gm) \times^{alg}_{\alpha} G / \gm$, which we will simply write as $\mathbb{C} \times^{alg}_{\alpha} G / \gm$.

Since $\mathbb{C}$ is unital the definitions of $B(\a, G, \gm)$ and $\mathbb{C} \times^{alg}_{\alpha} G / \gm$ coincide in this case. Moreover Definition \ref{def of crossed product by hecke pair} reads that $\mathbb{C} \times^{alg}_{\alpha} G / \gm$ is the set of functions $f: G /\gm \to \mathbb{C}$ satisfying the compatibility condition (\ref{compatibility cond}). Since the action $\overline{\alpha}$ is trivial, the compatibility condition simply says that $\mathbb{C} \times^{alg}_{\alpha} G / \gm$ consists of all the functions $f: G /\gm \to \mathbb{C}$ which are left $\gm$-invariant. Morever, the product and involution expressions become respectively
\begin{align*}
 (f_1 * f_2)(g \gm) & := \sum_{[h] \in G/ \gm} f_1(h \gm)\, f_2(h^{-1} g\gm)\,,\\
(f^*)\,(g \gm) & := \Delta(g^{-1})\;\overline{f(g^{-1} \gm)}\,.
\end{align*}
Hence, it is clear that $\mathbb{C} \times^{alg}_{\alpha} G / \gm$ is nothing but the Hecke algebra $\h(G, \gm)$.

It follows from this that the product $\gm g \gm * 1_{*\gm}* \gm s \gm$ is just the product of the double cosets $\gm g \gm$ and $\gm s \gm$ inside the Hecke algebra, since $1_{*\gm}$ is the identity element. It is interesting to note in this regard that the expression for this product described in Theorem \ref{product formula} is a familiar expression for the product $\gm g \gm * \gm s \gm$ in $\h(G, \gm)$. To see this, we note that the stabilizer $\mathcal{S}_{*}$ of $*$ is the whole group $G$, and therefore $E^*_{u,v}$ consists only of the class $[e]$. Moreover, the numbers $n_{u^{-1}, v}^*$ and $d_{u^{-1},v}^*$, defined in (\ref{number n}) and (\ref{number d}), are equal, so that $N^*_{u^{-1}, v} = 1$. Thus, the expression described in Theorem  \ref{product formula} is just the usual expression
\begin{align*}
 \gm g \gm * \gm s \gm = \sum_{\substack{[u] \in \gm g^{-1} \gm / \gm \\ [v] \in \gm s \gm / \gm}} \frac{\Delta(g)}{L(u^{-1}v)} \;  \gm u^{-1} v \gm \,.
\end{align*}\\
\end{ex}

\begin{ex}
\label{action fixes all points tensor product example}
As a generalization of Example \ref{Hecke algebra as an algebraic crossed product example} we will now show that if the $G$-action fixes every element of the bundle $\a$, then $C_c(\a / \gm) \times^{alg}_{\alpha} G / \gm$ is isomorphic to the $^*$-algebraic tensor product of $C_c(\a/ \gm)$ and $\h(G, \gm)$. This result also has a known analogue in the theory of crossed products by groups.\\

\begin{prop}
If the $G$-action fixes every element of $\a$, then we have
\begin{align*}
 C_c(\a / \gm) \times^{alg}_{\alpha} G / \gm \;\; \cong \;\; C_c(\a / \gm) \odot \h(G, \gm)\,,
\end{align*}
where $\odot$ is the symbol that denotes the $^*$-algebraic tensor product.\\
\end{prop}

{\bf \emph{Proof:}} Given that we have canonical embeddings of $C_c(\a / \gm)$ and $\h(G, \gm)$ into $M(C_c(\a / \gm) \times^{alg}_{\alpha} G / \gm)$ we have a canonical linear map from $C_c(\a / \gm) \odot \h(G, \gm)$ to $M(C_c(\a / \gm) \times^{alg}_{\alpha} G / \gm)$ determined by
\begin{align}
\label{mapping from tensor product to crossed prod}
 f_1 \otimes f_2 \mapsto f_1*f_2\,,
\end{align}
where $f_1 \in C_c(\a / \gm)$ and $f_2 \in \h(G, \gm)$. Standard arguments can be used to show that this mapping is injective (since the mappings from both $C_c(\a / \gm)$ and $\h(G, \gm)$ into the multiplier algebra of the crossed product are injections). It is also clear that the image of the map determined by (\ref{mapping from tensor product to crossed prod}) is contained in $C_c(\a / \gm) \times^{alg}_{\alpha} G / \gm$. Let us now check that this mapping is surjective. First we will show that the elements of $C_c(\a / \gm)$ commute with elements of $\h(G, \gm)$ inside $M(C_c(\a / \gm) \times^{alg}_{\alpha} G / \gm)$. It follows from expressions (\ref{a x gm times gm g gm expression}) and (\ref{alt prod expression}) that
\begin{eqnarray*}
 [a]_{x\gm} * \gm g \gm & = & \sum_{[\gamma] \in \mathcal{S}_x \backslash \gm / \gm^g } [a]_{x\gm} * \gm g \gm * 1_{\mathbf{s}(x) \gamma g \gm}\\
& = & \sum_{[\gamma] \in \mathcal{S}_x \backslash \gm / \gm^g } 1_{\mathbf{r}(x)   \gm} * \gm g \gm * [\alpha_{g^{-1} \gamma^{-1}}(a)]_{x \gamma g \gm}\,.
\end{eqnarray*}
Since every point of $X$ is fixed by the associated $G$-action on $X$, we have that $\mathcal{S}_x = G$, and therefore $\mathcal{S}_x \backslash \gm / \gm^g$ consists only of the class $[e]$, so that we can write
\begin{eqnarray*}
& = &  1_{\mathbf{r}(x)   \gm} * \gm g \gm * [\alpha_{g^{-1}}(a)]_{x g \gm}\,.
\end{eqnarray*}
Moreover, since the $G$-actions on $\a$ and $X$ are trivial we can furthermore write
\begin{eqnarray*}
& = &  1_{\mathbf{r}(x) g^{-1}  \gm} * \gm g \gm * [a]_{x \gm}\,.
\end{eqnarray*}
Now, by the same reasoning as above and using expression (\ref{gm g gm times a x gm expression}) we have
\begin{eqnarray*}
& = & \sum_{[\gamma] \in \mathcal{S}_x \backslash \gm / \gm^{g^{-1}} }  1_{\mathbf{r}(x) \gamma g^{-1}  \gm} * \gm g \gm * [a]_{x \gm}\\
& = & \gm g \gm *[a]_{x \gm}\,.
\end{eqnarray*}
Thus we conclude that $[a]_{x\gm} * \gm g \gm =  \gm g \gm * [a]_{x\gm}$. By Theorem \ref{prop decomp of f in crossed prod} we know that elements of the form $[a]_{x \gm}*\gm g \gm * 1_{\mathbf{s}(x) g \gm}$ span $C_c(\a / \gm) \times^{alg}_{\alpha} G / \gm$, and from the commutation relation we just proved it follows that
\begin{eqnarray*}
 [a]_{x \gm}*\gm g \gm * 1_{\mathbf{s}(x) g \gm} & = & \gm g \gm * [a]_{x \gm} * 1_{\mathbf{s}(x) g \gm}\\
& = & \gm g \gm * [a]_{x \gm} * 1_{\mathbf{s}(x)  \gm}\\
& = & \gm g \gm * [a]_{x \gm} \\
& = & [a]_{x \gm}*\gm g \gm\,,
\end{eqnarray*}
so that $C_c(\a / \gm) \times^{alg}_{\alpha} G / \gm$ is spanned by elements of the form $a_{x \gm } * \gm g \gm$. We now conclude that the image of the map (\ref{mapping from tensor product to crossed prod}) is the whole $C_c(\a / \gm) \times^{alg}_{\alpha} G / \gm$.

The fact that this map is a $^*$-homomorphism also follows directly from the commutation relation proved above. \qed\\
\end{ex}

\subsection{Representation theory}
\label{representation theory of crossed products by Hecke pairs section}

In this section we develop the representation theory of crossed products by Hecke pairs. We will introduce the notion of a \emph{covariant pre-representation} and show that there is a bijective correspondence between covariant pre-representations and representations of the crossed product, in a similar fashion to the theory of crossed products by groups.

Recall from Proposition \ref{ext reps to MB(A)} that every nondegenerate $^*$-representation $\pi : C_c(\a / \gm)  \to B(\mathscr{H})$ extends uniquely to a $^*$-representation
\begin{align*}
 \widetilde{\pi}:M_B(C_c(\a / \gm)) \to B(\mathscr{H})\,.
\end{align*}
 We will use the notation $\widetilde{\pi}$ to denote this extension throughout this section, many times without any reference. Since $C_c(X^0 / \gm)$ is a $BG^*$-algebra (it is spanned by projections) we naturally have $C_c(X^0 / \gm) \subseteq M_B(C_c(\a / \gm))$.\\

\begin{df}
\label{covariant pre rep def}
 Let $\pi$ be a nondegenerate $^*$-representation of $C_c(\a / \gm)$ on a Hilbert space $\mathscr{H}$ and $\widetilde{\pi}$ its unique extension to a $^*$-representation of $M_B(C_c(\a / \gm))$. Let $\mu$ be a unital pre-$^*$-representation of $\h(G, \gm)$ on the inner product space $\mathscr{W}:=\pi(C_c(\a / \gm))\mathscr{H}$. We say that $(\pi, \mu)$ is a \emph{covariant pre-$^*$-representation} if the following equality
\begin{eqnarray}
& & \mu(\gm g \gm)\pi([a]_{x\gm})\mu(\gm s \gm) \;\; = \;\; \label{cov pre rep equality}\\ 
& = & \!\!\!\!  \sum_{\substack{[u] \in \gm g^{-1} \gm / \gm \\ [v] \in \gm s \gm / \gm}}\; \sum_{[\gamma] \in E_{u,v}^{\mathbf{s}(x)}} \!\! \frac{\Delta(g)N_{u^{-1}, v}^{\mathbf{s}(x)\gamma}}{L(u^{-1}v)} \; \widetilde{\pi}([\alpha_{u^{-1}\gamma^{-1}}(a)]_{x \gamma u \gm}) \, \mu(\gm u^{-1} v \gm) \, \widetilde{\pi}(1_{\mathbf{s}(x) \gamma v \gm})\,, \notag
\end{eqnarray}
holds on $L(\mathscr{W})$, for all $g, s \in G$ and $x \in X$.\\
\end{df}

 Condition (\ref{cov pre rep equality}) simply says that the pair $(\pi, \mu)$ must preserve the structure of products of the form $\gm g \gm *[a]_{x\gm}*\gm s \gm$, when expressed in terms of the canonical spanning set of elements of the form $[b]_{y\gm}*\gm d \gm*1_{\mathbf{s}(y)d\gm}$, as explicitly described in  Corollary \ref{cor prod in corossed prod}.

The reader should note the similarity between our definition of a covariant pre-$^*$-representation and the \emph{covariant pairs} of an Huef, Kaliszewski and Raeburn in \cite[Definition 1.1]{cov}. Their notion of covariant pairs served as a motivation for us and is actually a particular case of our Definition \ref{covariant pre rep def}, as we shall see in the sequel to this article.

The operators $\widetilde{\pi}([\alpha_{u^{-1}\gamma^{-1}}(a)]_{x \gamma u \gm}) \, \mu(\gm u^{-1} v \gm) \, \widetilde{\pi}(1_{\mathbf{s}(x) \gamma v \gm})$ in expression (\ref{cov pre rep equality}) are all bounded, as we will now show, and are therefore defined in the whole Hilbert space $\mathscr{H}$.\\

\begin{thm}
\label{operators in covariant def are bounded}
Let $\pi:C_(\a / \gm) \to B(\mathscr{H})$ be a nondegenerate $^*$-representation and $\mu:\h(G, \gm) \to L(\mathscr{W})$ a pre-$^*$-representation on the inner product space $\mathscr{W}:=\pi(C_c(\a / \gm))$. Every element of the form
\begin{align*}
\pi([a]_{x\gm})\mu(\gm g \gm) \widetilde{\pi}(1_{\mathbf{s}(x)g\gm})\,,
\end{align*}
is a bounded operator on $\mathscr{W}$ and therefore extends uniquely to the whole Hilbert space $\mathscr{H}$.\\
\end{thm}

We will need some preliminary facts and lemmas in order to prove Theorem \ref{operators in covariant def are bounded}. These auxiliary results will also be useful later in this section.

Let $\pi: C_c(\a / \gm) \to B(\mathscr{H})$ be a  nondegenerate $^*$-representation and $\widetilde{\pi}$ its extension to $M_B(C_c(\a / \gm))$. For any unit $u \in X^0$ the operator $\widetilde{\pi}(1_{u \gm}) \in B(\mathscr{H})$ is a projection, and therefore $\widetilde{\pi}(1_{u \gm})\mathscr{H}$ is a Hilbert subspace. The fiber $(\a / \gm)_{u \gm}$ is a $C^*$-algebra which we can naturally identify with the $^*$-subalgebra
\begin{align*}
 \{[a]_{u \gm} \in C_c(\a / \gm) : [a] \in (\a / \gm)_{u \gm} \} \subseteq C_c(\a / \gm)\,,
\end{align*}
under the identification given by
\begin{align*}
 (\a / \gm)_{u \gm} \ni [a]\; \longleftrightarrow \; [a]_{u \gm} \in C_c(\a / \gm)\,.
\end{align*}
The $^*$-representation $\widetilde{\pi}$ when restricted to $(\a / \gm)_{u \gm}$, under the above identification, leaves the subspace $\widetilde{\pi}(1_{u \gm})\mathscr{H}$ invariant, because
\begin{align*}
 \widetilde{\pi}([a]_{u \gm}) \widetilde{\pi}(1_{u \gm}) \xi = \widetilde{\pi}([a]_{u \gm}) \xi = \widetilde{\pi}(1_{u \gm}) \widetilde{\pi}([a]_{u \gm}) \xi\,.
\end{align*}
The following lemma assures that this restriction is nondegenerate.\\

\begin{lemma}
\label{rep of Au is nondeg}
 Let $\pi : C_c(\a / \gm)  \to B(\mathscr{H})$ be a nondegenerate $^*$-representation and $\widetilde{\pi}$ its unique extension to $M_B(C_c(\a / \gm))$. The $^*$-representation of $(\a / \gm)_{u \gm}$ on the Hilbert space $\widetilde{\pi}(1_{u \gm})\mathscr{H}$, as above, is nondegenerate.\\
\end{lemma}

{\bf \emph{Proof:}} Let $\widetilde{\pi}(1_{u \gm}) \xi$ be an element of $\widetilde{\pi}(1_{u \gm})\mathscr{H}$ such that
\begin{align*}
 \widetilde{\pi}([a]_{u \gm}) \widetilde{\pi}(1_{u \gm}) \xi = 0\,,
\end{align*}
for all $[a] \in (\a / \gm)_{u \gm}$. We want to prove that $\widetilde{\pi}(1_{u \gm}) \xi = 0$. To see this, let $x \in X$ and $[b] \in (\a / \gm)_{x \gm}$. We have two alternatives: either $\mathbf{s}(x) \gm \neq u \gm$ or $\mathbf{s}(x) \gm = u \gm$. In the first case we see that
\begin{align*}
 \widetilde{\pi}([b]_{x \gm}) \widetilde{\pi}(1_{u \gm}) \xi = \widetilde{\pi}([b]_{x \gm} \cdot 1_{u \gm}) \xi = 0\,,
\end{align*}
whereas for the second we see that
\begin{eqnarray*}
 \|\widetilde{\pi}([b]_{x \gm}) \widetilde{\pi}(1_{u \gm}) \xi\|^2 & = & \langle  \widetilde{\pi}( [b]_{x \gm}) \widetilde{\pi}(1_{u \gm}) \xi \, , \, \widetilde{\pi}( [b]_{x \gm}) \widetilde{\pi}(1_{u \gm}) \xi \rangle\\
 & = & \langle  \widetilde{\pi}(  [b^*b]_{\mathbf{s}(x) \gm}) \widetilde{\pi}(1_{u \gm}) \xi \, , \,  \widetilde{\pi}(1_{u \gm}) \xi \rangle\\
 & = & \langle  \widetilde{\pi}( [b^*b]_{u \gm}) \widetilde{\pi}(1_{u \gm}) \xi \, , \,  \widetilde{\pi}(1_{u \gm}) \xi \rangle\\
& = & 0\,,
\end{eqnarray*}
by assumption. Thus, in any case we have $\widetilde{\pi}([b]_{x \gm}) \widetilde{\pi}(1_{u \gm}) \xi = 0$ for all $x \in X$ and $[b] \in (\a / \gm)_{x \gm}$. By nondegeneracy of $\pi$, this implies that $\widetilde{\pi}(1_{u \gm}) \xi = 0$, as we wanted to prove. \qed\\

\begin{lemma}
\label{pi CcA H= tilde pi CcX0 H lemma}
Let $\pi$ be a nondegenerate $^*$-representation of $C_c(\a / \gm)$ on a Hilbert space $\mathscr{H}$. We have that $\pi(C_c(\a / \gm)) \mathscr{H} = \widetilde{\pi}(C_c(X^0 / \gm))\mathscr{H}$.\\
\end{lemma}

{\bf \emph{Proof:}} It is clear that $\pi(C_c(\a / \gm)) \mathscr{H} \subseteq \widetilde{\pi}(C_c(X^0 / \gm))\mathscr{H}$ since for any element of the form $[a]_{x\gm}$ in $C_c(\a / \gm)$ and $\xi \in \mathscr{H}$ we have $\pi([a]_{x\gm})\xi = \pi(1_{\mathbf{r}(x)\gm} [a]_{x\gm}) \xi = \widetilde{\pi}(1_{\mathbf{r}(x)\gm})\pi([a]_{x\gm})\xi$.

Let us now prove that $\widetilde{\pi}(C_c(X^0 / \gm))\mathscr{H} \subseteq \pi(C_c(\a / \gm)) \mathscr{H}$. Let $u\gm \in X^0 / \gm$ and $\xi \in \mathscr{H}$. We know, by Lemma \ref{rep of Au is nondeg}, that $\pi$ gives a nondegenerate $^*$-representation of $(\a / \gm)_{u\gm}$ on $\widetilde{\pi}(1_{u\gm})\mathscr{H}$. Since $(\a / \gm)_{u \gm}$ is a $C^*$-algebra we have, by the general version of Cohen's factorization theorem (\cite[Theorem 5.2.2]{palmer11}), that there exists $[c] \in (\a / \gm)_{u \gm}$ and $\eta \in \widetilde{\pi}(1_{u\gm}) \mathscr{H}$ such that
\begin{align*}
\widetilde{\pi}(1_{u\gm}) \xi = \pi([c]_{u\gm})\eta\,,
\end{align*}
which means that $\widetilde{\pi}(1_{u\gm}) \xi \in \pi(C_c(\a / \gm)) \mathscr{H}$. This finishes the proof. \qed\\

{\bf \emph{Proof of Theorem \ref{operators in covariant def are bounded}:}} The operator $\pi([a]_{x\gm}) \mu(\gm g \gm)\widetilde{\pi}(1_{\mathbf{s}(x)g \gm})$ is clearly defined on the inner product space $\pi(C_c(\a / \gm)) \mathscr{H}$. By Lemma \ref{pi CcA H= tilde pi CcX0 H lemma} this operator is then defined on the space $\widetilde{\pi}(C_c(X^0 / \gm))\mathscr{H}$. Since
\begin{align*}
\pi([a]_{x\gm}) \mu(\gm g \gm)\widetilde{\pi}(1_{\mathbf{s}(x)g \gm}) = \pi([a]_{x\gm}) \mu(\gm g \gm)\widetilde{\pi}(1_{\mathbf{s}(x)g \gm})\widetilde{\pi}(1_{\mathbf{s}(x)g \gm})\,,
\end{align*}
it follows that the operator $\pi([a]_{x\gm}) \mu(\gm g \gm)\widetilde{\pi}(1_{\mathbf{s}(x)g \gm})$ is actually defined in the whole Hilbert space $\mathscr{H}$ (or in other words, it extends canonically to  $\mathscr{H}$).

A similar argument shows that $\widetilde{\pi}(1_{\mathbf{s}(x)g \gm}) \mu((\gm g \gm)^*)\pi([a^*]_{x^{-1}\gm})$ is also defined in the whole Hilbert space $\mathscr{H}$ and it is easy to see that $\pi([a]_{x\gm}) \mu(\gm g \gm)\widetilde{\pi}(1_{\mathbf{s}(x)g \gm})$ is an adjointable operator on $\mathscr{H}$, whose adjoint is $\widetilde{\pi}(1_{\mathbf{s}(x)g \gm}) \mu((\gm g \gm)^*)\pi([a^*]_{x^{-1}\gm})$. Since adjointable operators on a Hilbert space are necessarily bounded (see \cite[Proposition 9.1.11]{palmer}), it follows that $\pi([a]_{x\gm}) \mu(\gm g \gm)\widetilde{\pi}(1_{\mathbf{s}(x)g \gm})$ is a bounded operator. \qed\\

The striking feature that we actually have to consider pre-representations of $\h(G, \gm)$, and not just representations, was not present in the theory of crossed products by groups because a group algebra $\mathbb{C}(G)$ of a discrete group is always a $BG^*$-algebra and therefore all of its pre-representations come from true representations (see further Remark \ref{cov pre reps when gm is normal}).

It will be useful to distinguish between covariant pre-$^*$-representations and covariant $^*$-representations, so we will treat them in separate definitions. As will be discussed below we will see covariant $^*$-representations as a particular type of covariant pre-$^*$-representations.\\

\begin{df}
\label{covariant rep def}
 Let $\pi$ be a nondegenerate $^*$-representation of $C_c(\a / \gm)$ on a Hilbert space $\mathscr{H}$ and $\mu$ a unital $^*$-representation of $\h(G, \gm)$ on $\mathscr{H}$. We say that $(\pi, \mu)$ is a \emph{covariant $^*$-representation} if equality (\ref{cov pre rep equality}) holds in $B(\mathscr{H})$ for all $g, s \in G$ and $x \in X$.\\
\end{df}

\begin{lemma}
\label{mu leaves W invariant}
 Let $(\pi, \mu)$ be a covariant $^*$-representation on a Hilbert space $\mathscr{H}$. Then $\mu$ leaves the subspace $\mathscr{W}:= \pi(C_c(\a / \gm)) \mathscr{H}$ invariant.\\
\end{lemma}

{\bf \emph{Proof:}} Consider elements of the form $\pi([a]_{x\gm})\xi$, whose span gives $\mathscr{W}$. Using the fact that $\mu$ is unital and the covariance relation (\ref{cov pre rep equality}) we see that
\begin{eqnarray*}
& &\!\!\!\! \mu(\gm g\gm) \pi([a]_{x\gm})\xi \;=\\
& = &\!\!\!\! \mu(\gm g\gm) \pi([a]_{x\gm}) \mu(\gm)\xi\\
& = &\!\!\!\!\! \sum_{[u] \in \gm g^{-1} \gm / \gm}\; \sum_{[\gamma] \in E_{u,e}^{\mathbf{s}(x)}} \!\! \frac{\Delta(g)N_{u^{-1}, e}^{\mathbf{s}(x)\gamma}}{L(u^{-1})} \; \widetilde{\pi}([\alpha_{u^{-1} \gamma^{-1}}(a)]_{x \gamma u \gm}) \, \mu(\gm u^{-1} \gm) \, \widetilde{\pi}(1_{\mathbf{s}(x) \gamma \gm})\xi\,.
\end{eqnarray*}
Hence, $ \mu(\gm g\gm) \pi([a]_{x\gm})\xi \in \mathscr{W}$, and consequently $\mu(\gm g\gm)$ leaves $\mathscr{W}$ invariant. This finishes the proof. \qed\\

From a covariant $^*$-representation $(\pi, \mu)$ one can obtain canonically a covariant pre-$^*$-representation $(\pi, \mu)$, just by restricting $\mu$ to the dense subspace $\mathscr{W}:= \pi(C_c(\a / \gm))\mathscr{H}$ (which is an invariant subspace by Lemma \ref{mu leaves W invariant}). So we can regard covariant $^*$-representations as a special kind of covariant pre-$^*$-representations: they are exactly those for which $\mu$ is normed. As we shall see later in Example \ref{example of int form normed but mu not normed}, there are covariant pre-$^*$-representations which are not covariant $^*$-representations, thus in general the latter form a proper subclass of the former. We shall also see examples where they actually coincide.\\

\begin{rem}
 Equivalently, one could define covariant (pre-)$^*$-representation using any other of the equalities in Corollary \ref{cor prod in corossed prod} and substituting with the appropriate (pre-)$^*$-representations. It is easy to see, using completely analogous arguments as in the proof of Corollary \ref{cor prod in corossed prod} or Theorem \ref{product formula}, that all three expressions yield the same result.\\
\end{rem}

\begin{rem}
\label{cov pre reps when gm is normal}

 Even though it might not be entirely clear from the start, when $\gm$ is a normal subgroup of $G$ the definition of a covariant pre-representation is nothing but the usual definition of covariant representation of the system $(C_c(\a / \gm), G / \gm)$. We recall that a covariant representation of $(C_c(\a / \gm), G / \gm)$ is a pair $(\pi, U)$ consisting of a nondegenerate $^*$-representation $\pi$ of $C_c(\a / \gm)$ and a unitary representation $U$ of $G / \gm$ satisfying the relation
\begin{align*}
 \pi(\overline{\alpha}_{g\gm}(f)) = U_{g \gm} \pi(f) U_{g^{-1}\gm}\,,
\end{align*}
for all $f \in C_c(\a / \gm)$ and $g\gm \in G / \gm$. Now, as it is well known, every unitary representation $U$ of $G / \gm$ is associated in a canonical way to a unital $^*$-representation $\mu$ of the group algebra $\mathbb{C}(G / \gm)$, so that we can write the covariance condition as $\pi(\overline{\alpha}_{g \gm}(f)) = \mu(g \gm) \pi(f) \mu(g^{-1}\gm)$. As a consequence we have that for any $g \gm, s \gm \in G / \gm$, $x \in X$ and $a \in \a_x$:
\begin{align*}
 \mu(g\gm) \pi([a]_{x\gm}) \mu(s\gm) = \pi([\alpha_g(a)]_{xg^{-1} \gm}) \mu(g^{-1} s\gm)\,.
\end{align*}
We want to check that covariant representations of the system $(C_c(\a / \gm), G / \gm)$ are the same as covariant pre-$^*$-representations as in Definition \ref{covariant pre rep def}.

Given a covariant pre-$^*$-representation $(\pi, \mu)$ on some Hilbert space $\mathscr{H}$ in the sense of Definition \ref{covariant pre rep def}, we have that $\mu$ is a pre-$^*$-representation of $\mathbb{C}(G / \gm)$, which is normed since any group algebra of a discrete group is a $BG^*$-algebra, and thus we can see $\mu$ as a true $^*$-representation on $\mathscr{H}$. We then have that
\begin{eqnarray*}
\!\!& & \mu(g \gm) \pi([a]_{x\gm}) \mu(g^{-1} \gm)\\ & = & \mu(\gm g \gm) \pi([a]_{x\gm}) \mu(\gm g^{-1} \gm)\\
& = & \!\!\!\!\! \sum_{\substack{[u] \in \gm g^{-1} \gm / \gm \\ [v] \in \gm g^{-1} \gm / \gm}}\; \sum_{[\gamma] \in E_{u,v}^{\mathbf{s}(x)}} \!\! \frac{\Delta(g)N_{u^{-1}, v}^{\mathbf{s}(x)\gamma}}{L(u^{-1}v)} \; \widetilde{\pi}([\alpha_{u^{-1}\gamma^{-1}}(a)]_{x \gamma u \gm}) \, \mu(\gm u^{-1} v \gm) \, \widetilde{\pi}(1_{\mathbf{s}(x) \gamma v \gm})\\
& = &\!\!\!\! \sum_{[\gamma] \in E_{g^{-1},g^{-1}}^{\mathbf{s}(x)}} N_{g, g^{-1}}^{\mathbf{s}(x)\gamma} \; \widetilde{\pi}([\alpha_{g}(a)]_{x  g^{-1} \gm}) \, \mu( g g^{-1} \gm) \, \widetilde{\pi}(1_{\mathbf{s}(x) g^{-1} \gm})\\
& = &\!\!\!\! \sum_{[\gamma] \in E_{g^{-1},g^{-1}}^{\mathbf{s}(x)}} N_{g, g^{-1}}^{\mathbf{s}(x)\gamma} \; \widetilde{\pi}([\alpha_g(a)]_{x  g^{-1} \gm} \cdot 1_{\mathbf{s}(x) g^{-1} \gm})\\
& = &\!\!\!\! \sum_{[\gamma] \in E_{g^{-1},g^{-1}}^{\mathbf{s}(x)}} N_{g, g^{-1}}^{\mathbf{s}(x)\gamma} \; \pi([\alpha_g(a)]_{x  g^{-1} \gm})\,.
\end{eqnarray*}
It is clear from the normality of $\gm$ that $E_{g^{-1},g^{-1}}^{\mathbf{s}(x)}$ consists only of the class $[e]$ and moreover $N_{g, g^{-1}}^{\mathbf{s}(x)} = 1$, so that
\begin{eqnarray*}
\mu(g \gm) \pi([a]_{x\gm}) \mu(g^{-1} \gm) & = &  \pi([\alpha_g(a)]_{x  g^{-1} \gm}) \,.
\end{eqnarray*}
By linearity it follows that $\mu(g \gm) \pi(f) \mu(g^{-1} \gm)  =   \pi(\overline{\alpha}_{g\gm}(f))$ for any $f \in C_c(\a / \gm)$. Thus, with $U$ being the unitary representation of $G / \gm$ associated to $\mu$, we see that $(\pi, U)$ is covariant representation of the system $(C_c(\a / \gm), G / \gm)$.

For the other direction, let $(\pi, U )$ be a covariant representation of the system $(C_c(\a / \gm), G / \gm)$ and let $\mu$ be the $^*$-representation of $\mathbb{C}(G / \gm)$ associated to $U$, which we restrict to the inner product space $\pi(C_c(\a / \gm)) \mathscr{H}$. We want to prove that $(\pi, \mu)$ is a covariant pre-$^*$-representation in the sense of Definition \ref{covariant pre rep def}. We have
\begin{eqnarray*}
 & & \mu(g \gm) \pi([a]_{x\gm}) \mu(s \gm)\\ & = & \mu(g \gm) \pi([a]_{x\gm}) \mu(g^{-1} \gm) \mu(gs \gm)\\
& = &  \pi([\alpha_g(a)]_{x  g^{-1} \gm}) \mu(gs\gm)\\
& = &\!\!\!\!\!\! \sum_{\substack{[u] \in \gm g^{-1} \gm / \gm \\ [v] \in \gm s \gm / \gm}}\; \sum_{[\gamma] \in E_{u,v}^{\mathbf{s}(x)}}\!\! \frac{\Delta(g)N_{u^{-1}, v}^{\mathbf{s}(x)\gamma}}{L(u^{-1}v)} \; \widetilde{\pi}([\alpha_{u^{-1}\gamma^{-1}}(a)]_{x \gamma u \gm}) \, \mu(\gm u^{-1} v \gm) \, \widetilde{\pi}(1_{\mathbf{s}(x) \gamma v \gm})\,,
\end{eqnarray*}
where the last equality is obtained following analogous computations as those above. Thus, $(\pi, \mu)$ is a covariant pre-$^*$-representation in the sense of Definition  \ref{covariant pre rep def}.\\
\end{rem}

The following result makes it clear that some of the relations we have inside the crossed product (see Proposition \ref{invol in crossed prod}) are preserved upon taking covariant pre-$^*$-representations. This is expected since, as we stated before, we will prove that covariant pre-representations give rise to representations of the crossed product, and this result is the first step in that direction:\\

\begin{prop}
\label{analogue strange eq covariance}
 Let $(\pi, \mu)$ be a covariant pre-$^*$-representation. The following two equalities hold:
\begin{align}
\label{first strange eq covariance}
\widetilde{\pi}( 1_{\mathbf{r}(x) \gm}) \mu(  \gm g \gm ) \widetilde{\pi}(  [\alpha_{g^{-1}}(a)]_{xg\gm})\; =\;  \widetilde{\pi}([a]_{x \gm}) \mu(\gm g \gm) \widetilde{\pi}( 1_{\mathbf{s}(x) g \gm})\,.
\end{align}

\begin{align}
\label{second strange eq covariance}
\mu(\gm g \gm) \widetilde{\pi}([a]_{x \gm})  = \sum_{[\gamma] \in E^{\mathbf{s}(x)}_{g^{-1}, e}} \widetilde{\pi}(1_{\mathbf{r}(x) \gamma g^{-1} \gm}) \mu( \gm g  \gm ) \widetilde{\pi}( [a]_{x \gm})\,.
\end{align}\\

\end{prop}

{\bf \emph{Proof:}} Since $(\pi, \mu)$ is a covariant pre-$^*$-representation we have
\begin{eqnarray*}
 \mu(\gm g \gm) \widetilde{\pi}([a]_{x \gm}) & = & \mu(\gm g \gm) \widetilde{\pi}([a]_{x \gm}) \mu(\gm)\\
& = & \sum_{[\gamma] \in E^{\mathbf{s}(x)}_{g^{-1}, e}} N^{\mathbf{s}(x) \gamma}_{g,e}\, \widetilde{\pi}( [\alpha_{g\gamma^{-1}}(a)]_{x \gamma g^{-1} \gm})  \mu( \gm g  \gm ) \widetilde{\pi}(1_{\mathbf{s}(x) \gamma \gm})\\
& = & \sum_{[\gamma] \in E^{\mathbf{s}(x)}_{g^{-1}, e}} \widetilde{\pi}( [\alpha_{g\gamma^{-1}}(a)]_{x \gamma g^{-1} \gm})  \mu( \gm g  \gm ) \widetilde{\pi}(1_{\mathbf{s}(x) \gm})\,,
\end{eqnarray*}
where the last equality comes from the fact that $n_{g,e}^{\mathbf{s}(x) \gamma} = 1 = d_{g, e}^{\mathbf{s}(x) \gamma}$, and thus $N^{\mathbf{s}(x) \gamma}_{g,e}=1$. From this it follows that
\begin{eqnarray*}
 & & \widetilde{\pi}( 1_{\mathbf{r}(x) g^{-1} \gm}) \mu(  \gm g \gm ) \widetilde{\pi}(  [a]_{x\gm}) \;\; = \;\;\\
  & = & \sum_{[\gamma] \in E^{\mathbf{s}(x)}_{g^{-1}, e}} \widetilde{\pi}( 1_{\mathbf{r}(x) g^{-1} \gm}) \widetilde{\pi}( [\alpha_{g\gamma^{-1}}(a)]_{x \gamma g^{-1} \gm})  \mu( \gm g  \gm ) \widetilde{\pi}(1_{\mathbf{s}(x) \gm})\\
 & = & \sum_{[\gamma] \in E^{\mathbf{s}(x)}_{g^{-1}, e}} \widetilde{\pi}( 1_{\mathbf{r}(x) g^{-1} \gm} \cdot [\alpha_{g\gamma^{-1}}(a)]_{x \gamma g^{-1} \gm})  \mu( \gm g  \gm ) \widetilde{\pi}(1_{\mathbf{s}(x) \gm})\,.
\end{eqnarray*}
Now the product $1_{\mathbf{r}(x) g^{-1}\gm} \cdot [\alpha_{g\gamma^{-1}}(a)]_{x \gamma g^{-1} \gm}$ is nonzero only when $\mathbf{r}(x) g^{-1} \gm = \mathbf{r}(x) \gamma g^{-1} \gm$, from which one readily concludes that $\mathbf{r}(x) \gamma \in \mathbf{r}(x) g^{-1} \gm g$. Since one trivially has $\mathbf{r}(x) \gamma \in \mathbf{r}(x) \gm$ we conclude that
\begin{align*}
 \mathbf{r}(x) \gamma \in \mathbf{r}(x) \gm \cap \mathbf{r}(x) g^{-1} \gm g\,,
\end{align*}
and by the $\gm$-intersection property we have $\mathbf{r}(x) \gamma \in \mathbf{r}(x) \gm^{g^{-1}}$. From Proposition \ref{xKH bijection with Sx K H} this means that $[\gamma] = [e]$ in $E^{\mathbf{r}(x)}_{g^{-1},e}$. We recall that $E^{\mathbf{r}(x)}_{g^{-1},e} = \mathcal{S}_{\mathbf{r}(x)} \backslash \gm / \gm^{g^{-1}}$, and since $\gm^{g^{-1}} \subseteq \gm$ we have by Proposition \ref{double coset spaces prop} that $[\gamma] \to [\gamma]$ defines a canonical bijection between  $E^{\mathbf{r}(x)}_{g^{-1},e}$ and $(\mathcal{S}_{\mathbf{r}(x)} \cap \gm) \backslash \gm / \gm^{g^{-1}}$. Since the $G$-action is $\gm$-good we necessarily have $\mathcal{S}_{\mathbf{s}(x)} \cap \gm = \mathcal{S}_x \cap \gm = \mathcal{S}_{\mathbf{r}(x)} \cap \gm$, and therefore using Proposition \ref{double coset spaces prop} one more time we can say that $E^{\mathbf{r}(x)}_{g^{-1},e} = E^{\mathbf{s}(x)}_{g^{-1},e}$. Hence, we can say that $[\gamma] = [e]$ in $E^{\mathbf{s}(x)}_{g^{-1},e}$. We conclude that
\begin{eqnarray*}
 \widetilde{\pi}( 1_{\mathbf{r}(x) g^{-1} \gm}) \mu(  \gm g \gm ) \widetilde{\pi}(  [a]_{x\gm}) & = & \widetilde{\pi}( 1_{\mathbf{r}(x) g^{-1} \gm} \cdot [\alpha_g(a)]_{x g^{-1} \gm})  \mu( \gm g  \gm ) \widetilde{\pi}(1_{\mathbf{s}(x) \gm})\\
& = & \widetilde{\pi}([\alpha_g(a)]_{x g^{-1} \gm})  \mu( \gm g  \gm ) \widetilde{\pi}(1_{\mathbf{s}(x) \gm})\,.
\end{eqnarray*}
Since the last expression is valid for any $x \in X$ and $[a] \in (\a / \gm)_{x\gm}$, if we take $x$ to be $xg$ and $[a]$ to be $[\alpha_{g^{-1}}(a)]$ we obtain the desired equality (\ref{first strange eq covariance}):
\begin{align*}
 \widetilde{\pi}( 1_{\mathbf{r}(x) \gm}) \mu(  \gm g \gm ) \widetilde{\pi}(  [\alpha_{g^{-1}}(a)]_{x g \gm}) = \widetilde{\pi}([a]_{x \gm})  \mu( \gm g  \gm ) \widetilde{\pi}(1_{\mathbf{s}(x) g \gm})\,.
\end{align*}
Let us now prove equality  (\ref{second strange eq covariance}). Using the equality in beginning of this proof and equality (\ref{first strange eq covariance}) which we have just proven, we get precisely
\begin{eqnarray*}
 \mu(\gm g \gm) \widetilde{\pi}([a]_{x \gm}) & = & \sum_{[\gamma] \in E^{\mathbf{s}(x)}_{g^{-1}, e}} \widetilde{\pi}( [\alpha_{g \gamma^{-1}}(a)]_{x \gamma g^{-1} \gm})  \mu( \gm g  \gm ) \widetilde{\pi}(1_{\mathbf{s}(x) \gm})\\
& = & \sum_{[\gamma] \in E^{\mathbf{s}(x)}_{g^{-1}, e}} \widetilde{\pi}( [\alpha_{g \gamma^{-1}}(a)]_{x \gamma g^{-1} \gm})  \mu( \gm g  \gm ) \widetilde{\pi}(1_{\mathbf{s}(x \gamma g^{-1}) g \gm})\\
& = & \sum_{[\gamma] \in E^{\mathbf{s}(x)}_{g^{-1}, e}} \widetilde{\pi}(1_{\mathbf{r}(x) \gamma g^{-1}  \gm})  \mu( \gm g  \gm ) \widetilde{\pi}( [\alpha_{ \gamma^{-1}}(a)]_{x \gamma \gm}) \\
& = & \sum_{[\gamma] \in E^{\mathbf{s}(x)}_{g^{-1}, e}} \widetilde{\pi}(1_{\mathbf{r}(x) \gamma g^{-1}  \gm})  \mu( \gm g  \gm )  \widetilde{\pi}( [a]_{x \gm}) \,.
\end{eqnarray*}
This finishes the proof. \qed\\

The passage from a covariant pre-representation $(\pi, \mu)$ to a representation of $C_c(\a / \gm) \times_{\alpha}^{alg} G / \gm$ is done via the so-called \emph{integrated form} $\pi \times \mu$, which we now describe:\\

\begin{df}
 Let $(\pi , \mu)$ be a covariant pre-$^*$-representation on a Hilbert space $\mathscr{H}$. We define the \emph{integrated form} of $(\pi, \mu)$ as the function $\pi \times \mu: C_c(\a / \gm) \times_{\alpha}^{alg} G / \gm \to B(\mathscr{H})$\index{pi times mu@$\pi \times \mu$} defined by
\begin{align*}
 [\pi \times \mu] (f) := \sum_{[g] \in \gm \backslash G / \gm}\; \sum_{x\gm^g \in X / \gm^g} \widetilde{\pi}\Big( \big[f(g\gm) (x) \big]_{x \gm} \Big) \,  \mu(\gm g \gm) \, \widetilde{\pi}(1_{\mathbf{s}(x) g \gm})\,.
\end{align*}\\
\end{df}

\begin{rem}
\label{remark about int form on blocks}
 For $f$ of the form $f = a_{x\gm} * \gm g \gm * 1_{\mathbf{s}(x)g\gm}$ we have
\begin{align*}
 [\pi \times \mu] (f) = \widetilde{\pi}([a]_{x\gm}) \, \mu(\gm g \gm) \, \widetilde{\pi}( 1_{\mathbf{s}(x)g\gm})\,.
\end{align*}
Moreover, from equality (\ref{first strange eq covariance}), for $f'$ of the form $f' = 1_{\mathbf{r}(x)\gm} * \gm g \gm * [\alpha_{g^{-1}}(a)]_{x g\gm} $ we have
\begin{align*}
 [\pi \times \mu] (f') = \widetilde{\pi}(1_{\mathbf{r}(x)\gm})\, \mu (\gm g \gm ) \, \widetilde{\pi}( [\alpha_{g^{-1}}(a)]_{x g\gm}) \,.\\
\end{align*}
\end{rem}

\begin{prop}
 The integrated form $\pi \times \mu$ of a covariant pre-$^*$-representation $(\pi, \mu)$ is a well-defined nondegenerate $^*$-representation.\\
\end{prop}

{\bf \emph{Proof:}} First we need to check that the expression that defines $[\pi \times \mu] (f)$ for a given $f \in C_c(\a / \gm) \times_{\alpha}^{alg} G / \gm$ is well-defined. This is proven in an entirely analogous way as in the proof that the expression (\ref{decomp of f in crossed prod}) in Proposition \ref{prop decomp of f in crossed prod} is well-defined. Secondly, we need to show that $[\pi \times \mu] (f)$ makes sense as an element of $B(\mathscr{H})$. From Theorem \ref{operators in covariant def are bounded} we have that
\begin{align*}
 \widetilde{\pi}\Big( \big[f(g\gm) (x) \big]_{x \gm} \Big) \mu(\gm g \gm) \widetilde{\pi}(1_{\mathbf{s}(x)g \gm}) \in B(\mathscr{W})\,,
\end{align*}
 thus, it follows that $[\pi \times \mu] (f) \in B(\mathscr{W})$, and therefore $[\pi \times \mu] (f)$ admits a unique extension to $B(\mathscr{H})$.

Now, it is obvious that $\pi \times \mu$ is a linear transformation. Let us check that it preserves the involution. It is then enough to check it for elements of the form $f = [a]_{x\gm} * \gm g \gm * 1_{\mathbf{s}(x)g\gm}$. Since $(\pi , \mu)$ is a covariant pre-$^*$-representation we have, by Propositions \ref{analogue strange eq covariance} and \ref{invol in crossed prod},
\begin{eqnarray*}
 \big([\pi \times \mu] (f) \big)^* & = & \Delta(g)\; \widetilde{\pi}(1_{\mathbf{s}(x) g \gm})\, \mu(\gm g^{-1}  \gm)\, \widetilde{\pi}([a^*]_{x^{-1}\gm})\\
& = & \Delta(g)\; \widetilde{\pi}(1_{\mathbf{r}(x^{-1}) g \gm})\, \mu(\gm g^{-1}  \gm)\, \widetilde{\pi}([a^*]_{x^{-1} g g^{-1}\gm})\\
& = & \Delta(g)\; \widetilde{\pi}([\alpha_{g^{-1}}(a^*)]_{x^{-1} g \gm}) \, \mu(\gm g^{-1}  \gm)\, \widetilde{\pi}(1_{\mathbf{s}(x^{-1}) g g^{-1} \gm})\\
& = & \Delta(g)\; \widetilde{\pi}([\alpha_{g^{-1}}(a^*)]_{x^{-1} g \gm}) \, \mu(\gm g^{-1}  \gm)\, \widetilde{\pi}(1_{\mathbf{s}(x^{-1}) \gm})\\
& = & [\pi \times \mu]\, (\Delta(g)\;[\alpha_{g^{-1}}(a^*)]_{x^{-1} g \gm} * \gm g^{-1}  \gm * 1_{\mathbf{s}(x^{-1}) \gm})\\
& = & [\pi \times \mu]\, (f^*)\,.
\end{eqnarray*}

Let us now prove that $\pi \times \mu$ preserves products. We will start by proving that
\begin{align}
\label{prod in integrated form}
 [\pi \times \mu] (f_1 * f_2) = [\pi \times \mu] (f_1 )\,[\pi \times \mu] (f_2)\,,
\end{align}
for $f_1:= [a]_{x\gm} * \gm g \gm * 1_{\mathbf{s}(x) g \gm}$ and $f_2:= [b]_{y\gm} * \gm s \gm * 1_{\mathbf{s}(y) s \gm}$. Let us compute the expression on the left side of (\ref{prod in integrated form}). First, we notice that for the product $f_1 * f_2$ to be non-zero one must have $\mathbf{r}(y) \in \mathbf{s}(x)g \gm$, and in this case we obtain

\begin{align*}
f_1 * f_2  =  [a]_{x\gm} * \gm g \gm * [b]_{y\gm} * \gm s \gm * 1_{\mathbf{s}(y) s \gm}
\end{align*}
which by Corollary \ref{cor prod in corossed prod} gives
\begin{eqnarray*}
& = & \!\!\!\!\!\!\!\!\! \sum_{\substack{[u] \in \gm g^{-1} \gm / \gm \\ [v] \in \gm s \gm / \gm\\ [\gamma] \in E_{u,v}^{\mathbf{s}(y)}}} \!\!\!\frac{\Delta(g)N_{u^{-1}, v}^{\mathbf{s}(y)\gamma}}{L(u^{-1}v)} \; [a]_{x\gm} *[\alpha_{u^{-1}\gamma^{-1}}(b)]_{y \gamma u \gm} * \gm u^{-1} v \gm * 1_{\mathbf{s}(y) \gamma v \gm} * 1_{\mathbf{s}(y) s \gm}\\
& = &\!\!\!\!\!\!\!\!\! \sum_{\substack{[u] \in \gm g^{-1} \gm / \gm \\ [v] \in \gm s \gm / \gm \\ [\gamma] \in E_{u,v}^{\mathbf{s}(y)}\\ \mathbf{s}(y) s\gm = \mathbf{s}(y) \gamma v \gm}}\!\!\! \frac{\Delta(g)N_{u^{-1}, v}^{\mathbf{s}(y)\gamma}}{L(u^{-1}v)} \; [a]_{x\gm} *[\alpha_{u^{-1}\gamma^{-1}}(b)]_{y \gamma u \gm} * \gm u^{-1} v \gm * 1_{\mathbf{s}(y) \gamma v \gm}\\
& = &\!\!\!\!\!\!\!\!\! \sum_{\substack{[u] \in \gm g^{-1} \gm / \gm \\ [v] \in \gm s \gm / \gm\\ [\gamma] \in E_{u,v}^{\mathbf{s}(y)}\\ \mathbf{s}(y) s\gm = \mathbf{s}(y) \gamma v \gm}}\!\!\! \frac{\Delta(g)N_{u^{-1}, v}^{\mathbf{s}(y)\gamma}}{L(u^{-1}v)} \; [a]_{x\gm} *[\alpha_{u^{-1}\gamma^{-1}}(b)]_{y \gamma u \gm} * \gm u^{-1} v \gm * 1_{\mathbf{s}(y \gamma u) u^{-1} v \gm}
\end{eqnarray*}
The product $[a]_{x\gm} *[\alpha_{u^{-1}\gamma^{-1}}(b)]_{y \gamma u \gm}$ is always either zero or of the form $[c]_{(x \theta)(y \gamma u) \gm}$, for some $\theta \in \gm$ and $c \in \a_{(x \theta)(y \gamma u)}$. The point is that $\mathbf{s}\big((x \theta)(y \gamma u) \big) = \mathbf{s}(y \gamma u)$, so that each non-zero summand in the last sum above is actually of the form
\begin{align*}
 [c]_{z\gm}* \gm d \gm * 1_{\mathbf{s}(z) d \gm}\,,
\end{align*}
for appropriate $[c] \in (\a / \gm)_{z\gm}$, $z \in X$ and $d \in G$. Thus, by linearity of $\pi \times \mu$ and Remark \ref{remark about int form on blocks} we obtain
\begin{eqnarray*}
& & [\pi \times \mu] (f_1 * f_2) = \\
& = &\!\!\!\!\!\!\!\! \sum_{\substack{[u] \in \gm g^{-1} \gm / \gm \\ [v] \in \gm s \gm / \gm \\ [\gamma] \in E_{u,v}^{\mathbf{s}(y)} \\ \mathbf{s}(y) s\gm = \mathbf{s}(y) \gamma v \gm}} \frac{\Delta(g)N_{u^{-1}, v}^{\mathbf{s}(y)\gamma}}{L(u^{-1}v)} \; \widetilde{\pi}([a]_{x\gm} \cdot [\alpha_{u^{-1}\gamma^{-1}}(b)]_{y \gamma u \gm})\, \mu( \gm u^{-1} v \gm)\, \widetilde{\pi}( 1_{\mathbf{s}(y) \gamma v \gm})\,.
\end{eqnarray*}

Let us now compute the expression on the right side of (\ref{prod in integrated form}). We have
\begin{align*}
 [\pi \times \mu](f_1)\,[\pi \times \mu] (f_2) & = & \widetilde{\pi}([a]_{x\gm})\, \mu(\gm g \gm) \, \widetilde{\pi}(1_{\mathbf{s}(x) g \gm}) \, \widetilde{\pi}([b]_{y\gm})\, \mu( \gm s \gm) \, \widetilde{\pi}(1_{\mathbf{s}(y) s \gm})\,.
\end{align*}

For $1_{\mathbf{s}(x) g \gm} \cdot [b]_{y \gm}$ to be non-zero we must have $\mathbf{r}(y) \in \mathbf{s}(x)g \gm$, and in this case we obtain, using the definition of a covariant pre-$^*$-representation, 
\begin{eqnarray*}
& & [\pi \times \mu](f_1)\,[\pi \times \mu] (f_2)  =\\
& = &   \widetilde{\pi}([a]_{x\gm})\, \mu(\gm g \gm) \, \widetilde{\pi}([b]_{y\gm})\, \mu( \gm s \gm) \, \widetilde{\pi}(1_{\mathbf{s}(y) s \gm})\\
& = &\!\!\!\!\!\!\!\!\!\!\! \sum_{\substack{[u] \in \gm g^{-1} \gm / \gm \\ [v] \in \gm s \gm / \gm \\ [\gamma] \in E_{u,v}^{\mathbf{s}(y)}}} \!\!\!\!\!\! \frac{\Delta(g)N_{u^{-1}, v}^{\mathbf{s}(y)\gamma}}{L(u^{-1}v)}  \widetilde{\pi}([a]_{x\gm}[\alpha_{u^{-1}\gamma^{-1}}(b)]_{y \gamma u \gm}) \mu(\gm u^{-1} v \gm)  \widetilde{\pi}(1_{\mathbf{s}(y) \gamma v \gm}) \widetilde{\pi}(1_{\mathbf{s}(y)s\gm}\!)\\
& = &\!\!\!\!\!\!\!\!\!\! \sum_{\substack{[u] \in \gm g^{-1} \gm / \gm \\ [v] \in \gm s \gm / \gm\\ [\gamma] \in E_{u,v}^{\mathbf{s}(y)} \\ \mathbf{s}(y) s\gm = \mathbf{s}(y) \gamma v \gm}}\!\!\!\! \frac{\Delta(g)N_{u^{-1}, v}^{\mathbf{s}(y)\gamma}}{L(u^{-1}v)}  \widetilde{\pi}([a]_{x\gm} \cdot [\alpha_{u^{-1}\gamma^{-1}}(b)]_{y \gamma u \gm})\, \mu(\gm u^{-1} v \gm) \, \widetilde{\pi}(1_{\mathbf{s}(y) \gamma v \gm}\!)\,.
\end{eqnarray*}
Hence, we have proven equality (\ref{prod in integrated form}) for the special case of $f_1$ and $f_2$ being $f_1 := [a]_{x\gm} * \gm g \gm * 1_{\mathbf{s}(x) g \gm}$ and $f_2:=[b]_{y\gm} * \gm s \gm * 1_{\mathbf{s}(y) s \gm}$. Using this we will now show that equality (\ref{prod in integrated form}) holds for any $f_1, f_2 \in C_c(\a / \gm) \times_{\alpha}^{alg} G / \gm$.  In fact, by Proposition \ref{prop decomp of f in crossed prod}, $f_1$ and $f_2$ can be written as sums 
\begin{align*}
 f_1 =  \sum_i v_i \,, \qquad\qquad\qquad f_1 = \sum_j w_j\,,
\end{align*}
where each $v_i$ and $w_j$ is of the form $[a]_{x\gm} * \gm g \gm * 1_{\mathbf{s}(x) g \gm}$, for some $g\gm \in G / \gm$, $x \in X$ and $a \in \a_x$. Since $\pi \times \mu$ is a linear mapping we have
\begin{eqnarray*}
 [\pi \times \mu](f_1 * f_2) & = & [\pi \times \mu]\Big( \big( \sum_i v_i \big) * \big( \sum_j w_j \big) \Big)\\
& = & [\pi \times \mu] \Big( \sum_{i,j} v_i * w_j \Big)\\
& = & \sum_{i,j} [\pi \times \mu] (v_i * w_j)\,,
\end{eqnarray*}
and by the special case of equality (\ref{prod in integrated form}) we have just proven we get
\begin{eqnarray*}
[\pi \times \mu](f_1 * f_2) & = & \sum_{i,j} [\pi \times \mu] (v_i )[\pi \times \mu] ( w_j)\\
& = & \Big( \sum_i [\pi \times \mu] (v_i) \Big) \Big(\sum_j [\pi \times \mu] (w_j)\Big)\\
& = & [\pi \times \mu] \big( \sum_i  v_i \big)  [\pi \times \mu]\big(\sum_j  w_j\big)\\
& = & [\pi \times \mu](f_1)[\pi \times \mu](f_2)\,.
\end{eqnarray*}
Hence, $\pi \times \mu$ is a $^*$-representation. To finish the proof we now only need to show that $\pi \times \mu$ is nondegenerate. The restriction of $\pi \times \mu$ to the $^*$-subalgebra $C_c(\a / \gm)$ is precisely the representation $\pi$. Since $\pi$ is assumed to be nondegenerate it follows that $\pi \times \mu$ must be nondegenerate as well. \qed\\

The next result shows how from a representation of the crossed product one can naturally form a covariant pre-representation.\\

\begin{prop}
\label{rep of crossed product to cov pre rep}
 Let $\Phi: C_c(\a / \gm) \times_{\alpha}^{alg} G / \gm \to B(\mathscr{H})$ be a nondegenerate $^*$-representation. Consider the pair $(\Phi|, \omega_{\Phi})$ defined by
\begin{itemize}
 \item $\Phi|$ is the restriction of $\Phi$ to $C_c(\a / \gm)$.
 \item Let $\widetilde{\Phi}$ be the extension of $\Phi$ to a pre-$^*$-representation (via Proposition \ref{ext hom *-alg}) of $M(C_c(\a / \gm) \times_{\alpha}^{alg} G / \gm)$ on the inner product space $\Phi (C_c(\a / \gm) \times_{\alpha}^{alg} G / \gm) \mathscr{H}$. We define $\omega_{\Phi}$\index{mu phi@$\mu_{\Phi}$} to be the restriction of $\widetilde{\Phi}$ to $\h(G, \gm)$.
\end{itemize}
The pair $(\Phi|, \omega_{\Phi})$ is a covariant pre-$^*$-representation.\\
\end{prop}

We will need some preliminary lemmas in order to prove Proposition \ref{rep of crossed product to cov pre rep}.

\begin{lemma}
\label{restict is nondegenerate too}
 If $\Phi: C_c(\a / \gm) \times_{\alpha}^{alg} G / \gm \to B(\mathscr{H})$ is a nondegenerate $^*$-representation, then its restriction to $C_c(\a / \gm)$ is also nondegenerate.\\
\end{lemma}

{\bf \emph{Proof:}} Let $\xi \in \mathscr{H}$ be such that $\Phi(C_c(\a / \gm)) \,\xi = \{0\}$. We want to show that $\xi = 0$. Since $\Phi$ is nondegenerate, it is then enough to prove that $\Phi(C_c(\a / \gm) \times_{\alpha}^{alg} G / \gm)\, \xi = \{0\}$. Thus, by virtue of Proposition \ref{invol in crossed prod}, it suffices to prove that $\Phi(1_{\mathbf{r}(x) \gm} * \gm g \gm * [\alpha_{g^{-1}}(a)]_{x g \gm})\xi =0$ for all $g \in G$, $x \in X$, $a \in \a_x$. We have
\begin{eqnarray*}
& & \| \Phi(1_{\mathbf{r}(x) \gm} * \gm g \gm * [\alpha_{g^{-1}}(a)]_{x g \gm})\xi\|^2 \;=\\
& = & \!\!\!\! \Delta(g) \langle \Phi( [\alpha_{g^{-1}}(a^*)]_{x^{-1} g \gm} * \gm g^{-1} \gm *  1_{\mathbf{r}(x) \gm}*1_{\mathbf{r}(x) \gm} * \gm g \gm * [\alpha_{g^{-1}}(a)]_{x g \gm}) \xi \,, \,\xi \rangle\\
& = &\!\!\!\! \Delta(g) \langle \Phi( [\alpha_{g^{-1}}(a^*)]_{x^{-1} g \gm}) \Phi( \gm g^{-1} \gm * 1_{\mathbf{r}(x) \gm} * \gm g \gm * [\alpha_{g^{-1}}(a)]_{x g \gm}) \xi \,, \,\xi \rangle\\
& = &\!\!\!\! \Delta(g) \langle \Phi( \gm g^{-1} \gm * 1_{\mathbf{r}(x) \gm} * \gm g \gm * [\alpha_{g^{-1}}(a)]_{x g \gm}) \xi \,, \, \Phi( [\alpha_{g^{-1}}(a)]_{x g \gm}) \xi \rangle\\
& = &\!\!\!\! 0\,.
\end{eqnarray*}
Hence $\xi = 0$ and therefore $\Phi$ restricted to $C_c(\a / \gm)$ is nondegenerate. \qed\\

\begin{lemma}
 \label{two ways of getting reps of CcU}
Let $\Phi:C_c(\a / \gm) \times_{\alpha}^{alg} G / \gm \to B(\mathscr{H})$ be a nondegenerate $^*$-representation and  $\widetilde{\Phi}$ its unique extension to $M_B(C_c(\a / \gm) \times_{\alpha}^{alg} G / \gm)$ (via Proposition \ref{ext reps to MB(A)}). Let $\Phi|$ be the restriction of $\Phi$ to $C_c(\a / \gm)$ and $\widetilde{\Phi|}$ its unique extension to $M_B(C_c(\a / \gm))$. We have that
\begin{align*}
 \widetilde{\Phi}(f) = \widetilde{\Phi|}(f) \,,
\end{align*}
for all $f \in C_c(X^0 / \gm)$. In other words, the two $^*$-representations $\widetilde{\Phi}$ and  $\widetilde{\Phi|}$ are the same in $C_c(X^0/ \gm)$.\\
\end{lemma}

{\bf \emph{Proof:}} By Lemma \ref{restict is nondegenerate too} the subspace $\Phi(C_c(\a / \gm))\mathscr{H}$ is dense in $\mathscr{H}$, so that it is enough to check that $\widetilde{\Phi}(f) \Phi(f_2) \xi = \widetilde{\Phi|}(f) \Phi(f_2) \xi$, for all $f_2 \in C_c(\a / \gm)$ and $\xi \in \mathscr{H}$. By definition of the extension $\widetilde{\Phi}$ (see Proposition \ref{ext reps to MB(A)}) we have
\begin{eqnarray*}
 \widetilde{\Phi}(f) \Phi(f_2) \xi & = & \Phi(f * f_2) \xi\,,
\end{eqnarray*}
where $f*f_2$ is the product of $f$ and $f_2$, which lies inside $C_c(\a / \gm) \times_{\alpha}^{alg} G / \gm$. Since both $f$ and $f_2$ are elements of $B(\a, G, \gm)$ we see the product $f*f_2$ as taking place in $B(\a, G, \gm)$. By definition of the embeddings of $C_c(X^0 / \gm)$ and $C_c(\a / \gm)$ in $B(\a, G, \gm)$ we have that $f*f_2$ is nothing but the element $f\cdot f_2$, where the product is just the product of $f$ and $f_2$ inside $M(C_c(\a))$. As we observed in Section \ref{big algrebra MCcA section}, this product is exactly same as the product of $f$ and $f_2$ in $M(C_c(\a / \gm))$. Thus, the following computation makes sense:
\begin{eqnarray*}
 \widetilde{\Phi}(f) \Phi(f_2) \xi & = & \Phi(f * f_2)\xi \;\; = \;\; \Phi(f \cdot f_2)\xi\\
& = & \Phi|(f \cdot f_2)\xi \;\; =\;\; \widetilde{\Phi|}(f)\Phi|(f_2)\xi\,.
\end{eqnarray*}
This finishes the proof. \qed\\

\begin{lemma}
\label{restict is nondegenerate too formula}
 Let $\Phi: C_c(\a / \gm) \times_{\alpha}^{alg} G / \gm \to B(\mathscr{H})$ be a nondegenerate $^*$-representation. We have that
\begin{align*}
 \Phi(C_c(\a / \gm)) \mathscr{H} \; = \; \Phi (C_c(\a / \gm) \times_{\alpha}^{alg} G / \gm) \mathscr{H}\,.\\
\end{align*}
\end{lemma}

{\bf \emph{Proof:}}  The inclusion $\Phi(C_c(\a / \gm)) \mathscr{H} \subseteq \Phi (C_c(\a / \gm) \times_{\alpha}^{alg} G / \gm) \mathscr{H}$ is obvious. To check the converse inclusion it is enough to prove that
\begin{align*}
 \Phi([a]_{x \gm} * \gm g \gm * 1_{\mathbf{s}(x) \gm}) \xi \in \Phi(C_c(\a / \gm)) \mathscr{H}\,,
\end{align*}
for all $x \in X$, $a \in \a_x$,  $ g \in G$ and $\xi \in \mathscr{H}$.  Let $\widetilde{\Phi}: M_B(C_c(\a / \gm) \times_{\alpha}^{alg} G / \gm) \to B(\mathscr{H})$ be the unique extension of $\Phi$ to a $^*$-representation of $M_B(C_c(\a / \gm) \times_{\alpha}^{alg} G / \gm)$, as in Proposition \ref{ext reps to MB(A)}. We then get
\begin{eqnarray*}
 \Phi([a]_{x \gm} * \gm g \gm * 1_{\mathbf{s}(x) g \gm}) \xi & = & \Phi(1_{\mathbf{r}(x) \gm} *[a]_{x \gm} * \gm g \gm * 1_{\mathbf{s}(x) g \gm}) \xi\\
& = & \widetilde{\Phi}(1_{\mathbf{r}(x) \gm}) \Phi([a]_{x \gm} * \gm g \gm * 1_{\mathbf{s}(x) g \gm})  \xi\,.
\end{eqnarray*}
Denoting by $\Phi|$ the restriction of $\Phi$ to $C_c(\a / \gm)$ we have, by Lemma \ref{two ways of getting reps of CcU}, that
\begin{eqnarray*}
 & = &  \widetilde{\Phi|}(1_{\mathbf{r}(x) \gm}) \Phi([a]_{x \gm} * \gm g \gm * 1_{\mathbf{s}(x) g \gm})  \xi\,,
\end{eqnarray*}
i.e. $\Phi([a]_{x \gm} * \gm g \gm * 1_{\mathbf{s}(x) g \gm}) \xi \in \widetilde{\Phi|}(C_c(X / \gm)) \mathscr{H}$. By Lemma \ref{pi CcA H= tilde pi CcX0 H lemma} it then follows that $\Phi([a]_{x \gm} * \gm g \gm * 1_{\mathbf{s}(x) g \gm}) \xi \in \Phi|(C_c(\a / \gm)) \mathscr{H}$. \qed\\

{\bf \emph{Proof of Proposition \ref{rep of crossed product to cov pre rep}:}} First of all, by Lemma \ref{restict is nondegenerate too}, $\Phi|$ is indeed a nondegenerate $^*$-representation of $C_c(\a / \gm)$. Secondly, from Lemma \ref{restict is nondegenerate too formula}, we have
\begin{align*}
 \Phi(C_c(\a / \gm)) \mathscr{H} \; = \; \Phi (C_c(\a / \gm) \times_{\alpha}^{alg} G / \gm) \mathscr{H}\,.
\end{align*}
Thus, $\omega_{\Phi}$ is a pre-$^*$-representation of $\h(G, \gm)$ on $\mathscr{W}:= \Phi(C_c(\a / \gm)) \mathscr{H}$. We now only need to check covariance. We have
\begin{eqnarray*}
&  & \omega_{\Phi}(\gm g \gm)\Phi|([a]_{x\gm})\omega_{\Phi}(\gm s \gm) \; =\\
& = &\!\! \widetilde{\Phi}(\gm g \gm)\widetilde{\Phi}([a]_{x\gm})\widetilde{\Phi}(\gm s \gm)\\
& = &\!\! \widetilde{\Phi}(\gm g \gm*[a]_{x\gm}*\gm s \gm)\\
& = &\!\! \widetilde{\Phi}\Big( \sum_{\substack{[u] \in \gm g^{-1} \gm / \gm \\ [v] \in \gm s \gm / \gm}} \sum_{[\gamma] \in E_{u,v}^{\mathbf{s}(x)}} \frac{\Delta(g)N_{u^{-1}, v}^{\mathbf{s}(x)\gamma}}{L(u^{-1}v)} \; [\alpha_{u^{-1}\gamma^{-1}}(a)]_{x \gamma u \gm} *\gm u^{-1} v \gm *1_{\mathbf{s}(x) \gamma v \gm} \Big)\\
& = &\!\!\!\!\!\!\!\!\!\! \sum_{\substack{[u] \in \gm g^{-1} \gm / \gm \\ [v] \in \gm s \gm / \gm}} \sum_{[\gamma] \in E_{u,v}^{\mathbf{s}(x)}} \!\!\frac{\Delta(g)N_{u^{-1}, v}^{\mathbf{s}(x)\gamma}}{L(u^{-1}v)} \; \widetilde{\Phi}([\alpha_{u^{-1}\gamma^{-1}}(a)]_{x \gamma u \gm})\widetilde{\Phi}(\gm u^{-1} v \gm )\widetilde{\Phi}(1_{\mathbf{s}(x) \gamma v \gm} )\,.
\end{eqnarray*}
Denoting by $\widetilde{\Phi|}$ the unique extension of $\Phi|$ to $M_B(C_c(\a/ \gm))$ we have, by Lemma \ref{two ways of getting reps of CcU}, that
\begin{eqnarray*}
& = &\!\!\!\!\!\!\!\!\!\! \sum_{\substack{[u] \in \gm g^{-1} \gm / \gm \\ [v] \in \gm s \gm / \gm}} \sum_{[\gamma] \in E_{u,v}^{\mathbf{s}(x)}}\!\! \frac{\Delta(g)N_{u^{-1}, v}^{\mathbf{s}(x)\gamma}}{L(u^{-1}v)} \; \widetilde{\Phi|}([\alpha_{u^{-1}\gamma^{-1}}(a)]_{x \gamma u \gm}) \omega_{\Phi}(\gm u^{-1} v \gm) \widetilde{\Phi|}(1_{\mathbf{s}(x) \gamma v \gm})\,.
\end{eqnarray*}
 This finishes the proof. \qed\\

\begin{thm}
\label{bijective correspondence cov rep - rep}
 There is a bijective correspondence between nondegenerate $^*$-representations of $C_c(\a / \gm) \times_{\alpha}^{alg} G / \gm$ and covariant pre-$^*$-representations. This bijection is given by $ (\pi, \omega) \longmapsto \pi \times \mu$,
with inverse given by $\Phi \longmapsto (\Phi|, \omega_{\Phi})$.\\
\end{thm}

{\bf \emph{Proof:}} We have to prove that the composition of these maps, in both orders, is the identity.

Let $(\pi, \mu)$ be a covariant pre-$^*$-representation and $\pi \times \mu$ its integrated form. We want to show that
\begin{align*}
 \big((\pi \times \mu)|,\, \omega_{\pi \times \mu} \big) = (\pi, \mu)\,.
\end{align*}
By definition of the integrated form we readily have $(\pi \times \mu)| = \pi$. This also implies, via Lemma \ref{restict is nondegenerate too}, that the inner product spaces on which $\mu$ and $\omega_{\pi \times \mu}$ are defined are actually the same. Thus, it remains to be checked that $\omega_{\pi \times \mu} = \mu$. Let $\pi([a]_{x\gm})\xi$ be one of the generators of $\pi(C_c(\a / \gm)) \mathscr{H}$. We have
\begin{eqnarray*}
& & \omega_{\pi \times \mu} (\gm g\gm) \, \pi([a]_{x\gm})\xi \; =\\
& = & \widetilde{[\pi \times \mu]}(\gm g \gm) \, \pi([a]_{x\gm})\xi\\
& = & [\pi \times \mu](\gm g \gm * [a]_{x\gm})\xi\,,
\end{eqnarray*}
and using Proposition \ref{invol in crossed prod}, Remark \ref{remark about int form on blocks} and Proposition \ref{analogue strange eq covariance} we obtain
\begin{eqnarray*}
& = & [\pi \times \mu] \big(\sum_{[\gamma] \in E^{\mathbf{s}(x)}_{g^{-1},e}} 1_{\mathbf{r}(x)\gamma g\gm} * \gm g \gm * [a]_{x\gm} \big)\xi\\
& = & \sum_{[\gamma] \in E^{\mathbf{s}(x)}_{g^{-1},e}} \widetilde{\pi}(1_{\mathbf{r}(x)\gamma g\gm}) \mu(\gm g \gm)  \widetilde{\pi} ([a]_{x\gm})\xi\\
& = & \mu(\gm g \gm)\,  \pi ([a]_{x\gm})\xi
\end{eqnarray*}

Hence, we conclude that $\omega_{\pi \times \mu} = \mu$.

Now let $\Phi$ be a $^*$-representation of $C_c(\a / \gm) \times_{\alpha}^{alg} G / \gm$ and $(\Phi|, \omega_{\Phi})$ its corresponding covariant pre-$^*$-representation. We want to prove that
\begin{align*}
 \Phi| \times \omega_{\Phi} = \Phi\,.
\end{align*}
Let $ 1_{\mathbf{r}(x) \gm} * \gm g \gm * [\alpha_{g^{-1}}(a)]_{x g\gm}$ be one of the spanning elements of $C_c(\a / \gm) \times_{\alpha}^{alg} G / \gm$ and $\xi \in \mathscr{H}$. We have
\begin{eqnarray*}
  [\Phi| \times \omega_{\Phi}]\,(1_{\mathbf{r}(x) \gm} * \gm g \gm * [\alpha_{g^{-1}}(a)]_{x g\gm})\, \xi \!\!\!\! & = &\!\!\!\! \widetilde{\Phi|}(1_{\mathbf{r}(x) \gm}) \omega_{\Phi} (\gm g \gm) \widetilde{\Phi|}([\alpha_{g^{-1}}(a)]_{x g\gm})\, \xi\,,
\end{eqnarray*}
which by Lemma \ref{two ways of getting reps of CcU} gives that
\begin{eqnarray*}
& = & \widetilde{\Phi}(1_{\mathbf{r}(x) \gm}) \widetilde{\Phi}(\gm g \gm) \widetilde{\Phi} ([\alpha_{g^{-1}}(a)]_{x g\gm})\, \xi\\
& = & \Phi(1_{\mathbf{r}(x) \gm} * \gm g \gm * [\alpha_{g^{-1}}(a)]_{x g\gm})\, \xi\,.
\end{eqnarray*}
Thus, $\Phi| \times \omega_{\Phi} = \Phi$. \qed\\

\subsection{More on covariant pre-$^*$-representations}

In the previous section we introduced the notion of covariant pre-$^*$-representations of $C_c(\a / \gm) \times_{\alpha}^{alg} G / \gm $ (Definition \ref{covariant pre rep def}) and a particular instance of these which we called covariant $^*$-representations (Definition \ref{covariant rep def}).

In this section we will see that the class of covariant pre-$^*$-representations is in general strictly larger than the class of covariant $^*$-representations. It is thus unavoidable, in general, to consider pre-representations of the Hecke algebra in the representation theory of crossed products by Hecke pairs. We shall also see, nevertheless, that in many interesting situations every covariant pre-$^*$-representation is actually a covariant $^*$-representation.\\

\begin{ex}
\label{example of int form normed but mu not normed}
 Let $(G, \gm)$ be a Hecke pair such that its corresponding Hecke algebra $\h(G, \gm)$ does not have an enveloping $C^*$-algebra (it is well known that such pairs exist, as for example $(G, \gm) = (SL_2(\mathbb{Q}_p) , SL_2(\mathbb{Z}_p))$ as discussed in \cite{hall}). The fact that the Hecke algebra does not have an enveloping $C^*$-algebra implies that there is a sequence of $^*$-representations $\{\mu_n\}_{n \in \mathbb{N}}$ of $\h(G, \gm)$ on Hilbert spaces $\{\mathscr{H}_n\}_{n \in \mathbb{N}}$ and an element $f \in \h(G, \gm)$  such that $\|\mu_n(f) \| \to \infty$. Let $\mathscr{V}$ be the inner product space $\mathscr{V}:= \bigoplus_{n \in \mathbb{N}} \mathscr{H}_n$ and $\mu: \h(G, \gm) \to L(\mathscr{V})$ the diagonal pre-$^*$-representation
\begin{align*}
 \mu := \bigoplus_{n \in \mathbb{N}} \mu_n\,,
\end{align*}
which of course is not normed. Let $X = \{x_1, x_2, \dots \}$ be an infinite countable set, with the trivial groupoid structure, i.e. $X$ is just a set. We consider the Fell bundle $\a$ over $X$ whose fibers are the complex numbers, i.e. $\a_x = \mathbb{C}$ for every $x \in X$, and we consider the trivial action of $G$ on $\a$, i.e. the action that fixes every element of $\a$. Thus, the action is $\gm$-good and has the $\gm$-intersection property. We also have that
\begin{align*}
 C_c(\a / \gm) = C_c(X) = C_c(X^0 / \gm)\,.
\end{align*}

Let $\pi: C_c(X) \to B(\overline{\mathscr{V}})$ be the $^*$-representation on the Hilbert space completion $\overline{\mathscr{V}}$ of $\mathscr{V}$ such that $\pi(1_{x_n})$ is the projection onto the subspace $\mathscr{H}_n$.

We claim that $(\pi, \mu)$ is a covariant pre-$^*$-representation of $C_c(X) \times^{alg}_{\alpha} G/ \gm$. To see this, first we notice that $\pi$ is obviously nondegenerate and moreover $\pi(C_c(X)) \overline{\mathscr{V}} = \mathscr{V}$, which is the inner product space where $\mu$ is defined. Next we notice that for every $x_n \in X$ and $g \in G$, the operators $\pi(1_{x_n})$ and  $\mu(\gm g \gm)$ commute. Moreover, we have
\begin{align*}
 \pi(1_{x_n}) \mu(\gm g \gm) \pi(1_{x_n}) = \mu_n(\gm g \gm)\,,
\end{align*}
on the subspace $\mathscr{H}_n$. Also we have
\begin{eqnarray*}
& & \mu(\gm g \gm) \pi(1_{x_n}) \mu(\gm s \gm) \;\; = \;\;\\
 & = & \mu(\gm g \gm) \mu(\gm s \gm) \pi(1_{x_n}) \\
& = & \sum_{\substack{[u] \in \gm g^{-1} \gm / \gm \\ [v] \in \gm s \gm / \gm}} \frac{\Delta(g)}{L(u^{-1}v)} \;  \, \mu(\gm u^{-1} v \gm) \, \pi(1_{x_n})\\
& = & \sum_{\substack{[u] \in \gm g^{-1} \gm / \gm \\ [v] \in \gm s \gm / \gm}} \frac{\Delta(g)}{L(u^{-1}v)} \;  \, \pi(1_{x_n})\mu(\gm u^{-1} v \gm) \, \pi(1_{x_n})\\
& = & \sum_{\substack{[u] \in \gm g^{-1} \gm / \gm \\ [v] \in \gm s \gm / \gm}} \sum_{[\gamma] \in E^{x_n}_{u,v}} \frac{\Delta(g)N^{x_n \gamma}_{u^{-1}, v}}{L(u^{-1}v)} \;  \, \pi(1_{x_n\gamma u})\mu(\gm u^{-1} v \gm) \, \pi(1_{x_n\gamma v})\,,
\end{eqnarray*}
where the last equality comes from the fact that since $\mathcal{S}_{x_n} = G$ we must have that $E^{x_n}_{u,v}$ consists only of the class $[e]$,  $N_{u^{-1}, v}^{x_n} = 1$ and also that $1_{x_n \gamma u} = 1_{x_n} = 1_{x_n \gamma v}$.

 So we have established that $(\pi, \mu)$ is indeed a covariant pre-$^*$-representation. Nevertheless, $\mu$ is not normed, so that $(\pi, \mu)$ is not a covariant $^*$-representation.

It is worth noting that here we are in the conditions of Example \ref{action fixes all points tensor product example}, so that $C_c(X) \times^{alg}_{\alpha} G / \gm \;\; \cong \;\; C_c(X) \odot \h(G, \gm)$.\\

\end{ex}

Example \ref{example of int form normed but mu not normed} shows that there can be more covariant pre-$^*$-representations than covariant $^*$-representations. Nevertheless, the two classes actually coincide in many cases. One such case is when $C_c(\a / \gm)$ has an identity element:\\

\begin{prop}
 If the crossed product $C_c(\a / \gm) \times^{alg}_{\alpha} G /\gm$ has an identity element (equivalently, if $C_c(\a / \gm)$ has an identity element), then every covariant pre-$^*$-representation is a covariant $^*$-representation.\\
\end{prop}

{\bf \emph{Proof:}} Let us assume that $C_c(\a / \gm) \times^{alg}_{\alpha} G /\gm$ has an identity element (equivalently, $C_c(\a / \gm)$ has an identity element).

Let $(\pi, \mu)$ be a covariant pre-$^*$-representation. As it was shown in Theorem \ref{bijective correspondence cov rep - rep}, the integrated form $\pi \times \mu$ is a $^*$-representation of $C_c(\a / \gm) \times^{alg}_{\alpha} G /\gm$ such that $\mu = \omega_{\pi \times \mu}$, where $\omega_{\pi \times \mu}$ is the pre-$^*$-representation which is obtained by extending $\pi \times \mu$ to the multiplier algebra $M(C_c(\a / \gm) \times^{alg}_{\alpha} G /\gm)$ and then restricting it to $\h(G, \gm)$. Since the crossed product $C_c(\a / \gm) \times^{alg}_{\alpha} G /\gm$ has an identity element, we have
\begin{align*}
 M(C_c(\a / \gm) \times^{alg}_{\alpha} G /\gm) = C_c(\a / \gm) \times^{alg}_{\alpha} G /\gm\,,
\end{align*}
and therefore $\omega_{\pi \times \mu}$ is just the restriction of $\pi \times \mu$ to the the Hecke algebra $\h(G, \gm)$. Hence, $\mu = \omega_{\pi \times \mu}$ is a true $^*$-representation. \qed\\

Another interesting situation where covariant pre-$^*$-representations coincide with covariant $^*$-representations is when $\h(G, \gm)$ is a $BG^*$-algebra. This is known to be the case for many classes of Hecke pairs $(G, \gm)$ as we proved in \cite{palma}. Actually, most of the classes of Hecke pairs for which a full Hecke $C^*$-algebra is known to exist are such that $\h(G, \gm)$ is $BG^*$-algebra.\\

\begin{prop}
 If $\h(G, \gm)$ is a $BG^*$-algebra, then every covariant pre-$^*$-representation is a covariant $^*$-representation.\\
\end{prop}

{\bf \emph{Proof:}} If $\h(G, \gm)$ is a $BG^*$-algebra, then every pre-$^*$-representation of $\h(G, \gm)$ is automatically normed and hence arises from a true $^*$-representation. \qed\\

\subsection{Crossed product in the case of free actions}

In this section we will see that when the associated $G$-action on $X$ is free the expressions for the products of the form $\gm g \gm * [a]_{x\gm} * \gm s \gm$, described in Corollary \ref{cor prod in corossed prod}, as well as the definition of a covariant pre-$^*$-representation become much simpler and even more similar to the notion of \emph{covariant pairs} of \cite{cov}.\\

\begin{thm}
\label{covariance in case of free actions}
 If the action of $G$ on $X$ is free, then
\begin{align}
 \gm g \gm * 1_{y \gm} * \gm s \gm = \sum_{\substack{[u] \in \gm g^{-1} \gm / \gm \\ [v] \in \gm s \gm / \gm}} 1_{yu\gm} * \gm u^{-1} v \gm * 1_{yv\gm}
\end{align}
and similarly,
\begin{align}
 \gm g \gm * [a]_{x \gm} * \gm s \gm = \sum_{\substack{[u] \in \gm g^{-1} \gm / \gm \\ [v] \in \gm s \gm / \gm}} [\alpha_{u^{-1}}(a)]_{xu\gm} * \gm u^{-1} v \gm * 1_{\mathbf{s}(x)v\gm}\,.
\end{align}\\
\end{thm}

\begin{lemma}
\label{numbers n and d for free actions}
 If the action of $G$ on $X$ is free, then
\begin{align*}
 n_{w, v}^{y}=1 \qquad\text{and}\qquad d_{w, v}^{y } = [\gm^{wv} : \gm^{wv} \cap w \gm w^{-1}]\,.\\
\end{align*}
\end{lemma}

{\bf \emph{Proof:}} We have
\begin{eqnarray*}
 \mathfrak{n}_{w, v}^{y} & = &  \big\{[r] \in \gm w \gm / \gm: r^{-1}w v\gm \subseteq \gm v\gm\; \text{and} \;\; y   w^{-1} \in y \gm r^{-1}\big\}\\
& = &  \big\{[r] \in \gm w \gm / \gm: r^{-1}w v\gm \subseteq \gm v\gm\; \text{and} \;\;  w^{-1} \in  \gm r^{-1}\big\}\\
& = & \big\{[r] \in \gm w \gm / \gm: r^{-1}w v\gm \subseteq \gm v\gm\; \text{and} \;\;    r\gm =  w \gm \big\}\\
& = & \{w\gm\}\,.
\end{eqnarray*}
Thus, $n_{w, v}^{y}=1$. Also,
\begin{eqnarray*}
 \mathfrak{d}_{w, v}^{y} & = &  \big\{[r] \in \gm w \gm / \gm: r^{-1}w v\gm \subseteq \gm v\gm\; \text{and} \;\; y w^{-1} \in y \gm r^{-1} \gm^{wv} \big\}\\
& = &  \big\{[r] \in \gm w \gm / \gm: r^{-1}w v\gm \subseteq \gm v\gm\; \text{and} \;\;  w^{-1} \in  \gm r^{-1}\gm^{wv} \big\}\,.
\end{eqnarray*}
Now we notice that in the above set the condition $r^{-1}w v\gm \subseteq \gm v\gm$ is automatically satisfied from the second condition $w^{-1} \in  \gm r^{-1}\gm^{wv}$, because the latter means that $r^{-1} = \theta_1 w^{-1} \theta_2$ for some $\theta_1 \in \gm$ and $\theta_2 \in \gm^{wv}$. Thus, we get
\begin{eqnarray*}
 \mathfrak{d}_{w, v}^{y} & = & \big\{[r] \in \gm w \gm / \gm:  w^{-1} \in  \gm r^{-1}\gm^{wv} \big\}\\
& = &  \big\{[r] \in \gm w \gm / \gm:  r \in  \gm^{wv} w\gm \big\}\\
& = & \gm^{wv} w \gm /\gm\,.
\end{eqnarray*}
Thus, we obtain $d_{w, v}^{y} =  \big|\gm^{wv} w \gm /\gm \big| = [\gm^{wv} : \gm^{wv} \cap w \gm w^{-1}]$. \qed\\

{\bf \emph{Proof of Theorem \ref{covariance in case of free actions}:}} We have seen in Theorem \ref{product formula} that
\begin{align*}
 \gm g \gm * 1_{y \gm} * \gm s \gm = \sum_{\substack{[u] \in \gm g^{-1} \gm / \gm \\ [v] \in \gm s \gm / \gm}}\; \sum_{[\gamma] \in E_{u,v}^y} \frac{\Delta(g)N_{u^{-1}, v}^{y\gamma}}{L(u^{-1}v)} \; \big(1_{y \gamma u \gm} * \gm u^{-1} v \gm * 1_{y \gamma v \gm}\big)
\end{align*}

It follows from Lemma \ref{numbers n and d for free actions} that
\begin{align*}
 N_{u^{-1}, v}^{y\gamma} = \frac{1}{[\gm^{u^{-1}v} : \gm^{u^{-1}v} \cap u^{-1} \gm u]}\,.
\end{align*}
Moreover, freeness of the action also implies that
\begin{eqnarray*}
 E_{u, v}^y & = & S_y \backslash \gm / (v \gm v^{-1} \cap u \gm u^{-1})\\
 & = & \gm / (v \gm v^{-1} \cap u \gm u^{-1})\,.
\end{eqnarray*}
Now, we have the following well-defined bijective correspondence
\begin{align*}
\gm / (\gm^u \cap \gm^v)\; &\longrightarrow\; \gm / (v \gm v^{-1} \cap u \gm u^{-1})\\
[\theta]\; &\mapsto\; [\theta]\,,
\end{align*}
given by Proposition \ref{double coset spaces prop}. Note that $\gm^u \cap \gm^v$ is simply the subgroup $u\gm u^{-1} \cap v \gm v^{-1} \cap \gm$, but in the following we will take preference on the notation $\gm^u \cap \gm^v$ for being shorter.

Consider now the action of $\gm$ on $G / \gm \times G / \gm$ by left multiplication and denote by $\mathcal{O}_{h_1, h_2}$ the orbit of the element $(h_1 \gm, h_2 \gm) \in G / \gm \times G / \gm$. It is easy to see that the map
\begin{align*}
\gm / (\gm^{h_1} \cap \gm^{h_2})\; &\longrightarrow\; \mathcal{O}_{h_1, h_2}\\
 [\theta] \; &\mapsto \; (\theta h_1 \gm, \theta h_2 \gm)
\end{align*}
is also  well-defined and is a bijection. We will denote by $\mathcal{C}$ the set of all orbits contained in $\gm g^{-1} \gm / \gm \times \gm s \gm / \gm$ (note that this set is $\gm$-invariant, so that it is a union of orbits). We then have
\begin{eqnarray*}
&  & \gm g \gm * 1_{y \gm} * \gm s \gm  =\\
& = & \sum_{\substack{[u] \in \gm g^{-1} \gm / \gm \\ [v] \in \gm s \gm / \gm}}\; \sum_{[\gamma] \in E_{u,v}^y} \frac{\Delta(g)N_{u^{-1}, v}^{y\gamma}}{L(u^{-1}v)} \; \big(1_{y \gamma u \gm} * \gm u^{-1} v \gm * 1_{y \gamma v \gm}\big)\\
& = & \sum_{\substack{[u] \in \gm g^{-1} \gm / \gm \\ [v] \in \gm s \gm / \gm}}\; \sum_{[\gamma] \in \gm / (\gm^u \cap \gm^v)} \frac{\Delta(g)N_{u^{-1}, v}^{y\gamma}}{L(u^{-1}v)} \; \big(1_{y \gamma u \gm} * \gm u^{-1} v \gm * 1_{y \gamma v \gm}\big)\\
& = & \sum_{\substack{[u] \in \gm g^{-1} \gm / \gm \\ [v] \in \gm s \gm / \gm}}\; \sum_{[\gamma] \in \gm / (\gm^u \cap \gm^v)} \frac{\Delta(g)N_{u^{-1} \gamma^{-1}, \gamma v}^y}{L(u^{-1}\gamma^{-1} \gamma v)} \; \big(1_{y \gamma u \gm} * \gm u^{-1} \gamma^{-1} \gamma v \gm * 1_{y \gamma v \gm}\big)
\end{eqnarray*}
where the last equality comes from the fact that $N^{y \gamma}_{u^{-1}, v} = N^y_{u^{-1} \gamma^{-1}, \gamma v}$, which is a consequence of Lemma \ref{lemma numbers n and d} $iii)$, or simply by Lemma \ref{numbers n and d for free actions}. Using now the bijection between $\gm / (\gm^u \cap \gm^v)$ and the orbit space $\mathcal{O}_{u,v}$ as described above, we obtain
\begin{eqnarray*}
 & = & \sum_{\substack{[u] \in \gm g^{-1} \gm / \gm \\ [v] \in \gm s \gm / \gm}}\; \sum_{([r], [t]) \in \mathcal{O}_{u,v}} \frac{\Delta(g)N_{r^{-1}, t}^{y}}{L(r^{-1}t)} \; \big(1_{y r \gm} * \gm r^{-1} t \gm * 1_{y t \gm}\big)\\
& = & \sum_{\mathcal{O} \in \mathcal{C}}\; \sum_{([u],[v]) \in \mathcal{O}}\; \sum_{([r], [t]) \in \mathcal{O}_{u,v}} \frac{\Delta(g)N_{r^{-1}, t}^{y}}{L(r^{-1}t)} \; \big(1_{y r \gm} * \gm r^{-1} t \gm * 1_{y t \gm}\big)\\
& = & \sum_{\mathcal{O} \in \mathcal{C}}\; \sum_{([u],[v]) \in \mathcal{O}}\; \sum_{([r], [t]) \in \mathcal{O}} \frac{\Delta(g)N_{r^{-1}, t}^{y}}{L(r^{-1}t)} \; \big(1_{y r \gm} * \gm r^{-1} t \gm * 1_{y t \gm}\big)\\
& = & \sum_{\mathcal{O} \in \mathcal{C}}\; \sum_{([r],[t]) \in \mathcal{O}}\; \frac{\#\mathcal{O}\,\Delta(g)N_{r^{-1}, t}^{y}}{L(r^{-1}t)} \; \big(1_{y r \gm} * \gm r^{-1} t \gm * 1_{y t \gm}\big)\,,
\end{eqnarray*}
where $\# \mathcal{O}$ denotes the total number of elements of the given orbit $\mathcal{O}$. Changing the names of the variables ($r$ to $u$ and $t$ to $v$) we get
\begin{eqnarray*}
& = & \sum_{\mathcal{O} \in \mathcal{C}}\; \sum_{([u],[v]) \in \mathcal{O}}\; \frac{\#\mathcal{O}\,\Delta(g)N_{u^{-1}, v}^{y}}{L(u^{-1}v)} \; \big(1_{y u \gm} * \gm u^{-1} v \gm * 1_{y v \gm}\big)\\
& = & \sum_{\substack{[u] \in \gm g^{-1} \gm / \gm \\ [v] \in \gm s \gm / \gm}}\; \frac{\#\mathcal{O}_{u,v}\,\Delta(g)N_{u^{-1}, v}^{y}}{L(u^{-1}v)} \; \big(1_{y u \gm} * \gm u^{-1} v \gm * 1_{y v \gm}\big)\,.
\end{eqnarray*}
We are now going to prove that the coefficients satisfy
\begin{align*}
 \frac{\#\mathcal{O}_{u,v}\,\Delta(g)N_{u^{-1}, v}^{y}}{L(u^{-1}v)} = 1\,.
\end{align*}
This follows from the following computation:
\begin{eqnarray*}
 \frac{\#\mathcal{O}_{u,v}}{L(u^{-1}v)} \,N_{u^{-1}, v}^{y} \, \Delta(g) & = & \frac{[\gm : \gm^u \cap \gm^v]}{[\gm : \gm^{u^{-1}v}]} \cdot \frac{1}{[\gm^{u^{-1}v} : \gm^{u^{-1}v} \cap u^{-1} \gm u]} \cdot \frac{[\gm : \gm^{u^{-1}}]}{[\gm : \gm^u]}\\
& = & \frac{[\gm : \gm^u \cap \gm^v]\, [\gm : \gm^{u^{-1}}]}{[\gm  : \gm^{u^{-1}v} \cap u^{-1} \gm u][\gm : \gm^u]}\\
& = & \frac{[\gm^u : \gm^u \cap \gm^v]\, [\gm : \gm^{u^{-1}}]}{[\gm  : \gm^{u^{-1}v} \cap u^{-1} \gm u]}\\
& = & \frac{[\gm^u : \gm^u \cap \gm^v]\, [u\gm u^{-1} : \gm^u]}{[\gm : \gm^{u^{-1}v} \cap u^{-1} \gm u]}\\
& = & \frac{[u\gm u^{-1} : \gm^u \cap \gm^v]}{[\gm  : \gm^{u^{-1}v} \cap u^{-1} \gm u]}\\
& = & \frac{[u\gm u^{-1} : \gm^u \cap \gm^v]}{[u\gm u^{-1}  : \gm^u \cap \gm^v]}\\
& = & 1\,.
\end{eqnarray*}

This finishes the first claim of the theorem. The second claim, concerning the product $\gm g \gm * [a]_{x \gm} * \gm s \gm$, is proven in a completely similar fashion. \qed\\

\begin{prop}
\label{covariant pre-representation free actions}
Let $\pi: C_c(\a / \gm) \to B(\mathscr{H})$ be a nondegenerate $^*$-representation, $\mu: \h(G, \gm) \to L(\pi(C_c(\a / \gm) \mathscr{H})$ a unital pre-$^*$-representation, and let us assume that the associated $G$-action on $X$ is free. The pair $(\pi, \mu)$ is a covariant pre-$^*$-representation if and only if the following equality
\begin{align}
 \mu(\gm g \gm)\pi( [a]_{x \gm} ) \mu( \gm s \gm) = \sum_{\substack{[u] \in \gm g^{-1} \gm / \gm \\ [v] \in \gm s \gm / \gm}} \pi([\alpha_{u^{-1}}(a)]_{xu\gm} ) \mu(\gm u^{-1} v \gm) \widetilde{\pi} (1_{\mathbf{s}(x)v\gm})\,.
\end{align}
holds for all $g, s \in G$, $x \in X$ and $a \in \a_x$.\\
\end{prop}

{\bf \emph{Proof:}} $(\Longrightarrow)$ Assume that $(\pi, \mu)$ is a covariant pre-$^*$-representation. Then we have
\begin{eqnarray*}
\mu(\gm g \gm)\pi( [a]_{x \gm} ) \mu( \gm s \gm) & = & [\pi \times \mu] (\gm g \gm * [a]_{x \gm} * \gm s \gm)\\
& = & [\pi \times \mu] \Big(\sum_{\substack{[u] \in \gm g^{-1} \gm / \gm \\ [v] \in \gm s \gm / \gm}} [\alpha_{u^{-1}}(a)]_{xu\gm} * \gm u^{-1} v \gm * 1_{\mathbf{s}(x)v\gm} \Big)\\
& = & \sum_{\substack{[u] \in \gm g^{-1} \gm / \gm \\ [v] \in \gm s \gm / \gm}} \pi([\alpha_{u^{-1}}(a)]_{xu\gm} ) \mu(\gm u^{-1} v \gm) \widetilde{\pi} (1_{\mathbf{s}(x)v\gm})\,.
\end{eqnarray*}

$(\Longleftarrow)$ In order to prove equality (\ref{cov pre rep equality}) one just needs to show that
\begin{eqnarray*}
& & \sum_{\substack{[u] \in \gm g^{-1} \gm / \gm \\ [v] \in \gm s \gm / \gm}}\; \sum_{[\gamma] \in E_{u,v}^{\mathbf{s}(x)}} \frac{\Delta(g)N_{u^{-1}, v}^{\mathbf{s}(x)\gamma}}{L(u^{-1}v)} \; \widetilde{\pi}([\alpha_{u^{-1}}(a)]_{x \gamma u \gm}) \, \mu(\gm u^{-1} v \gm) \, \widetilde{\pi}(1_{\mathbf{s}(x) \gamma v \gm})\\
& = & \sum_{\substack{[u] \in \gm g^{-1} \gm / \gm \\ [v] \in \gm s \gm / \gm}} \widetilde{\pi}([\alpha_{u^{-1}}(a)]_{xu\gm} ) \mu(\gm u^{-1} v \gm) \widetilde{\pi} (1_{\mathbf{s}(x)v\gm})\,,
\end{eqnarray*}
and  this is proven in a completely analogous way as in the proof of Theorem \ref{covariance in case of free actions}. \qed\\

\end{document}